\documentclass[10pt]{amsart}

\usepackage[text={420pt,660pt},centering]{geometry}
\usepackage{graphicx}
\usepackage{amssymb}
\usepackage{mathtools}
\usepackage[colorlinks=true, pdfstartview=Fit, linkcolor=blue, citecolor=blue, urlcolor=blue,pagebackref=false]{hyperref}
\hypersetup{
    pdftitle={On Hamilton-Jacobi equations and perceptron models},
    pdfkeywords={spin glasses, bipartite models, perceptrons, Hamilton--Jacobi equations}
}
\usepackage{orcidlink}
\usepackage[expansion=false]{microtype}

\usepackage{bm}

\usepackage{mathrsfs}
\usepackage{xcolor}

\newtheorem{proposition}{Proposition}
\newtheorem{theorem}[proposition]{Theorem}
\newtheorem{lemma}[proposition]{Lemma}

\theoremstyle{remark}
\newtheorem{remark}[proposition]{Remark}

\theoremstyle{definition}
\newtheorem{definition}[proposition]{Definition}

\numberwithin{equation}{section}
\numberwithin{proposition}{section}
\numberwithin{figure}{section}
\numberwithin{table}{section}

\newcommand{\N}{\mathbb{N}}

\newcommand{\R}{\mathbb{R}}

\newcommand{\E}{\mathbb{E}}

\newcommand{\eps}{\varepsilon}

\renewcommand{\le}{\leqslant}
\renewcommand{\ge}{\geqslant}
\renewcommand{\leq}{\leqslant}
\renewcommand{\geq}{\geqslant}

\renewcommand{\bar}{\overline}
\newcommand{\Ll}{\left}
\newcommand{\Rr}{\right}
\renewcommand{\d}{\mathrm{d}}
\newcommand{\dr}{\partial}
\newcommand{\mcl}{\mathcal}

\DeclareMathOperator{\supp}{supp}

\newcommand{\sP}{\mathscr{P}}

\newcommand{\la}{\left\langle}
\newcommand{\ra}{\right\rangle}

\renewcommand{\H}{\mathsf{H}}
\newcommand{\cH}{\mathcal H}

\newcommand{\one}{\mathbf{1}}
\newcommand{\bip}{\mathrm{bip}}

\newcommand{\fR}{\mathfrak{R}}

\renewcommand{\Vec}{{\mathrm{vec}}}

\newcommand{\ext}{{\mathrm{ext}}}
\newcommand{\sph}{{\mathrm{sph}}}

\newcommand{\Sph}{\mathbb{S}}

\newcommand{\nabsub}[1]{\nabla_{\!#1}}
\newcommand{\nablatx}{\nabsub{t,x}}

\begin{document}

\author[H.-B.\ Chen]{Hong-Bin Chen\,\orcidlink{0000-0001-6412-0800}}
\address[Hong-Bin Chen]{NYU-ECNU Institute of Mathematical Sciences, NYU Shanghai, China}

\author[J.-C.\ Mourrat]{Jean-Christophe Mourrat\,\orcidlink{0000-0002-2980-725X}}
\address[Jean-Christophe Mourrat]{Department of Mathematics, ENS Lyon and CNRS, Lyon, France}

\title{On Hamilton--Jacobi equations and perceptron models}

\begin{abstract}
The solution to a finite-dimensional Hamilton--Jacobi equation is classically known to admit a variational representation in two situations: when the nonlinearity is convex or concave, by the Hopf--Lax formula; and when the initial condition is convex or concave, by the Hopf formula. We identify a third situation: the space variable is split into two groups of coordinates, the nonlinearity is the scalar product of the two corresponding partial gradients, and the initial condition is the sum of a convex function of the first group of coordinates and of an arbitrary Lipschitz function of the second group. This structure is precisely that of the Hamilton--Jacobi equations describing the limit free energy of perceptron models. We illustrate this point by identifying the limit free energy of a class of perceptron models that can be rewritten as bipartite spin glasses, with Ising spins in one layer and bounded soft spins in the other. We also discuss how to extend the approach to several other classes of models, including general Ising and spherical perceptrons, possibly with labelled patterns.
\end{abstract}

\maketitle

\section{Introduction}
\label{s.intro}

\subsection{A variational formula for Hamilton--Jacobi equations}
\label{s.intro.hj}

Let $D \ge 1$ be an integer. We consider Hamilton--Jacobi equations of the form
\begin{equation}
\label{e.hj.general}
\dr_t u - \H(\nabla u) = 0 \quad \text{on } (0,\infty) \times \R^D, \qquad u(0,\cdot) = \psi \quad \text{on } \R^D,
\end{equation}
where the nonlinearity $\H : \R^D \to \R$ and the initial condition $\psi : \R^D \to \R$ are given, and where the unknown is $u : \R_+ \times \R^D \to \R$, with $\R_+ := [0,\infty)$. Throughout the paper, we write $\nabla$ for the gradient with respect to the space variable, and $\dr_t$ for the derivative with respect to the time variable. Solutions to~\eqref{e.hj.general} are understood in the viscosity sense; see Definition~\ref{d.vs_general} below. When $\H$ is locally Lipschitz continuous and $\psi$ is Lipschitz continuous, the equation~\eqref{e.hj.general} admits a unique Lipschitz continuous viscosity solution; see for instance~\cite[Corollary~3.7]{HJbook}.

There are two classical situations in which this solution admits a variational representation. The first one is when the nonlinearity $\H$ is convex. In this case, the Hopf--Lax formula states that, for every $t > 0$ and $x \in \R^D$,
\begin{equation}
\label{e.hopf.lax}
u(t,x) = \sup_{y \in \R^D} \inf_{p \in \R^D} \{\psi(y) + p \cdot (x-y) + t \H(p) \} =  \sup_{y \in \R^D} \Ll\{ \psi(y) - t \H^*\Ll(\frac{y-x}{t}\Rr) \Rr\} ,
\end{equation}
where $\H^*(y) := \sup_{p \in \R^D} \{p \cdot y - \H(p)\}$ denotes the convex conjugate of $\H$. The case of a concave nonlinearity reduces to this one by considering $-u$ in place of $u$, and yields a formula with $\inf_y \sup_p$ in place of $\sup_y \inf_p$. The second classical situation is when the initial condition $\psi$ is convex. In this case, the Hopf formula states that, for every $t \in \R_+$ and $x \in \R^D$,
\begin{equation}
\label{e.hopf.intro}
u(t,x) = \sup_{p \in \R^D} \inf_{y \in \R^D} \Ll\{ \psi(y) + p \cdot (x-y) + t \H(p) \Rr\} = \sup_{p \in \R^D} \Ll\{ p \cdot x + t \H(p) - \psi^*(p) \Rr\},
\end{equation}
and a symmetric formula, with $\inf_p \sup_y$ in place of $\sup_p \inf_y$, holds when $\psi$ is concave. We refer to~\cite{bardi1984hopf, HJbook, evans2010partial, lions1986hopf} for proofs of these formulas.

The mechanism behind the Hopf formula~\eqref{e.hopf.intro} can be described as follows. If the initial condition is affine, say $\psi(x) = p \cdot x - c$ for some $p \in \R^D$ and $c \in \R$, then the solution to~\eqref{e.hj.general} is simply $u(t,x) = p \cdot x + t\H(p) - c$. A convex function is the supremum of the affine functions lying below it, namely $\psi(x) = \sup_{p \in \R^D} \{p \cdot x - \psi^*(p)\}$, and the Hopf formula asserts that the solution started from $\psi$ is the supremum of the solutions started from these affine functions. Notice also that in order for \eqref{e.hopf.intro} to be the solution to \eqref{e.hj.general}, the convexity of $\psi$ is necessary, as one can see by evaluating \eqref{e.hopf.intro} at $t = 0$. 

In this paper, we present a new circumstance in which the solution to~\eqref{e.hj.general} admits a variational formula. We take $D = 2d$ for some integer $d \ge 1$, and we write $x = (x_1,x_2) \in \R^d \times \R^d$ for the space variable, and $\nabla = (\nabsub{x_1}, \nabsub{x_2})$ for the gradient. We consider the bilinear nonlinearity $\H(p_1,p_2) := p_1 \cdot p_2$, and initial conditions of the separable form $\psi(x_1,x_2) = \psi_1(x_1) + \psi_2(x_2)$. The equation~\eqref{e.hj.general} then becomes
\begin{equation}
\label{e.hj.bip}
\dr_t u - \nabsub{x_1} u \cdot \nabsub{x_2} u = 0 \quad \text{on } (0,\infty) \times \R^d \times \R^d, \qquad u(0,x_1,x_2) = \psi_1(x_1) + \psi_2(x_2).
\end{equation}
We assume that $\psi_1$ is convex, but we make no convexity assumption on $\psi_2$. Since the nonlinearity~$\H$ is neither convex nor concave, and since the initial condition $\psi$ is in general neither convex nor concave, neither of the two classical formulas recalled above applies to this situation. Nevertheless, we have the following result.

\begin{theorem}[Variational formula]
\label{t.hopf.half}
Let $\psi_1, \psi_2 : \R^d \to \R$ be Lipschitz continuous, with Lipschitz constants $L_1$ and $L_2$ respectively, and assume that $\psi_1$ is convex. Let $\psi_1^* : \R^d \to \R \cup \{+\infty\}$ denote the convex conjugate of $\psi_1$, that is, $\psi_1^*(p) := \sup_{y \in \R^d} \{p \cdot y - \psi_1(y)\}$. For every $(t,x_1,x_2) \in \R_+ \times \R^d \times \R^d$, let
\begin{equation}
\label{e.def.u.half}
u(t,x_1,x_2) := \sup_{p \in \R^d} \Ll\{ p \cdot x_1 + \psi_2(x_2 + tp) - \psi_1^*(p) \Rr\}.
\end{equation}
The supremum in~\eqref{e.def.u.half} is attained, and it can be restricted to the closed ball of radius $L_1$ centered at the origin. The function $u$ is the unique Lipschitz continuous viscosity solution to~\eqref{e.hj.bip}. Moreover, for every $t,s \in \R_+$ and $x_1,x_2,y_1,y_2 \in \R^d$, we have
\begin{equation}
\label{e.lip.u.half}
|u(t,x_1,x_2) - u(s,y_1,y_2)| \le L_1 |x_1 - y_1| + L_2 |x_2 - y_2| + L_1 L_2 |t-s|.
\end{equation}
\end{theorem}

The heuristic explanation for Theorem~\ref{t.hopf.half} is similar to that for the Hopf formula, but it only uses convexity in the first group of coordinates. Indeed, suppose first that $\psi_1$ is affine, say $\psi_1(x_1) = p \cdot x_1 - c$ for some $p \in \R^d$ and $c \in \R$, and let $\psi_2$ be an arbitrary Lipschitz continuous function. Then the solution to~\eqref{e.hj.bip} is explicit:
\begin{equation}
\label{e.sol.affine}
u(t,x_1,x_2) = p \cdot x_1 + \psi_2(x_2 + tp) - c.
\end{equation}
Indeed, the gradient of this function with respect to $x_1$ is constant equal to $p$, so that for this function, the equation in~\eqref{e.hj.bip} reduces to the transport equation $\dr_t u - p \cdot \nabsub{x_2} u = 0$, which is solved by translation. Writing $\psi_1(x_1) = \sup_{p \in \R^d} \{p \cdot x_1 - \psi_1^*(p)\}$, we see that the formula~\eqref{e.def.u.half} asserts that the solution started from $\psi_1(x_1) + \psi_2(x_2)$ is the supremum over $p \in \R^d$ of the solutions in~\eqref{e.sol.affine} with $c = \psi_1^*(p)$. The fact that this supremum is a viscosity subsolution follows directly from this observation. The proof that it is a viscosity supersolution requires a bit more work, and also relies on the convexity of $\psi_1$. We prove Theorem~\ref{t.hopf.half} in Section~\ref{s.finite}.

\begin{remark}
\label{r.hopf.half}
The formula~\eqref{e.def.u.half} can be rewritten as
\begin{equation}
\label{e.def.u.half.supinf}
u(t,x_1,x_2) = \sup_{p \in \R^d} \inf_{y \in \R^d} \Ll\{ \psi_1(y) + \psi_2(x_2 + tp) + p \cdot (x_1 - y) \Rr\}.
\end{equation}
If $\psi_2$ is also convex, then the Hopf formula applies to~\eqref{e.hj.bip}, and gives
\begin{equation}
\label{e.maybe.Hopf}
u(t,x_1,x_2) = \sup_{p_1,p_2 \in \R^d} \Ll\{ p_1 \cdot x_1 + p_2 \cdot x_2 + t \, p_1 \cdot p_2 - \psi_1^*(p_1) - \psi_2^*(p_2) \Rr\}.
\end{equation}
Taking the supremum over $p_2$ first and using that $\psi_2^{**} = \psi_2$, we recover~\eqref{e.def.u.half}. The point of Theorem~\ref{t.hopf.half} is that~\eqref{e.def.u.half} remains valid without any convexity assumption on $\psi_2$. 
\end{remark}

\begin{remark}
\label{r.hopf.half.symmetries}
The roles of the two groups of coordinates in Theorem~\ref{t.hopf.half} can of course be exchanged. Moreover, if $u$ solves~\eqref{e.hj.bip}, then the function $(t,x_1,x_2) \mapsto -u(t,x_1,-x_2)$ solves the same equation, with initial condition $-\psi_1(x_1) - \psi_2(-x_2)$. The conclusions of Theorem~\ref{t.hopf.half} can therefore be transposed to the case in which one of the two functions $\psi_1$ and $\psi_2$ is concave, and the other one is arbitrary. For instance, when $\psi_1$ is concave, the solution is given by~\eqref{e.def.u.half.supinf} with $\inf_p \sup_y$ in place of $\sup_p \inf_y$.
\end{remark}

For the application to perceptron models, we need a version of Theorem~\ref{t.hopf.half} in which the space~$\R^d$ is replaced by a closed convex cone, with no boundary condition imposed, in the spirit of~\cite{chen2022cones}, and ultimately by an infinite-dimensional cone of increasing paths. These extensions are the objects of Sections~\ref{s.finite} and~\ref{s.paths} respectively.

\subsection{Perceptron models and bipartite spin glasses}
\label{s.intro.perceptron}

While the structural assumptions of Theorem~\ref{t.hopf.half} may seem rather narrow, they are in fact satisfied by the Hamilton--Jacobi equations that are relevant to the study of the free energy of perceptron models, as we now explain. Let $N, M \ge 1$ be integers, let $(g_{ij})_{i \le N, j \le M}$ be independent standard Gaussian random variables, let $P_N$ be a probability measure on $\R^N$, and let $U : \R \to \R$ be a measurable function, which we call the activation function. The free energy of the associated perceptron model is
\begin{equation}
\label{e.def.perc.general}
-\frac{1}{N} \E \log \int \exp \Ll( \sum_{j=1}^{M} U\Ll( \frac{1}{\sqrt N} \sum_{i=1}^N g_{ij} \sigma_i \Rr) \Rr) \, \d P_N(\sigma),
\end{equation}
where $\E$ denotes the expectation with respect to the Gaussian random variables $(g_{ij})$. We are interested in the limit of this quantity as $N$ tends to infinity, with $M/N$ converging to some $\lambda \in (0,\infty)$. The most studied cases are the Ising perceptron, in which $P_N$ is the uniform measure on $\{-1,1\}^N$, and the spherical perceptron, in which $P_N$ is the uniform measure on the sphere of radius $\sqrt N$. 

In the interpretation of this model that goes back to~\cite{gardner1988, gardnerderrida1988}, the vector $\sigma$ encodes the synaptic weights of a single neuron, and each vector $(g_{ij})_{i \le N}$ is a pattern presented to this neuron, multiplied by the output in $\{-1,1\}$ that the neuron is required to produce in response; by symmetry, this multiplication does not change the law of the pattern. The neuron responds with the sign of its pre-activation, so the quantity $N^{-1/2} \sum_{i} g_{ij} \sigma_i$ is positive if and only if the response to the $j$-th pattern is correct, and $U$ assigns a reward to its value. For instance, the choice $U = -\beta \one_{(-\infty,\kappa)}$ penalizes each pattern that is not stored with margin~$\kappa$, while the hard constraint $U = \log \one_{[\kappa,\infty)}$, which takes the value $-\infty$, forbids such patterns altogether; the storage capacity is the largest value of $\lambda$ for which this constraint can typically be satisfied. 

Predictions for the free energy based on the replica method were developed in~\cite{gardner1988, gardnerderrida1988, gyorgyi2001techniques, gyorgyi1997parisi, gyorgyi2000beyond, krauthmezard1989}; a cavity derivation of the replica-symmetric formulas was given in~\cite{mezard1989space}. Related predictions in the context of empirical risk minimization were developed in~\cite{montanari2026replica}. 
Rigorous formulas for the Ising perceptron were first obtained in regimes of high temperature or small~$\lambda$ in \cite{talagrand2000halfspaces, talagrand2002}. For the hard constraint $U=\log\one_{[\kappa,\infty)}$, the replica-symmetric formula at small~$\lambda$ was established in \cite[Chapter~9]{Tbook2}; a different proof, covering a broader class of activation functions, was subsequently given in~\cite{bolthausen2022gardner}. For the spherical perceptron with the same hard constraint and $\kappa\geq0$, the Gardner formula was proved in~\cite{shcherbina2003rigorous}.

Much effort has also been devoted to determining the storage capacity and to understanding the geometry of the set of solutions; see for instance~\cite{kim1998covering, nakajimasun2022, talagrand1999intersecting, xu2021sharp} for the Ising perceptron, and~\cite{abbeli2021, aubin2019storage, perkinsxu2024} for its symmetric variants, such as $U = \log \one_{[-\kappa,\kappa]}$. For the zero-margin Ising perceptron, the capacity predicted in~\cite{krauthmezard1989} was established in~\cite{dingsun2025, huang2024capacity, shmalo2026capacity}.

Variational formulas for the limit free energy, valid for every~$\lambda$ and at every temperature, were recently obtained for the Ising perceptron in~\cite{openai2026ising} and for the spherical perceptron in~\cite{aronow2026negative, openai2026spherical}. In~\cite{fuSohn2026perceptrons}, variational formulas are obtained for both Ising and spherical perceptrons with i.i.d.\ random activation functions, independent of the Gaussian patterns, under a common Lipschitz bound and an integrability assumption at the origin. We also refer to~\cite{camilli2026perceptron} for variational bounds on the free energy of perceptrons with log-concave spin priors and concave activation functions.

In this paper, we relate perceptron models and Hamilton--Jacobi equations by following the approach developed for spin glasses in~\cite{chen2024ms, chen2026ising, chen2025free, chen2026spherical, HJbook, mourrat2021nonconvex, mourrat2022parisi, mourrat2023free, mourrat2020extending}. That is, we enrich the model in~\eqref{e.def.perc.general} by adding to its Hamiltonian Gaussian fields indexed by the leaves of a Poisson--Dirichlet cascade. One of these fields acts on the spins $\sigma_i$, and another one is added to the pre-activations $N^{-1/2} \sum_{i} g_{ij} \sigma_i$ inside $U$; the covariance structure of each of these fields is encoded by an increasing path. We also introduce a time parameter $t$, which interpolates between the original pre-activations and independent Gaussian random variables. The enriched free energy is then a function of $t$ and of the pair of paths. By Gaussian integration by parts, this function satisfies, up to an error term which vanishes if the overlaps concentrate conditionally on the overlap of the cascade, a Hamilton--Jacobi equation whose nonlinearity is the $L^2$ scalar product of the two partial derivatives with respect to the paths. This is the infinite-dimensional analogue of the bilinear nonlinearity in~\eqref{e.hj.bip}. At $t = 0$, the two groups of variables decouple, so the initial condition is the sum of a term that depends only on the path associated with the spins, and of a term that depends only on the path associated with the pre-activations. This is the analogue of the separable initial condition in~\eqref{e.hj.bip}. The first term is the cascade transform of the reference measure for one spin, and it does not depend on the activation function $U$; for Ising spins, this term was shown to be convex in~\cite[Proposition~2.2]{chen2026ising} and in~\cite{ho2026concavity}; for spherical spins, this term was shown to be convex in~\cite[Lemma~2.9]{chen2026spherical}. The second term depends on $U$, and in general it is not convex. This is exactly the structure of Theorem~\ref{t.hopf.half}, and the infinite-dimensional analogue of the formula~\eqref{e.def.u.half} is therefore expected to describe the limit free energy of perceptron models in full generality, the choice of the activation function $U$ affecting only the non-convex part of the initial condition. We refer to Section~\ref{s.general_perceptron} for a more detailed discussion of general Ising and spherical perceptrons, including the derivation of the candidate Hamilton--Jacobi equation and of the resulting candidate variational formula.

We illustrate this point with full rigor on a subclass of perceptron models which can conveniently be rewritten as bipartite spin glasses. We fix a probability measure $\nu$ with compact support, and consider activation functions of the form
\begin{equation}
\label{e.u_t,nu=}
u_{t,\nu}(x) := \log \int e^{\sqrt{2t} x v - t v^2} \, \d \nu(v), \qquad t \in \R_+.
\end{equation}
The parameter $t$ in~\eqref{e.u_t,nu=} plays the role of an inverse temperature, and it will also be the time variable of the Hamilton--Jacobi equation. For such activation functions, the free energy in~\eqref{e.def.perc.general} can be rewritten as that of a model with an auxiliary layer of $\nu$-distributed spins $v = (v_j)_{j \le M}$, and with a Gaussian Hamiltonian which is bilinear in $(\sigma,v)$; see~\eqref{e.def.FN.delta0} below. Related bipartite models arise in the study of restricted Boltzmann machines; see for instance~\cite{tubiana2017emergence}. This restriction is meant to simplify the presentation: it allows us to borrow cavity arguments and comparison results already developed for multi-species spin glasses~\cite{chen2024ms, chen2026ising, chen2025free, mourrat2021nonconvex, mourrat2023free}, and to focus on the new ingredient, which is the variational formula. For general activation functions, the Hamiltonian in~\eqref{e.def.perc.general} is not a Gaussian process, and these arguments need to be redone.

Our rigorous results on perceptron models are less general than those presented in~\cite{fuSohn2026perceptrons, openai2026ising, openai2026spherical}. Our aim with this paper is to present a systematic and clear structure for the study of the free energy of perceptron and related models. While we only identify the limit free energy of a restricted class of models with full rigor, we discuss in Sections~\ref{s.moving_domain} and~\ref{s.general_perceptron} how the approach can plausibly be generalized to the cases covered in \cite{fuSohn2026perceptrons, openai2026ising, openai2026spherical} and more (the entirety of Section~\ref{s.moving_domain} and several results in Section~\ref{s.general_perceptron} are fully rigorous, and the missing ingredients are clearly identified). To the best of our understanding, the approaches proposed in \cite{fuSohn2026perceptrons, openai2026ising, openai2026spherical} do not match the one we present here: the definition of the enriched free energy and the relevant Hamilton--Jacobi equation in \cite{openai2026ising, openai2026spherical} seem to be different from ours, while \cite{fuSohn2026perceptrons} does not seem to rely on Hamilton--Jacobi equations. The approach of \cite{aronow2026negative} seems closer to ours.

\subsection{The bipartite model and the main rigorous result}
\label{s.intro.main}

We now give the precise definition of the model that we study in full rigor. We fix $\lambda\in(0,\infty)$ and let $(M_N)_{N\in\N}$ be a sequence of positive integers satisfying $\lim_{N\to\infty} \frac{M_N}{N} = \lambda$. We fix a probability measure $\nu$ supported on $[-1,1]$; see Remark~\ref{r.bounded_support} for the extension to an arbitrary compact support. Let $(g_{ij})_{i,j\in\N}$ be a collection of i.i.d.\ standard Gaussian random variables. For $t \in \R_+$, we consider the Ising perceptron model with Hamiltonian $H_N^\mathrm{perc}:\{-1,1\}^N\to\R$ given by
\begin{equation*}
    H_N^\mathrm{perc}(\sigma) := \sum_{j=1}^{M_N} u_{t,\nu}\Ll(\frac{1}{\sqrt{N}}\sum_{i=1}^{N}g_{ij}\sigma_i\Rr),
\end{equation*}
where $u_{t,\nu}$ is defined in~\eqref{e.u_t,nu=}. The associated free energy is
\begin{equation}
\label{e.def.FN.perc}
    \bar F_N^\mathrm{perc}:= -\frac{1}{N}\E \log \sum_{\sigma\in\{-1,1\}^N}2^{-N}\exp\Ll(H^\mathrm{perc}_N(\sigma)\Rr).
\end{equation}
As announced, this model can be rewritten as a bipartite model with Ising spins in one layer and $\nu$-distributed spins in the other. We define the bipartite Hamiltonian $H_N:\{-1,1\}^N\times[-1,1]^{M_N}\to\R$ by
\begin{equation}
\label{e.def.HN.bip}
    H_N(\sigma,v):= \frac{1}{\sqrt{N}}\sum_{\substack{i\in\{1,\dots,N\}\\j\in\{1,\dots,M_N\}}}g_{ij} \sigma_i v_j,
\end{equation}
where we have written $\sigma=(\sigma_1,\dots,\sigma_N)$ and $v=(v_1,\dots,v_{M_N})$. Its covariance is given by
\begin{equation*}
    \E \Ll[H_N(\sigma,v)H_N(\sigma',v')\Rr] = \frac{1}{N} (\sigma\cdot\sigma')( v\cdot v'),
\end{equation*}
and thus $\E \Ll[H_N(\sigma,v)^2\Rr]=|v|^2$. Using~\eqref{e.u_t,nu=} and writing $\mu_{\mcl I}:= \frac{1}{2}(\delta_{-1}+\delta_1)$ for the Ising spin distribution, we have
\begin{equation}\label{e.def.FN.delta0}
    \bar F_N(t,0):= \bar F^\mathrm{perc}_N=-\frac{1}{N}\E\log \int \exp\Ll(\sqrt{2t}H_N(\sigma,v)-t|v|^2\Rr)\ \d \Ll( \mu_{\mcl I}^{\otimes N}\otimes \nu^{\otimes M_N}\Rr)(\sigma,v).
\end{equation}
Notice that the correction term $-t|v|^2$ is $-\frac{1}{2}$ times the variance of $\sqrt{2t}H_N(\sigma,v)$.

\smallskip
\noindent\textbf{Enriched models.}
We let $\mcl Q$ denote the set of right-continuous increasing paths $q:[0,1)\to\R_+$ (here and throughout, we say that a path $q$ is increasing provided that $q(s)\leq q(s')$ for every $0\leq s\leq s'<1$). For $r\in [1,\infty]$, we write $\mcl Q_r := \mcl Q \cap L^r[0,1)$. For $q\in\mcl Q_\infty$, we set $q(1):=\lim_{r\uparrow1}q(r)$, which exists by monotonicity. We use $|\cdot|_{L^r}$ for the usual norm on $L^r([0,1))$. For $c\geq0$, we set
\begin{equation*}
    \mcl Q_{\infty,\leq c}:=\{q\in\mcl Q_\infty:q(r)\leq c,\ \forall r\in[0,1)\}.
\end{equation*}

Below we will extend $\bar F_N(t,\cdot)$ so that its second argument is any element of the set $\mcl Q_2^2=\mcl Q_2\times \mcl Q_2$; the $0$ appearing as the second argument of $\bar F_N(t,0)$ in \eqref{e.def.FN.delta0} is the pair of paths that are constant equal to zero.

Let $\fR$ be a continuous Poisson--Dirichlet cascade independent of $(g_{ij})$ whose overlap is uniformly distributed on $[0,1]$. This is a random probability measure on the unit sphere of a separable Hilbert space, with almost surely ultrametric support; see~\cite[Section~4]{chen2025free} for its construction. For $\alpha,\alpha'\in\supp\fR$, we write $\alpha\wedge\alpha'$ for their scalar product. Under $\E\fR^{\otimes2}$, the overlap $\alpha\wedge\alpha'$ is uniformly distributed on $[0,1]$.

For $q=(q_1,q_2)\in\mcl Q_\infty^2$, conditionally on $\fR$, let $W_N^q$ be a centered Gaussian process indexed by $(\sigma,v,\alpha)\in\{-1,1\}^N\times[-1,1]^{M_N}\times\supp\fR$, independent of $H_N$, with covariance
\begin{equation}\label{e.E[WW]=}
    \E\Ll[W_N^q(\sigma,v,\alpha)W_N^q(\sigma',v',\alpha')\mid\fR\Rr]
    =q_1(\alpha\wedge\alpha')\sigma\cdot\sigma'+q_2(\alpha\wedge\alpha')v\cdot v'.
\end{equation}
The field $W_N^q$ is the sum of independent scalar cascade fields, one for each spin coordinate; see~\cite[Proposition~4.1 and~(5.1)]{chen2025free} for the construction. We set $P_N:=\mu_{\mcl I}^{\otimes N}\otimes\nu^{\otimes M_N}$. For $t\in\R_+$ and $q\in\mcl Q_\infty^2$, we define $H_N^{t,q}:\{-1,1\}^N\times[-1,1]^{M_N}\times\supp\fR\to\R$ and $\bar F_N:\R_+\times\mcl Q_\infty^2\to\R$ by
\begin{align}
    H_N^{t,q}(\sigma,v,\alpha)&:=\sqrt{2t}H_N(\sigma,v)+\sqrt2W_N^q(\sigma,v,\alpha)-Nq_1(1)-(t+q_2(1))|v|^2,\label{e.H^t,q_N=}\\
    \bar F_N(t,q)&:=-\frac1N\E\log\iint\exp(H_N^{t,q}(\sigma,v,\alpha))\,\d P_N(\sigma,v)\,\d\fR(\alpha).\label{e.F_N(t,q)=}
\end{align}
The terms subtracted in~\eqref{e.H^t,q_N=} equal one-half of the variance of the Gaussian part.

For a compactly supported probability measure $\rho$ on $\R$, we define $\psi_\rho:\mcl Q_\infty\to\R$ by
\begin{equation}
\label{e.def.psi.rho}
    \psi_\rho(q):=-\E\log\iint\exp\Ll(\sqrt2w^q(\alpha)\tau-q(1)\tau^2\Rr)\,\d\rho(\tau)\,\d\fR(\alpha),
\end{equation}
where, conditionally on $\fR$, the process $w^q$ is centered Gaussian with covariance $\E[w^q(\alpha)w^q(\alpha')\mid\fR]=q(\alpha\wedge\alpha')$, as in~\eqref{e.E[WW]=}.
For piecewise-constant $q$, the quantity $\psi_\rho(q)$ can be computed by recursive Gaussian integration; see~\cite[Theorem~5.25]{HJbook}. An equivalent representation through the Parisi PDE is given in~\cite[Lemma~6.4 and Proposition~6.6]{HJbook}.
For the fixed measures $\nu$ and $\mu_{\mcl I}$, we can thus define $\psi_\nu,\, \psi_{\mu_{\mcl I}}:\mcl Q_\infty\to\R$, which can be extended to $\mcl Q_1$ by Lemma~\ref{l.psi^vec_smooth}.
We also define $\psi_{\mu_{\mcl I}}^*:\mcl Q_\infty\to\R\cup\{+\infty\}$ by
\begin{equation}\label{e.conjugate.psiI}
    \psi_{\mu_{\mcl I}}^*(p):=\sup_{q\in\mcl Q_\infty}\Ll\{\la p, q \ra_{L^2}-\psi_{\mu_{\mcl I}}(q)\Rr\},\qquad\forall p\in\mcl Q_\infty.
\end{equation}
Here and throughout, $\la p,q\ra_{L^2}:=\int_0^1p(s) q(s)\d s$ denotes the pairing of two real-valued paths.

\begin{theorem}[Free energy of the bipartite model]\label{t.main}
For every $t \ge 0$ and $q=(q_1,q_2) \in \mcl Q_2^2$, we have
\begin{equation}  
\label{e.main.1}
\lim_{N \to +\infty} \bar F_N(t,q) = \sup_{p \in \mcl Q_\infty}\Ll\{\la p, q_1\ra_{L^2}+\lambda \psi_\nu(q_2+tp)-\psi^*_{\mu_{\mcl I}}(p)\Rr\}.
\end{equation}
Moreover, denoting this limit by $f(t,q)$, the function $f$ is the unique Lipschitz viscosity solution, in the sense of Definition~\ref{d.vs}, to
\begin{equation}\label{e.main.hj}
\begin{cases}
    \partial_t f - \int_0^1 \partial_{q_1} f \,\partial_{q_2} f\,\d r=0 ,\qquad &\text{on }(0,\infty)\times\mcl Q^2_2,
    \\
    f(0,q) = \psi_{\mu_{\mcl I}}(q_1)+\lambda \psi_\nu(q_2) ,\qquad &\text{for } q \in \mcl Q^2_2.
\end{cases}
\end{equation}
\end{theorem}

When $f$ is differentiable at $(t,q)$, its partial derivative $\partial_{q_i}f$ is an element of $L^2([0,1))$, represented by $r\mapsto\partial_{q_i}f(t,q,r)$. The nonlinearity in~\eqref{e.main.hj} is the $L^2$ scalar product of the two partial derivatives, namely $\int_0^1\partial_{q_1}f(t,q,r)\partial_{q_2}f(t,q,r)\d r$, and is thus the infinite-dimensional analogue of the bilinear nonlinearity appearing in \eqref{e.hj.bip}.  

We recall from \eqref{e.def.FN.delta0} that the limit of the free energy $\bar F_N^\mathrm{perc}$ defined in \eqref{e.def.FN.perc} is obtained by setting $q = 0$ in \eqref{e.main.1}. As already mentioned, the limit free energy of perceptron models with arbitrary bounded measurable activation functions was already determined in \cite{openai2026ising}.

\begin{remark}[Comparison with the Hopf formula]\label{r.compare_hopf}
The variational formula on the right-hand side of~\eqref{e.main.1} can be rewritten as
\begin{equation}\label{e.new_var}
    \sup_{p_1 \in \mcl Q_\infty}\inf_{q'_1\in\mcl Q_\infty}\Ll\{\psi_{\mu_{\mcl I}}(q'_1)+\lambda\psi_\nu(q_2+tp_1)+\la p_1, q_1-q'_1\ra_{L^2}\Rr\},
\end{equation}
which is the analogue on the space of paths of the formula~\eqref{e.def.u.half.supinf}. When $\psi_\nu$ is convex, one could adapt the approach in~\cite[Theorem~1.1]{chen2026ising} to obtain the Hopf formula for the limit free energy, reading
\begin{equation}\label{e.hopf}
    \sup_{p \in \mcl Q_\infty^2}\inf_{q'\in\mcl Q_\infty^2}\Ll\{\psi_{\mu_{\mcl I}}(q'_1)+\lambda\psi_\nu(q'_2)+\la p, q-q'\ra_{L^2}+t\la p_1,p_2\ra_{L^2}\Rr\},
\end{equation}
where $p=(p_1,p_2)$, $q'=(q'_1,q'_2)$, and $\la p, q-q'\ra_{L^2}=\sum_{i=1}^2 \la p_i, q_i-q'_i\ra_{L^2}$. As in Remark~\ref{r.hopf.half}, the two formulas~\eqref{e.new_var} and~\eqref{e.hopf} formally coincide when $\psi_\nu$ is convex: optimizing first over $q'_2$ and then over $p_2$ in~\eqref{e.hopf} produces the term $\lambda\psi_\nu(q_2+tp_1)$. However, for the formula~\eqref{e.hopf} to be valid, the convexity of both $\psi_{\mu_{\mcl I}}$ and $\psi_\nu$ is required, as can be seen by considering the case $t = 0$. In the present setting, the function $\psi_\nu$ need not be convex. For instance, convexity fails when $\nu$ is a probability measure on $\{-1,1\}$ whose mean has absolute value in $(1/\sqrt 3,1)$; see~\cite[Exercise~6.7 and solution]{HJbook}. Theorem~\ref{t.main} shows that~\eqref{e.new_var} gives the limit free energy and the viscosity solution of~\eqref{e.main.hj} regardless of whether $\psi_\nu$ is convex.
\end{remark}

\subsection{Outline of the proof and organization of the paper}
\label{s.intro.outline}

To identify the limiting free energy of the bipartite model, we proceed in the following way. On the one hand, we use \cite{mourrat2021nonconvex} to obtain that, for $t>0$ and $q\in\mcl Q_2^2$, 
\begin{equation}
\label{e.lower}
f(t,q) \le \liminf_{N\to\infty}\bar F_N(t,q),
\end{equation}
where $f$ is the solution to \eqref{e.main.hj}. For the converse bound, we follow the approach developed in \cite{chen2026ising}, which builds upon the cavity computations of \cite{chen2024ms, chen2025free}, and we show that
\begin{align}\label{e.bounds}
\limsup_{N \to \infty} \bar F_N(t,q) \leq \psi_{\mu_{\mcl I}}(q_1+t\bar p_2)+\lambda\psi_\nu(q_2+t\bar p_1)-t\la \bar p_1,\bar p_2\ra_{L^2},
\end{align}
for some $\bar p\in\mcl Q_\infty^2$ satisfying
\begin{align}\label{e.critical_pt}
    \bar p_1=\partial_q\psi_{\mu_{\mcl I}}(q_1+t\bar p_2)
    \qquad\text{and}\qquad
    \bar p_2=\lambda\partial_q\psi_\nu(q_2+t\bar p_1).
\end{align}
In order to relate the two bounds thus obtained, we show in Section~\ref{s.paths} that the solution to \eqref{e.main.hj} admits the variational representation
\begin{align}\label{e.f=sup}
    f(t,q) = \sup_{p \in \mcl Q_\infty}\Ll\{\la p, q_1\ra_{L^2}+\lambda \psi_\nu(q_2+tp)-\psi^*_{\mu_{\mcl I}}(p)\Rr\}.
\end{align}
Using the first relation in~\eqref{e.critical_pt} and the convexity of $\psi_{\mu_{\mcl I}}$, we also observe that
\begin{align}\label{e.psi^*=}
    \psi^*_{\mu_{\mcl I}}(\bar p_1) =\la \bar p_1,q_1+t\bar p_2\ra_{L^2} - \psi_{\mu_{\mcl I}}(q_1+t\bar p_2),
\end{align}
and therefore
\begin{align*}
    \text{r.h.s.\ of~\eqref{e.bounds}} \stackrel{\eqref{e.psi^*=}}{=} \la \bar p_1, q_1\ra_{L^2}+\lambda \psi_\nu(q_2+t\bar p_1)-\psi^*_{\mu_{\mcl I}}(\bar p_1) \stackrel{\eqref{e.f=sup}}{\leq} f(t,q).
\end{align*}
This and \eqref{e.lower}--\eqref{e.bounds} imply the desired identity~\eqref{e.main.1}.

For multi-species models in which every species consists of centered Ising spins, a bound of the form~\eqref{e.bounds} is proved in~\cite[Proposition~8.1]{chen2026ising}. The model we study here is not literally covered by this statement, since the spins of the second layer are sampled from~$\nu$. However, as we explain in Section~\ref{s.bounds}, the argument given there applies with very little modification. Alternatively, one could first show that the free energy converges as $N$ tends to infinity, by adapting the results of~\cite{chen2026vectorformula, ho2026convergence} to our multi-species context, and then appeal to the critical-point representation of the limit obtained in~\cite{chen2024ms, chen2025free}.

\smallskip

The rest of the paper is organized as follows. In Section~\ref{s.finite}, we prove Theorem~\ref{t.hopf.half}, and we extend it to equations posed on finite-dimensional cones. Section~\ref{s.enriched_properties} collects regularity properties of the enriched free energy and of the cascade transforms $\psi_{\mu_{\mcl I}}$ and $\psi_\nu$. In Section~\ref{s.paths}, we pass to the limit along dyadic discretizations, and obtain the variational formula on the space of paths. Section~\ref{s.bounds} completes the proof of Theorem~\ref{t.main}, following the outline above. In Section~\ref{s.moving_domain}, we show that the variational formula remains valid for equations posed on a domain of paths that shrinks as time increases. This is used in Section~\ref{s.general_perceptron}, where we discuss the extension of this approach to general Ising and spherical perceptrons, and to even more general models.

\smallskip

\noindent \textbf{Statement of AI use.} The ideas presented in this paper, except for the contents of Section~\ref{s.moving_domain}, were fully developed in early September 2026, prior to the release of \cite{aronow2026negative, openai2026ising, openai2026spherical}. Our starting point was the conjecture in \cite{gyorgyi1997parisi, gyorgyi2000beyond, montanari2026replica} for the spherical perceptron, which initially surprised us because we could not see how to reduce it to the classical cases under which we expected a variational formula to be valid. We observed that we could combine the arguments in \cite{chen2026ising} and in \cite{chen2026spherical} to obtain a variational formula for the large-$N$ limit of
\begin{equation*}
\frac 1 N \E \log \int \exp \Ll( \beta H_N(\sigma, v) \Rr) \, \d \mu_{\mcl I}^{\otimes N}(\sigma) \, \d \mu_{\mathbb{S}_{M_N}}(v),
\end{equation*}
where $H_N(\sigma,v)$ is as in \eqref{e.def.HN.bip}, and $\mu_{\mathbb{S}_{M_N}}$ is the uniform measure on $\mathbb{S}_{M_N} := \{v \in \R^{M_N} \ : \ |v|^2 = M_N\}$. Dividing this limit variational formula by $\beta$, and formally sending $\beta$ to infinity, we obtained a candidate formula for the large-$N$ limit of
\begin{equation*}
N^{-1}  \E \max_{\sigma \in \{-1,1\}^N, v \in \mathbb{S}_{M_N}} H_N(\sigma,v) = N^{-\frac 3 2} \, \E \max_{v \in \mathbb{S}_{M_N}} \sum_{i = 1}^N |g_{i\bullet} \cdot v|,
\end{equation*}
with the notation $g_{i\bullet} = (g_{ij})_{j \le M_N} \in \R^{M_N}$. The paper \cite{montanari2026replica} also proposes a variational formula for this limit (if we ignore the fact that the activation function is assumed bounded in \cite{montanari2026replica}). With assistance from ChatGPT 5.6, we could observe that the variational formula obtained as explained above and the one from \cite{montanari2026replica} are indeed equal. An inspection of the argument connecting the two formulas led us to wonder whether Theorem~\ref{t.hopf.half} could be correct. We then obtained a rigorous proof of Theorem~\ref{t.hopf.half} as well as the other results presented here, with assistance from ChatGPT 5.6. We have manually verified all statements and proofs, and we have polished the paper with assistance from ChatGPT 6 Astra and Fable 5.1.

\noindent \textbf{Acknowledgments.} We would like to thank Jeanne Boursier, Cédric Gerbelot and Victor Issa for stimulating discussions. HBC acknowledges funding from the NYU Shanghai Start-Up Fund and support from the NYU–ECNU Institute of Mathematical Sciences at NYU Shanghai. JCM acknowledges support from the MSCA grant SLOHD (101203974).

\section{Variational formula in finite dimensions}
\label{s.finite}

In this section, we prove Theorem~\ref{t.hopf.half}, and we then extend it to equations posed on finite-dimensional cones, in preparation for the passage to the space of paths in Section~\ref{s.paths}.

\subsection{Viscosity solutions}
\label{s.viscosity}

We use a single notion of viscosity solution in finite and infinite dimensions, following~\cite[Definitions~1.2 and~1.4]{chen2022hamilton}. Let $\cH$ be a separable Hilbert space, let $V\subseteq(0,\infty)\times\cH$, and let $\H:\cH\to\R$ be continuous. A function $\phi:V\to\R$ is called \textbf{smooth} if, for every $(t,x)\in V$, there is a unique element of $\R\times\cH$, denoted by $(\partial_t\phi(t,x),\nabsub{x}\phi(t,x))$, such that
\begin{equation*}
    \phi(s,y)-\phi(t,x)=\partial_t\phi(t,x)(s-t)+\la \nabsub{x}\phi(t,x),y-x\ra_{\cH}+O\Ll(|s-t|^2+|y-x|_{\cH}^2\Rr)
\end{equation*}
as $(s,y)\in V$ tends to $(t,x)$, and the map $(t,x)\mapsto(\partial_t\phi(t,x),\nabsub{x}\phi(t,x))$ from $V$ to $\R\times\cH$ is continuous.

\begin{definition}[Viscosity solutions]\label{d.vs_general}
A continuous function $u:V\to\R$ is a \textbf{viscosity subsolution} of $\partial_tu-\H(\nabsub{x}u)=0$ if, for every smooth $\phi:V\to\R$ and every local maximum $(t,x)\in V$ of $u-\phi$, we have
\begin{equation*}
    \partial_t\phi(t,x)-\H(\nabsub{x}\phi(t,x))\leq0.
\end{equation*}
A \textbf{viscosity supersolution} satisfies the reverse inequality at every local minimum of $u-\phi$. A \textbf{viscosity solution} is both a viscosity subsolution and a viscosity supersolution.
\end{definition}

If $V$ is open, then the differential of a smooth function is its Fr\'echet differential. If moreover $\cH$ is finite-dimensional, then every $C^2$ function on $V$ is smooth, and every smooth function is $C^1$. In this case, Definition~\ref{d.vs_general} agrees with the classical definition of viscosity solutions, in which the test functions are $C^1$, since the latter can equivalently be formulated with $C^\infty$ test functions; see for instance~\cite[Exercise~3.4 and solution]{HJbook}. If $V=(0,\infty)\times C$ for a closed convex cone $C\subseteq\cH$ satisfying $\overline{C-C}=\cH$, then the differential is unique also at boundary points; see~\cite[Remark~1.3]{chen2022hamilton}. 

\subsection{The case of the whole space}
\label{s.open}

Throughout this subsection, we place ourselves in the setting of Theorem~\ref{t.hopf.half}, and we write $g:=\psi_1^*$. We first record the properties of $g$ that will be used. The function $g$ is convex and lower semicontinuous, being a supremum of affine functions. It is proper: for every $p\in\R^d$, we have $g(p)\geq-\psi_1(0)>-\infty$, and if $p_0$ is a subgradient of the finite convex function $\psi_1$ at the origin, then $g(p_0)=-\psi_1(0)<+\infty$. Its effective domain $\operatorname{dom}g:=\{p\in\R^d:g(p)<+\infty\}$ is bounded: for every $p\in\R^d$ and $a\geq0$, the Lipschitz continuity of $\psi_1$ gives
\begin{equation*}
    g(p)\geq a|p|^2-\psi_1(ap)\geq a|p|\big(|p|-L_1\big)-\psi_1(0),
\end{equation*}
and letting $a$ tend to infinity shows that
\begin{equation}\label{e.dom.bound}
    \operatorname{dom}g\subseteq\{p\in\R^d:|p|\leq L_1\}.
\end{equation}
In particular, the supremum in~\eqref{e.def.u.half} can be restricted to the closed ball of radius $L_1$ centered at the origin. Finally, since $\psi_1$ is convex and continuous, the Fenchel--Moreau theorem gives
\begin{equation}\label{e.fenchel.moreau}
    g^*(x):=\sup_{p\in\R^d}\{p\cdot x-g(p)\}=\psi_1(x),\qquad\forall x\in\R^d;
\end{equation}
see for instance~\cite[Theorem~2.5]{HJbook}. For $p\in\operatorname{dom}g$ and $z=(t,x_1,x_2)\in\R_+\times\R^d\times\R^d$, we set
\begin{equation}\label{e.partial_objective}
    w_p(z):=p\cdot x_1+\psi_2(x_2+tp)-g(p),
\end{equation}
so that $w_p$ is the explicit solution~\eqref{e.sol.affine} with $c=g(p)$, and the function $u$ in~\eqref{e.def.u.half} is $u=\sup_{p\in\operatorname{dom}g}w_p$. We use Definition~\ref{d.vs_general} with ambient Hilbert space $\R^d\times\R^d$, with $V:=(0,\infty)\times\R^d\times\R^d$, and with the Hamiltonian $\H(a_1,a_2):=a_1\cdot a_2$. For a function $\phi$ of $z=(t,x_1,x_2)\in\R\times\R^d\times\R^d$, we write $\nablatx\phi:=(\dr_t\phi,\nabsub{x_1}\phi,\nabsub{x_2}\phi)$ for its full gradient, in contrast with $\nabla\phi=(\nabsub{x_1}\phi,\nabsub{x_2}\phi)$, which only involves the space variables.

\begin{proof}[Proof of Theorem~\ref{t.hopf.half}]
\smallskip
\noindent\emph{Step 1: Attainment, Lipschitz continuity, and the initial condition.}
We fix $z=(t,x_1,x_2)\in\R_+\times\R^d\times\R^d$. The function $p\mapsto p\cdot x_1+\psi_2(x_2+tp)-g(p)$, extended by $-\infty$ outside $\operatorname{dom}g$, is upper semicontinuous on the closed ball of radius $L_1$, and it is not identically $-\infty$ there. Hence $u(z)$ is finite, the supremum in~\eqref{e.def.u.half} is attained, and the set of maximizers
\begin{equation*}
    \mcl A(z):=\{p\in\operatorname{dom}g:w_p(z)=u(z)\}
\end{equation*}
is nonempty and compact. For every $p\in\operatorname{dom}g$, we have $|p|\leq L_1$ by~\eqref{e.dom.bound}, and thus, for every $(s,y_1,y_2)\in\R_+\times\R^d\times\R^d$,
\begin{equation*}
    |w_p(t,x_1,x_2)-w_p(s,y_1,y_2)|\leq L_1|x_1-y_1|+L_2|x_2-y_2|+L_1L_2|t-s|.
\end{equation*}
Taking the supremum over $p$ yields~\eqref{e.lip.u.half}. At $t=0$, the formula~\eqref{e.def.u.half} and the identity~\eqref{e.fenchel.moreau} give $u(0,x_1,x_2)=g^*(x_1)+\psi_2(x_2)=\psi_1(x_1)+\psi_2(x_2)$.

\smallskip
\noindent\emph{Step 2: The subsolution inequality.}
Let $z=(t,x_1,x_2)\in V$, and let $\phi$ be a smooth function such that $u-\phi$ has a local maximum at $z$. We fix $p\in\mcl A(z)$. Since $w_p\leq u$ and $w_p(z)=u(z)$, the function $w_p-\phi$ also has a local maximum at $z$, that is,
\begin{equation}\label{e.sub.touch}
    \phi(z')-\phi(z)\geq w_p(z')-w_p(z)
\end{equation}
for every $z'$ in a neighborhood of $z$. The function $w_p$ is affine in $x_1$ with slope $p$, and it is constant along the direction $(1,0,-p)$: for every $y\in\R^d$ and every $h\in\R$ sufficiently small,
\begin{equation*}
    w_p(t,x_1+hy,x_2)-w_p(z)=hp\cdot y
    \qquad\text{and}\qquad
    w_p(t+h,x_1,x_2-hp)-w_p(z)=0.
\end{equation*}
We apply~\eqref{e.sub.touch} with $z'=(t,x_1+hy,x_2)$ and with $z'=(t+h,x_1,x_2-hp)$, divide by $h$, and let $h$ tend to $0$ from above and from below. This yields $\nabsub{x_1}\phi(z)\cdot y=p\cdot y$ for every $y\in\R^d$, and $\partial_t\phi(z)-p\cdot\nabsub{x_2}\phi(z)=0$. Hence $\nabsub{x_1}\phi(z)=p$, and
\begin{equation*}
    \partial_t\phi(z)-\nabsub{x_1}\phi(z)\cdot \nabsub{x_2}\phi(z)=\partial_t\phi(z)-p\cdot\nabsub{x_2}\phi(z)=0.
\end{equation*}
In particular, $u$ is a viscosity subsolution of the equation in~\eqref{e.hj.bip}. We note that this step uses no regularity of $\psi_2$ beyond continuity.

\smallskip
\noindent\emph{Step 3: The supersolution inequality when $\psi_2\in C^1(\R^d)$.}
In this step, we assume that $\psi_2$ is continuously differentiable. For every $p\in\operatorname{dom}g$, the function $w_p$ is then continuously differentiable on $V$, with
\begin{equation}\label{e.partial_active_gradient}
    \nablatx w_p(t',x_1',x_2')=\big(p\cdot\nabla\psi_2(x_2'+t'p),\,p,\,\nabla\psi_2(x_2'+t'p)\big)\in\R\times\R^d\times\R^d,
\end{equation}
and the right-hand side depends continuously on $(p,t',x_1',x_2')$. We fix $z=(t,x_1,x_2)\in V$, and for $p\in\mcl A(z)$, we set $b_p:=\nabla\psi_2(x_2+tp)$ and $\gamma_p:=\nablatx w_p(z)=(p\cdot b_p,p,b_p)$.

We first derive the first-order optimality condition satisfied by $p\in\mcl A(z)$, namely
\begin{equation}\label{e.partial_subgradient}
    g(r)-g(p)\geq(x_1+tb_p)\cdot(r-p),\qquad\forall r\in\R^d,
\end{equation}
which states that $x_1+tb_p$ is a subgradient of $g$ at $p$. There is nothing to prove if $g(r)=+\infty$. Otherwise, for every $\eps\in(0,1]$, the point $p_\eps:=p+\eps(r-p)$ belongs to $\operatorname{dom}g$, by convexity, and $w_{p_\eps}(z)\leq u(z)=w_p(z)$. Together with the convexity inequality $g(p_\eps)\leq g(p)+\eps(g(r)-g(p))$, this gives
\begin{equation*}
    \eps\,(r-p)\cdot x_1+\psi_2\big(x_2+tp+\eps t(r-p)\big)-\psi_2(x_2+tp)\leq g(p_\eps)-g(p)\leq\eps\big(g(r)-g(p)\big).
\end{equation*}
Dividing by $\eps$ and letting $\eps\downarrow0$ yields~\eqref{e.partial_subgradient}. For $p,r\in\mcl A(z)$, we write~\eqref{e.partial_subgradient} as is and also with $p$ and $r$ interchanged, and add the two inequalities. Since $t>0$, this gives
\begin{equation}\label{e.partial_monotonicity}
    (p-r)\cdot(b_p-b_r)\geq0,\qquad\forall p,r\in\mcl A(z).
\end{equation}

Let now $\phi$ be a smooth function such that $u-\phi$ has a local minimum at $z$. We claim that
\begin{equation}\label{e.Dphi(z).v<}
    \nablatx\phi(z)\cdot v\leq\max_{p\in\mcl A(z)}\gamma_p\cdot v,\qquad\forall v\in\R\times\R^d\times\R^d.
\end{equation}
We fix $v\in\R\times\R^d\times\R^d$, and for every $h>0$ small enough that $z+hv\in V$, we pick $p_h\in\mcl A(z+hv)$. Since $u(z+hv)=w_{p_h}(z+hv)$ and $u(z)\geq w_{p_h}(z)$, the local minimum of $u-\phi$ at $z$ and the mean value theorem give
\begin{equation}\label{e.sup.mvt}
    \phi(z+hv)-\phi(z)\leq u(z+hv)-u(z)\leq w_{p_h}(z+hv)-w_{p_h}(z)=h\,\nablatx w_{p_h}(z+\theta_hhv)\cdot v
\end{equation}
for some $\theta_h\in(0,1)$. By~\eqref{e.dom.bound}, there is a sequence $h_n\downarrow0$ such that $p_{h_n}$ converges to some $p\in\R^d$. We verify that $p\in\mcl A(z)$. Writing $z_n:=z+h_nv=(t_n,x_{1,n},x_{2,n})$, the continuity of $\psi_2$ and $u$ gives
\begin{equation*}
    g(p_{h_n})=p_{h_n}\cdot x_{1,n}+\psi_2(x_{2,n}+t_np_{h_n})-u(z_n)\xrightarrow[n\to\infty]{}p\cdot x_1+\psi_2(x_2+tp)-u(z),
\end{equation*}
and the lower semicontinuity of $g$ yields $g(p)\leq p\cdot x_1+\psi_2(x_2+tp)-u(z)$, that is, $w_p(z)\geq u(z)$. Hence $p\in\mcl A(z)$. By the continuity of the right-hand side of~\eqref{e.partial_active_gradient}, we have $\lim_{n\to\infty}\nablatx w_{p_{h_n}}(z+\theta_{h_n}h_nv)=\gamma_p$. Dividing~\eqref{e.sup.mvt} by $h_n$ and letting $n\to\infty$, we obtain $\nablatx\phi(z)\cdot v\leq\gamma_p\cdot v$, which proves~\eqref{e.Dphi(z).v<}.

The set $\Gamma_z:=\{\gamma_p:p\in\mcl A(z)\}$ is compact, being the image of the compact set $\mcl A(z)$ under a continuous map, so its convex hull $\operatorname{co}\Gamma_z$ is compact and convex, by Carath\'eodory's theorem. If $\nablatx\phi(z)$ did not belong to $\operatorname{co}\Gamma_z$, then by the separating hyperplane theorem, there would exist $v\in\R\times\R^d\times\R^d$ such that $\nablatx\phi(z)\cdot v>\max_{\gamma\in\Gamma_z}\gamma\cdot v$, contradicting~\eqref{e.Dphi(z).v<}. Hence there exist $m\in\N$, $p_1,\dots,p_m\in\mcl A(z)$, and $\alpha_1,\dots,\alpha_m\geq0$ with $\sum_{i=1}^m\alpha_i=1$ such that
\begin{equation*}
    \nablatx\phi(z)=\sum_{i=1}^m\alpha_i\gamma_{p_i}=\sum_{i=1}^m\alpha_i(p_i\cdot b_{p_i},p_i,b_{p_i}).
\end{equation*}
Consequently, we have
\begin{align*}
    \partial_t\phi(z)-\nabsub{x_1}\phi(z)\cdot \nabsub{x_2}\phi(z)
    &=\sum_{i=1}^m\alpha_ip_i\cdot b_{p_i}-\Big(\sum_{i=1}^m\alpha_ip_i\Big)\cdot\Big(\sum_{j=1}^m\alpha_jb_{p_j}\Big)\\
    &=\frac12\sum_{i,j=1}^m\alpha_i\alpha_j(p_i-p_j)\cdot(b_{p_i}-b_{p_j})
    \stackrel{\eqref{e.partial_monotonicity}}{\geq}0.
\end{align*}
This is the viscosity supersolution inequality.

\smallskip
\noindent\emph{Step 4: Conclusion.}
We now remove the differentiability assumption on $\psi_2$. For $\eps>0$, let $\psi_{2,\eps}:=\psi_2*\rho_\eps$, where $\rho_\eps$ is a smooth probability density supported in the ball of radius $\eps$. The function $\psi_{2,\eps}$ is smooth and $L_2$-Lipschitz, and it satisfies $\sup_{\R^d}|\psi_{2,\eps}-\psi_2|\leq L_2\eps$. Let $u_\eps$ be defined by~\eqref{e.def.u.half} with $\psi_{2,\eps}$ in place of $\psi_2$. By Steps~1 to~3, the function $u_\eps$ is a viscosity solution to the equation in~\eqref{e.hj.bip}, and comparing~\eqref{e.def.u.half} for $\psi_2$ and for $\psi_{2,\eps}$ directly gives
\begin{equation}\label{e.unif.conv}
    \sup_{\R_+\times\R^d\times\R^d}|u_\eps-u|\leq\sup_{\R^d}|\psi_{2,\eps}-\psi_2|\leq L_2\eps.
\end{equation}
We conclude by a standard stability argument. Let $\phi$ be a smooth function such that $u-\phi$ has a local minimum at $z_0\in V$. Replacing $\phi$ by $z\mapsto\phi(z)-|z-z_0|^2$, which does not change $\nablatx\phi(z_0)$, we may assume that this minimum is strict, so that there is a closed ball $B\subseteq V$ centered at $z_0$ with $u-\phi>(u-\phi)(z_0)$ on $B\setminus\{z_0\}$. For every $\eps>0$, let $z_\eps$ be a minimum point of $u_\eps-\phi$ on~$B$. By~\eqref{e.unif.conv}, every limit point of $(z_\eps)_{\eps\downarrow0}$ is a minimum point of $u-\phi$ on $B$, and thus equals $z_0$. Hence $z_\eps$ tends to $z_0$, so for $\eps$ small, $z_\eps$ lies in the interior of $B$, and the supersolution inequality for $u_\eps$ at $z_\eps$ gives $\partial_t\phi(z_\eps)-\nabsub{x_1}\phi(z_\eps)\cdot\nabsub{x_2}\phi(z_\eps)\geq0$. Letting $\eps\downarrow0$ and using the continuity of $\nablatx\phi$, we obtain the supersolution inequality for $u$ at $z_0$. Together with Steps~1 and~2, this shows that $u$ is a Lipschitz continuous viscosity solution to~\eqref{e.hj.bip}. Finally, Lipschitz continuous viscosity solutions to~\eqref{e.hj.bip} are unique, by the comparison principle for Hamilton--Jacobi equations with locally Lipschitz nonlinearities; see for instance~\cite[Corollary~3.7]{HJbook}.
\end{proof}

\begin{remark}\label{r.general_g}
The proof above uses $\psi_1$ only through $g=\psi_1^*$, and the only properties of $g$ that were used are that it is proper, lower semicontinuous, and convex, with bounded effective domain. Conversely, if $g:\R^d\to\R\cup\{+\infty\}$ is any function with these properties, then $\psi_1:=g^*$ is finite, convex and Lipschitz continuous, with Lipschitz constant $\sup_{p\in\operatorname{dom}g}|p|$, and $\psi_1^*=g^{**}=g$, by the Fenchel--Moreau theorem. Theorem~\ref{t.hopf.half} can therefore be rephrased as follows: for every such~$g$ and every Lipschitz continuous $\psi_2:\R^d\to\R$, all the conclusions of Theorem~\ref{t.hopf.half} hold, with $L_1:=\sup_{p\in\operatorname{dom}g}|p|$, for the function
\begin{equation}\label{e.partial_hopf_E}
    (t,x_1,x_2)\mapsto\sup_{p\in\operatorname{dom}g}\big\{p\cdot x_1+\psi_2(x_2+tp)-g(p)\big\},
\end{equation}
which is thus the unique Lipschitz continuous viscosity solution to the equation in~\eqref{e.hj.bip} with initial condition $g^*(x_1)+\psi_2(x_2)$. We use this formulation in the proof of Proposition~\ref{p.partial_hopf_cone} below, and in Section~\ref{s.spherical}.
\end{remark}

\subsection{Extension to cones}
\label{s.cone}

Let $E$ be a finite-dimensional Euclidean space. For a closed convex cone $K$ in a Hilbert space $\cH$, we define the \textbf{dual cone} $K^*:=\{h\in\cH:\la h,x\ra_\cH\geq0,\ \forall x\in K\}$, and we say that a function $h$ defined on a subset $S$ of $\cH$ is \textbf{$K^*$-increasing} if $h(y)\geq h(x)$ whenever $x,y\in S$ and $y-x\in K^*$. Let $C\subseteq E$ be a pointed closed convex cone with nonempty interior; we recall that a cone $C$ is said to be pointed if $C \cap (-C) = \{0\}$. We have $(C^2)^*=(C^*)^2$, and we say that a function $v:\R_+\times C^2\to\R$ is $(C^*)^2$-increasing if $v(t,\cdot)$ is $(C^*)^2$-increasing for every $t\geq0$. Suppose that $\psi_1,\psi_2:C\to\R$ are $C^*$-increasing and Lipschitz, with constants $L_1,L_2$ respectively, and that $\psi_1$ is convex. We define $\psi_{1,C}^*:C\to\R\cup\{+\infty\}$ and $u:\R_+\times C^2\to\R$ by
\begin{align}
    \psi_{1,C}^*(p)&:=\sup_{y\in C}\{p\cdot y-\psi_1(y)\},\label{e.cone_conjugate}\\
    u(t,x_1,x_2)&:=\sup_{p\in C}\{p\cdot x_1+\psi_2(x_2+tp)-\psi_{1,C}^*(p)\}.\label{e.partial_hopf_cone}
\end{align}

\begin{proposition}\label{p.partial_hopf_cone}
The supremum in~\eqref{e.partial_hopf_cone} is attained. The function $u$ is $(C^*)^2$-increasing, and for every $t,s\geq0$ and $x_1,x_2,y_1,y_2\in C$, we have
\begin{equation}\label{e.lip.cone}
    |u(t,x_1,x_2)-u(s,y_1,y_2)|\leq L_1|x_1-y_1|+L_2|x_2-y_2|+L_1L_2|t-s|.
\end{equation}
Moreover, $u$ is a viscosity solution of
\begin{equation}\label{e.bipartite_E}
    \partial_tu-\nabsub{x_1}u\cdot \nabsub{x_2}u=0\qquad\text{on }(0,\infty)\times(\operatorname{int}C)^2,
\end{equation}
with initial condition $u(0,x_1,x_2)=\psi_1(x_1)+\psi_2(x_2)$ for all $x_1,x_2\in C$, and it is the only Lipschitz, $(C^*)^2$-increasing function $v:\R_+\times C^2\to\R$ with these two properties.
\end{proposition}

\begin{proof}
Let $g:E\to\R\cup\{+\infty\}$ be equal to $\psi_{1,C}^*$ on $C$ and to $+\infty$ outside of $C$. As a supremum of affine functions on the closed convex set $C$, this function is convex and lower semicontinuous. For every $p\in C$ and $a\geq0$, the choice of $y=ap$ in~\eqref{e.cone_conjugate} and the Lipschitz continuity of $\psi_1$ give $g(p)\geq a|p|^2-\psi_1(ap)\geq a|p|(|p|-L_1)-\psi_1(0)$. Letting $a\to\infty$ shows that
\begin{equation}\label{e.cone_conjugate_bound}
    \operatorname{dom}g\subseteq\{p\in C:|p|\leq L_1\}.
\end{equation}

We next show the Fenchel--Moreau identity on the cone,
\begin{equation}\label{e.cone_biconjugate}
    \psi_1(x)=\sup_{p\in C}\{p\cdot x-\psi_{1,C}^*(p)\}=g^*(x),\qquad\forall x\in C,
\end{equation}
which in particular shows that $g$ is proper. The inequality $\geq$ in~\eqref{e.cone_biconjugate} follows from the definition~\eqref{e.cone_conjugate}. Conversely, let first $x\in\operatorname{int}C$, and let $p\in E$ be a subgradient of $\psi_1$ at $x$ relative to~$C$, meaning that $\psi_1(y)\geq\psi_1(x)+p\cdot(y-x)$ for every $y\in C$. Such a subgradient exists since $\psi_1$ is convex and finite on a neighborhood of $x$; see~\cite[Proposition~2.11]{HJbook}. For every $h\in C^*$ and $\eps>0$ small enough that $x-\eps h\in C$, the monotonicity of $\psi_1$ and the subgradient inequality give
\begin{equation*}
    0\leq\psi_1(x)-\psi_1(x-\eps h)\leq\eps\, p\cdot h.
\end{equation*}
Thus $p\in(C^*)^*=C$, by the bipolar identity for closed convex cones; see~\cite[Section~1.1.1]{chen2022cones}. The subgradient inequality also shows that the supremum in~\eqref{e.cone_conjugate} is attained at $y=x$, that is, $\psi_{1,C}^*(p)=p\cdot x-\psi_1(x)$. Hence $\psi_1(x)=p\cdot x-\psi_{1,C}^*(p)\leq g^*(x)$, and~\eqref{e.cone_biconjugate} holds on $\operatorname{int}C$. Both sides of~\eqref{e.cone_biconjugate} are Lipschitz on $C$, the right-hand side because of~\eqref{e.cone_conjugate_bound}, and $\operatorname{int}C$ is dense in $C$, so the identity extends to $C$.

We now extend $\psi_2$ to $E$ by setting
\begin{equation*}
    \tilde\psi_2(x):=\inf_{y\in C}\big\{\psi_2(y)+L_2|x-y|\big\},\qquad\forall x\in E,
\end{equation*}
which defines an $L_2$-Lipschitz function on $E$ that coincides with $\psi_2$ on $C$. By Theorem~\ref{t.hopf.half} in the form of Remark~\ref{r.general_g}, applied to $g$ and $\tilde\psi_2$, and by~\eqref{e.cone_conjugate_bound}, the function
\begin{equation*}
    \tilde u(t,x_1,x_2):=\sup_{p\in\operatorname{dom}g}\big\{p\cdot x_1+\tilde\psi_2(x_2+tp)-g(p)\big\},\qquad(t,x_1,x_2)\in\R_+\times E^2,
\end{equation*}
is a viscosity solution of $\partial_t\tilde u-\nabsub{x_1}\tilde u\cdot\nabsub{x_2}\tilde u=0$ on $(0,\infty)\times E^2$, with initial condition $\tilde u(0,x_1,x_2)=g^*(x_1)+\tilde\psi_2(x_2)$; moreover, the supremum defining $\tilde u$ is attained, and $\tilde u$ satisfies the estimate~\eqref{e.lip.u.half}. Since $\operatorname{dom}g\subseteq C$ and $C$ is a convex cone, for every $x_2\in C$, $t\geq0$ and $p\in\operatorname{dom}g$, we have $x_2+tp\in C$, and thus $\tilde\psi_2(x_2+tp)=\psi_2(x_2+tp)$. As $\psi_{1,C}^*(p)=+\infty$ for $p\in C\setminus\operatorname{dom}g$, comparing with~\eqref{e.partial_hopf_cone} shows that $\tilde u=u$ on $\R_+\times C^2$. Hence the supremum in~\eqref{e.partial_hopf_cone} is attained, the estimate~\eqref{e.lip.cone} holds, and by~\eqref{e.cone_biconjugate}, we have $u(0,x_1,x_2)=\psi_1(x_1)+\psi_2(x_2)$ for every $x_1,x_2\in C$. Since the viscosity property is local, and $(0,\infty)\times(\operatorname{int}C)^2$ is an open subset of $(0,\infty)\times E^2$ on which $u$ and $\tilde u$ coincide, the function $u$ is a viscosity solution of~\eqref{e.bipartite_E}.

If $x_1,x_2,y_1,y_2\in C$ are such that $y_1-x_1$ and $y_2-x_2$ belong to $C^*$, then for every $p\in C$, we have $p\cdot y_1\geq p\cdot x_1$ and $\psi_2(y_2+tp)\geq\psi_2(x_2+tp)$, so that $u(t,y_1,y_2)\geq u(t,x_1,x_2)$. Hence $u$ is $(C^*)^2$-increasing.

Finally, the Hamiltonian $\H:E^2\to\R$, $\H(a_1,a_2):=a_1\cdot a_2$, is locally Lipschitz, and its restriction to $C^2$ is $(C^*)^2$-increasing: if $a,b\in C^2$ and $a-b\in(C^*)^2$, then
\begin{equation*}
    \H(a)-\H(b)=(a_1-b_1)\cdot a_2+b_1\cdot(a_2-b_2)\geq0.
\end{equation*}
The comparison principle of~\cite[Theorem~1.2(1)--(2)]{chen2022cones}, applied in the Hilbert space $E^2$ with the cone $C^2$, therefore gives the uniqueness of Lipschitz, $(C^*)^2$-increasing viscosity solutions of~\eqref{e.bipartite_E} with prescribed initial condition.
\end{proof}

\section{Regularity properties of the enriched free energy}
\label{s.enriched_properties}

We return to the bipartite model of the introduction and collect the regularity and compactness properties used below, as well as the convexity of $\psi_{\mu_{\mcl I}}$.

The $N$ Ising spins and the $M_N$ $\nu$-distributed spins form two species with one-dimensional single-spin laws $\mu_{\mcl I}$ and $\nu$. The covariance of $H_N$ is $N$ times the product of their overlaps $N^{-1}\sigma\cdot\sigma'$ and $N^{-1}v\cdot v'$, whose absolute values are bounded by $1$ and $M_N/N$, respectively. Thus the Gaussian interpolation estimates of~\cite[Proposition~4.1]{chen2024ms}, with normalization by $N$ instead of $N+M_N$, yield the following result.

\begin{proposition}[{\cite[Proposition~4.1]{chen2024ms}}]\label{p.F_N_smooth}
For every $N\in\N$, $t,t'\in\R_+$, and $q,q'\in\mcl Q_\infty^2$, we have
\begin{equation}\label{e.F_lip}
    |\bar F_N(t,q)-\bar F_N(t',q')|\leq\frac{M_N}{N}|t-t'|+|q_1-q_1'|_{L^1}+\frac{M_N}{N}|q_2-q_2'|_{L^1}.
\end{equation}
Consequently, $\bar F_N$ extends uniquely to $\R_+\times\mcl Q_1^2$.
\end{proposition}

In the notation of Section~\ref{s.cone}, the dual cone of $\mcl Q_2$ in $L^2([0,1))$ is
\begin{equation}\label{e.dual_Q2}
    \mcl Q_2^*:=\{h\in L^2([0,1)):\la h,q\ra_{L^2}\geq0,\ \forall q\in\mcl Q_2\}.
\end{equation}
A function $g:\mcl Q_2\to\R$ is \textbf{$\mcl Q_2^*$-increasing} if
\begin{equation}\label{e.increasing_Q2}
    g(q')\geq g(q)\quad\text{whenever }q,q'\in\mcl Q_2\text{ and }q'-q\in\mcl Q_2^*.
\end{equation}

We use the notion of Gateaux differentiability from~\cite[Definition~2.2]{chen2025free}. For $q\in\mcl Q_2$, we define the set of admissible directions by
\begin{equation*}
    \operatorname{Adm}(\mcl Q_2,q):=\left\{e\in L^2([0,1)):\exists a>0,\ \forall a'\in[0,a],\ q+a'e\in\mcl Q_2\right\}.
\end{equation*}
A function $g:\mcl Q_2\to\R$ is \textbf{Gateaux differentiable} at $q\in\mcl Q_2$ if, for every $e\in\operatorname{Adm}(\mcl Q_2,q)$, the limit
\begin{equation*}
    g'(q,e):=\lim_{a\downarrow0}\frac{g(q+ae)-g(q)}{a}
\end{equation*}
exists and there is a unique $y\in L^2([0,1))$ such that
\begin{equation*}
    g'(q,e)=\la y,e\ra_{L^2},\qquad\forall e\in\operatorname{Adm}(\mcl Q_2,q).
\end{equation*}
We call $y$ the \textbf{Gateaux derivative} of $g$ at $q$ and write $\partial_qg(q):=y$. We recall the regularity properties needed below.

\begin{lemma}[{\cite[Corollary~5.2]{chen2025free}}]\label{l.psi^vec_smooth}
Let $\rho$ be a probability measure with support on $[-1,1]$. The function $\psi_\rho$ defined in \eqref{e.def.psi.rho} extends uniquely to $\mcl Q_1$, with $\psi_\rho(0)=0$, and
\begin{equation}\label{e.psi_lip}
    |\psi_\rho(q)-\psi_\rho(q')|\leq|q-q'|_{L^1},\qquad\forall q,q'\in\mcl Q_1.
\end{equation}
Its restriction to $\mcl Q_2$ is Gateaux differentiable and $\mcl Q_2^*$-increasing, with
\begin{equation}\label{e.bound.der.psi}
    \partial_q\psi_\rho(q)\in\mcl Q_{\infty,\leq1},\qquad\forall q\in\mcl Q_2.
\end{equation}
The derivative extends continuously to $\mcl Q_1$, and
\begin{equation}\label{e.continuity.der.psi}
    |\partial_q\psi_\rho(q)-\partial_q\psi_\rho(q')|_{L^1}\leq16|q-q'|_{L^1},\qquad\forall q,q'\in\mcl Q_1.
\end{equation}
\end{lemma}

We will also use the following compactness consequence of~\cite[Lemma~3.4]{chen2025free}.
\begin{lemma}\label{l.compact_embed_paths}
For every $c\geq0$ and $r\in[1,\infty)$, the set $\mcl Q_{\infty,\leq c}$ is compact in $L^r([0,1))$.
\end{lemma}

The cascade recursion separates over independent spin coordinates; see~\cite[Lemma~4.11]{chen2024ms}. Thus, defining $\psi:\mcl Q_1^2\to\R$ by $\psi(q):=\psi_{\mu_{\mcl I}}(q_1)+\lambda\psi_\nu(q_2)$, we have
\begin{equation}\label{e.psi=sumlambdapsi}
    \bar F_N(0,q)=\psi_{\mu_{\mcl I}}(q_1)+\frac{M_N}{N}\psi_\nu(q_2)
    \qquad\text{and}\qquad
    \lim_{N\to\infty}\bar F_N(0,q)=\psi(q).
\end{equation}
For our model, we use the convexity only of $\psi_{\mu_{\mcl I}}$, provided by the following result.

\begin{proposition}[{\cite[Proposition~2.2]{chen2026ising}}, \cite{ho2026concavity}]\label{p.convex.psi.circ}
The function $\psi_{\mu_{\mcl I}}:\mcl Q_1\to\R$ is convex: for every $q,q'\in\mcl Q_1$ and $\theta\in[0,1]$, we have 
\begin{equation*}  \psi_{\mu_{\mcl I}}((1-\theta)q+\theta q')\leq(1-\theta)\psi_{\mu_{\mcl I}}(q)+\theta\psi_{\mu_{\mcl I}}(q').
\end{equation*}
\end{proposition}

\section{Variational formula on the space of paths}
\label{s.paths}

In this section, we extend Theorem~\ref{t.hopf.half} to the space of paths, by passing to the limit in Proposition~\ref{p.partial_hopf_cone} along dyadic discretizations. The main result of this section is Theorem~\ref{t.partial_hopf_paths} below, which identifies the right-hand side of~\eqref{e.main.1} with the viscosity solution of~\eqref{e.main.hj}.

We apply Definition~\ref{d.vs_general} to $\mcl Q_2^2\subseteq L^2([0,1);\R^2)$. The identity $\overline{\mcl Q_2^2-\mcl Q_2^2}=L^2([0,1);\R^2)$
ensures uniqueness of the differential in Definition~\ref{d.vs_general}; see~\cite[Remark~1.3]{chen2022hamilton}.

We set $\xi:\R^2\to\R$ by $\xi(a):=a_1a_2$, and choose a Lipschitz regularization $\bar\xi:\R_+^2\to\R$ agreeing with $\xi$ on the ball of radius $2\vee(1+\lambda)$. We require $\bar\xi$ to be coordinatewise increasing and, for every $b\in\R_+^2$, the function $a\mapsto\bar\xi(a+b)-\bar\xi(a)$ to be coordinatewise increasing. Such a regularization exists by~\cite[Lemma~4.4 and Remark~4.5]{chen2022hamilton}. Let $\H:L^2([0,1);\R^2)\to\R$ be such that
\begin{equation}\label{e.def_H_spin_glass}
    \H(h):=\inf\left\{\int_0^1\bar\xi(r(u))\,\d u:r\in\mcl Q_2^2,\ r-h\in(\mcl Q_2^*)^2\right\}.
\end{equation}
By~\cite[Lemma~4.6]{chen2022hamilton}, the map $\H$ is Lipschitz and $(\mcl Q_2^*)^2$-increasing, and for every $r\in\mcl Q_2^2$, we have $\H(r)=\int_0^1\bar\xi(r(u))\,\d u$.
We use the notion of viscosity solution from~\cite[Definition~4.3]{chen2022hamilton}.

\begin{definition}\label{d.vs}
A continuous function $f:\R_+\times\mcl Q_2^2\to\R$ is a \textbf{viscosity solution of~\eqref{e.main.hj}} if $f(0,\cdot)=\psi$ and it solves $\partial_tf-\H(\partial_qf)=0$ in the sense of Definition~\ref{d.vs_general} on $(0,\infty)\times\mcl Q_2^2$, with the regularized Hamiltonian $\H$ from~\eqref{e.def_H_spin_glass}.
\end{definition}

Existence, uniqueness in the Lipschitz class, and convergence of finite-dimensional approximations follow from~\cite[Theorem~4.7(1)]{chen2022hamilton}. In particular, the solution in Definition~\ref{d.vs} does not depend on the choice of $\bar\xi$; see~\cite[Remark~4.10]{chen2022hamilton}.\footnote{To use the unit Lipschitz normalization in~\cite[Theorem~4.7]{chen2022hamilton}, it suffices to divide the unknown function and its initial condition by $1+\lambda$, and to replace $\xi$ by $a\mapsto\xi((1+\lambda)a)/(1+\lambda)$. This is legitimate since~\eqref{e.psi_lip} gives an $L^1$-Lipschitz constant at most $1+\lambda$ for $\psi$.} The scalar two-species setting is treated by the same argument; see also~\cite[Remark~7.7]{chen2026ising}.

For $j\in\N$, we equip $E_j:=\R^{2^j}$ with the scalar product $\la x,y\ra_j:=2^{-j}\sum_{k=1}^{2^j}x_ky_k$, and set $C_j:=\{x\in E_j:0\leq x_1\leq\cdots\leq x_{2^j}\}$. We define the lift $\iota_j:E_j\to L^2([0,1))$ and the averaging map $\pi_j:L^2([0,1))\to E_j$ by
\begin{equation}\label{e.dyadic_maps}
    \iota_jx:=\sum_{k=1}^{2^j}x_k\one_{[(k-1)2^{-j},k2^{-j})}
    \qquad\text{and}\qquad
    (\pi_jq)_k:=2^j\int_{(k-1)2^{-j}}^{k2^{-j}}q(r)\,\d r,
    \quad k\in\{1,\dots,2^j\}.
\end{equation}
These maps preserve the increasing cones, $\iota_j$ is an isometry, and
\begin{equation}\label{e.dyadic_convergence}
    \lim_{j\to\infty}|\iota_j\pi_jq-q|_{L^2}=0,\qquad\forall q\in L^2([0,1)).
\end{equation}
For $\rho\in\{\mu_{\mcl I},\nu\}$, we define $\psi_{\rho,j}:C_j\to\R$ by $\psi_{\rho,j}(x):=\psi_\rho(\iota_jx)$. Its cone conjugate for $\rho=\mu_{\mcl I}$ is $\psi_{\mu_{\mcl I},j}^*:C_j\to\R\cup\{+\infty\}$, defined by
\begin{equation}\label{e.dyadic_conjugate}
    \psi_{\mu_{\mcl I},j}^*(p):=\sup_{x\in C_j}\{\la p,x\ra_j-\psi_{\mu_{\mcl I}}(\iota_jx)\}.
\end{equation}
We now consider the following finite-dimensional approximation of~\eqref{e.main.hj}:
\begin{equation}\label{e.dyadic_hj}
    \begin{cases}
        \partial_tf_j-\la \nabsub{x_1}f_j,\nabsub{x_2}f_j\ra_j=0,&\text{on }(0,\infty)\times(\operatorname{int}C_j)^2,\\
        f_j(0,x_1,x_2)=\psi_{\mu_{\mcl I},j}(x_1)+\lambda\psi_{\nu,j}(x_2),&\text{for }(x_1,x_2)\in C_j^2.
    \end{cases}
\end{equation}
Here the gradients $\nabsub{x_i}f_j$ are taken with respect to $\la\cdot,\cdot\ra_j$: at differentiability points, the derivative in direction $h\in E_j$ is $\la \nabsub{x_i}f_j,h\ra_j$. This problem is of the form considered in Proposition~\ref{p.partial_hopf_cone}, with the Euclidean space $E_j$, the cone $C_j$, and the functions $\psi_1=\psi_{\mu_{\mcl I},j}$ and $\psi_2=\lambda\psi_{\nu,j}$. The assumptions of this proposition are satisfied: the functions $\psi_{\rho,j}$ are $C_j^*$-increasing because $\iota_j(C_j^*)\subseteq\mcl Q_2^*$, as follows from the identity $\la\iota_jh,r\ra_{L^2}=\la h,\pi_jr\ra_j$; the function $\psi_{\mu_{\mcl I},j}$ is convex by Proposition~\ref{p.convex.psi.circ}; and Lemma~\ref{l.psi^vec_smooth} gives the Lipschitz bounds. Proposition~\ref{p.partial_hopf_cone} therefore ensures that~\eqref{e.dyadic_hj} admits a unique Lipschitz, $(C_j^*)^2$-increasing viscosity solution $f_j:\R_+\times C_j^2\to\R$, and that this solution is given by
\begin{equation}\label{e.dyadic_formula}
    f_j(t,x_1,x_2)=\sup_{p\in C_j}\{\la p,x_1\ra_j+\lambda\psi_{\nu,j}(x_2+tp)-\psi_{\mu_{\mcl I},j}^*(p)\}.
\end{equation}

\begin{theorem}\label{t.partial_hopf_paths}
The right-hand side of~\eqref{e.main.1} defines the unique Lipschitz viscosity solution $f:\R_+\times\mcl Q_2^2\to\R$ of~\eqref{e.main.hj}, in the sense of Definition~\ref{d.vs}. The supremum can be restricted to $\mcl Q_{\infty,\leq1}$ and is attained. For every $t,s\geq0$ and $q,y\in\mcl Q_2^2$, we have
\begin{equation}\label{e.formula_lip}
    |f(t,q)-f(s,y)|\leq|q_1-y_1|_{L^1}+\lambda|q_2-y_2|_{L^1}+\lambda|t-s|.
\end{equation}
\end{theorem}

\begin{proof}
\smallskip
\noindent\emph{Step 1: Restriction to $\mcl Q_{\infty,\leq1}$.}
Recall the conjugate $\psi_{\mu_{\mcl I}}^*$ defined in~\eqref{e.conjugate.psiI}. We first show that $\psi_{\mu_{\mcl I}}^*(p)<+\infty$ implies $p\in\mcl Q_{\infty,\leq1}$. If $p(r)>1$ for some $r\in[0,1)$, then~\eqref{e.psi_lip}, applied with $q_a:=a\one_{[r,1)}$ for $a\geq0$, gives
\begin{equation*}
    \psi_{\mu_{\mcl I}}^*(p)\geq a\int_r^1(p(u)-1)\,\d u,\qquad\forall a\geq0.
\end{equation*}
The integral is strictly positive by monotonicity of $p$, so $\psi_{\mu_{\mcl I}}^*(p)=+\infty$. Similarly, using $x=(0,\dots,0,a)\in C_j$ in~\eqref{e.dyadic_conjugate}, with $a\geq0$, shows that $\psi_{\mu_{\mcl I},j}^*(p)<+\infty$ implies $p_{2^j}\leq1$. Finally, $\psi_{\mu_{\mcl I}}^*(0)=\psi_{\mu_{\mcl I},j}^*(0)=0$, since $\psi_{\mu_{\mcl I}}(0)=0$ and~\eqref{e.bound.der.psi} gives $\psi_{\mu_{\mcl I}}(q)\geq0$ for every $q\in\mcl Q_\infty$.

\smallskip
\noindent\emph{Step 2: Identification of the limit of the dyadic solutions.}
By~\cite[Theorem~4.7(1)]{chen2022hamilton}, the unique viscosity solution $f$ of~\eqref{e.main.hj} satisfies
\begin{equation}\label{e.dyadic_solution_limit}
    f(t,q)=\lim_{j\to\infty}f_j(t,\pi_jq_1,\pi_jq_2),\qquad\forall(t,q)\in\R_+\times\mcl Q_2^2.
\end{equation}
We identify this limit with the right-hand side of~\eqref{e.main.1}. For every $p\in\mcl Q_\infty$, the definitions~\eqref{e.dyadic_conjugate} and~\eqref{e.conjugate.psiI}, together with the identity $\la\pi_jp,x\ra_j=\la p,\iota_jx\ra_{L^2}$ give
\begin{equation}\label{e.projected_conjugate_upper}
    \psi_{\mu_{\mcl I},j}^*(\pi_jp)
    =\sup_{x\in C_j}\{\la p,\iota_jx\ra_{L^2}-\psi_{\mu_{\mcl I}}(\iota_jx)\}
    \leq\psi_{\mu_{\mcl I}}^*(p).
\end{equation}
Using $\pi_jp\in C_j$ in~\eqref{e.dyadic_formula} and applying~\eqref{e.projected_conjugate_upper}, we obtain
\begin{equation*}
    f_j(t,\pi_jq_1,\pi_jq_2)
    \geq\la\pi_jp,\pi_jq_1\ra_j+\lambda\psi_\nu\big(\iota_j\pi_jq_2+t\iota_j\pi_jp\big)-\psi_{\mu_{\mcl I}}^*(p).
\end{equation*}
For $p$ with $\psi_{\mu_{\mcl I}}^*(p)<+\infty$, Step 1 gives $p\in\mcl Q_{\infty,\leq1}$. By~\eqref{e.dyadic_convergence}, we have $\lim_{j\to\infty}\la\pi_jp,\pi_jq_1\ra_j=\la p,q_1\ra_{L^2}$, while~\eqref{e.psi_lip} and~\eqref{e.dyadic_convergence} imply
\begin{equation*}
    \lim_{j\to\infty}\psi_\nu(\iota_j\pi_jq_2+t\iota_j\pi_jp)=\psi_\nu(q_2+tp).
\end{equation*}
Consequently, by~\eqref{e.dyadic_solution_limit}, we obtain
\begin{equation*}
    f(t,q)\geq\la p,q_1\ra_{L^2}+\lambda\psi_\nu(q_2+tp)-\psi_{\mu_{\mcl I}}^*(p).
\end{equation*}
Taking the supremum over $p\in\mcl Q_\infty$ gives $f(t,q)\geq\text{r.h.s.\ of \eqref{e.main.1}}$.

Conversely, Proposition~\ref{p.partial_hopf_cone} ensures that the supremum in~\eqref{e.dyadic_formula} is attained. For each $j\in\N$, we choose a maximizer $p_j\in C_j$ at $(t,\pi_jq_1,\pi_jq_2)$. Step 1 gives $0\leq(p_j)_1\leq\cdots\leq(p_j)_{2^j}\leq1$. Since the expression in braces in~\eqref{e.dyadic_formula} at $p=p_j$ is at least its value at $p=0$, and $\psi_{\mu_{\mcl I},j}^*(0)=0$, we have
\begin{equation*}
    \la p_j,\pi_jq_1\ra_j+\lambda\psi_\nu(\iota_j\pi_jq_2+t\iota_jp_j)-\psi_{\mu_{\mcl I},j}^*(p_j)
    \geq\lambda\psi_\nu(\iota_j\pi_jq_2).
\end{equation*}
Thus, using~\eqref{e.psi_lip}, $0\leq\iota_jp_j\leq1$, and $q_1\geq0$, we obtain
\begin{align}
    0\leq\psi_{\mu_{\mcl I},j}^*(p_j)
    \leq\la p_j,\pi_jq_1\ra_j+\lambda\Big(\psi_\nu(\iota_j\pi_jq_2+t\iota_jp_j)-\psi_\nu(\iota_j\pi_jq_2)\Big)\leq|q_1|_{L^1}+\lambda t.\label{e.dyadic_conjugate_bound}
\end{align}
By Lemma~\ref{l.compact_embed_paths}, we may pass to a subsequence for which $\lim_{j\to\infty}|\iota_jp_j-p|_{L^2}=0$ for some $p\in\mcl Q_{\infty,\leq1}$. For every $r\in\mcl Q_\infty$, taking $x=\pi_jr$ in~\eqref{e.dyadic_conjugate} gives
\begin{equation*}
    \psi_{\mu_{\mcl I},j}^*(p_j)\geq\la\iota_jp_j,\iota_j\pi_jr\ra_{L^2}-\psi_{\mu_{\mcl I}}(\iota_j\pi_jr).
\end{equation*}
Passing to the lower limit, using~\eqref{e.dyadic_convergence},~\eqref{e.psi_lip}, and the $L^2$ convergence of $\iota_jp_j$, and then taking the supremum over $r\in\mcl Q_\infty$, we obtain
\begin{equation}\label{e.projected_conjugate_lower}
    \liminf_{j\to\infty}\psi_{\mu_{\mcl I},j}^*(p_j)\geq\psi_{\mu_{\mcl I}}^*(p).
\end{equation}
In particular, $\psi_{\mu_{\mcl I}}^*(p)<+\infty$ by~\eqref{e.dyadic_conjugate_bound}. By~\eqref{e.dyadic_convergence} and the $L^2$ convergence of $\iota_jp_j$, we have
\begin{align*}
    \lim_{j\to\infty}\la p_j,\pi_jq_1\ra_j=\la p,q_1\ra_{L^2}\qquad
    \text{and}\qquad\lim_{j\to\infty}\psi_\nu(\iota_j\pi_jq_2+t\iota_jp_j)=\psi_\nu(q_2+tp),
\end{align*}
where the second identity uses the $L^1$-Lipschitz continuity in~\eqref{e.psi_lip}. Passing to the limit in~\eqref{e.dyadic_formula} therefore gives
\begin{equation*}
    f(t,q)\stackrel{\eqref{e.projected_conjugate_lower}}{\leq}
    \la p,q_1\ra_{L^2}+\lambda\psi_\nu(q_2+tp)-\psi_{\mu_{\mcl I}}^*(p)
    \leq\text{r.h.s.\ of \eqref{e.main.1}}.
\end{equation*}
Together with the reverse inequality, this identifies $f$ with the right-hand side of~\eqref{e.main.1}.

\smallskip
\noindent\emph{Step 3: Attainment and the Lipschitz estimate.}
On $\mcl Q_{\infty,\leq1}$, the function $\psi_{\mu_{\mcl I}}^*$ is lower semicontinuous in the $L^2$ topology, being a supremum of continuous affine functions. By~\eqref{e.psi_lip}, the map
\begin{equation*}
    p\mapsto\la p,q_1\ra_{L^2}+\lambda\psi_\nu(q_2+tp)-\psi_{\mu_{\mcl I}}^*(p)
\end{equation*}
is therefore upper semicontinuous on $\mcl Q_{\infty,\leq1}$. Since this set is compact by Lemma~\ref{l.compact_embed_paths} and the function is finite at $p=0$, the supremum is attained.

Finally, if $p\in\mcl Q_{\infty,\leq1}$ and $\psi_{\mu_{\mcl I}}^*(p)<+\infty$, then~\eqref{e.psi_lip} gives
\begin{align*}
    \big|\la p,q_1-y_1\ra_{L^2}+\lambda\psi_\nu(q_2+tp)-\lambda\psi_\nu(y_2+sp)\big|\leq|q_1-y_1|_{L^1}+\lambda|q_2-y_2|_{L^1}+\lambda|t-s|.
\end{align*}
Since this bound is uniform in~$p$, and $f(t,q)$ and $f(s,y)$ are suprema over the same set of~$p$, we obtain~\eqref{e.formula_lip}.
\end{proof}

\section{Proof of Theorem~\ref{t.main}}
\label{s.bounds}

In this section, we complete the proof of Theorem~\ref{t.main}, following the outline given in Section~\ref{s.intro.outline}. We first show in Proposition~\ref{p.HJ_bdd_ms} that the enriched free energy is asymptotically bounded from below by the solution to~\eqref{e.main.hj}. We then show in Proposition~\ref{p.crit_pt_bdd_ms} that it is asymptotically bounded from above by the value of an explicit functional at one of its critical points, following the argument developed in~\cite{chen2026ising} for multi-species models with centered Ising spins. Finally, we use the convexity of $\psi_{\mu_{\mcl I}}$ and the variational formula of Theorem~\ref{t.partial_hopf_paths} to show that this critical value is itself bounded from above by the solution to~\eqref{e.main.hj}, which completes the argument.

\begin{proposition}[Lower bound]\label{p.HJ_bdd_ms}
Let $f:\R_+\times\mcl Q_2^2\to\R$ be the unique Lipschitz viscosity solution to~\eqref{e.main.hj}, in the sense of Definition~\ref{d.vs}. For every $(t,q)\in\R_+\times\mcl Q_2^2$, we have
\begin{equation}\label{e.HJ_lower}
    f(t,q)\leq\liminf_{N\to\infty}\bar F_N(t,q).
\end{equation}
\end{proposition}

\begin{proof}
A bound of this form is proved in~\cite[Theorem~2.7]{mourrat2021nonconvex} for bipartite models in which the two layers have the same number of spins. There, the enriched free energy is parametrized by pairs of probability measures on $\R_+$ instead of pairs of paths; we pass from one description to the other by identifying a probability measure with its quantile function. In order to reduce our setting to that of~\cite{mourrat2021nonconvex}, we pad the two layers with zeros. We set $K_N:=N\vee M_N$ and $c_N:=K_N/N$, so that $c_N$ converges to $c:=1\vee\lambda$ as $N$ tends to infinity. We identify each $\sigma\in\{-1,1\}^N$ with the vector of $\R^{K_N}$ obtained by appending $K_N-N$ coordinates equal to zero, we proceed similarly for each $v\in[-1,1]^{M_N}$, and we denote by $\widehat P_N$ the image of the measure $P_N$ under this identification, which is a probability measure on $\R^{K_N}\times\R^{K_N}$. We let $\widehat F_N:\R_+\times\mcl Q_2^2\to\R$ be the enriched free energy of the bipartite model with $K_N$ spins in each layer and with reference measure $\widehat P_N$, as defined in~\cite[Section~2]{mourrat2021nonconvex}. In this model, the Gaussian Hamiltonian is normalized by $K_N^{-1/2}$ in place of $N^{-1/2}$, and the free energy by $K_N^{-1}$ in place of $N^{-1}$. Since the additional coordinates vanish, the Gaussian Hamiltonian of this model is $c_N^{-1/2}H_N(\sigma,v)$, while the fields indexed by the cascade are the same as in~\eqref{e.E[WW]=}. The enriched Hamiltonian of this model at $(s,q)$ is therefore the Hamiltonian in~\eqref{e.H^t,q_N=} at $(s/c_N,q)$, and we obtain that, for every $(s,q)\in\R_+\times\mcl Q_2^2$,
\begin{equation}\label{e.padding_identity}
    \widehat F_N(s,q)=\frac1{c_N}\,\bar F_N\Ll(\frac s{c_N},q\Rr).
\end{equation}

The measures $\widehat P_N$ are not product measures of the form assumed in the statement of~\cite[Theorem~2.7]{mourrat2021nonconvex}. However, as explained in~\cite[the paragraph following~(2.19)]{mourrat2021nonconvex}, the proof of this theorem only uses two properties of the reference measures: that the squared norm of each layer of spins is bounded by the number of spins in the layer, and that the free energy at time zero converges as the number of spins tends to infinity. Both properties are valid here. Indeed, under $\widehat P_N$, the squared norms of the two layers are equal to $N$ and bounded by $M_N$ respectively, and they are thus both bounded by $K_N$. Moreover, by~\eqref{e.padding_identity} and~\eqref{e.psi=sumlambdapsi}, we have for every $q\in\mcl Q_2^2$ that
\begin{equation*}
    \lim_{N\to\infty}\widehat F_N(0,q)=\widehat\psi(q),\qquad\text{where}\quad\widehat\psi(q):=\frac{\psi(q)}{c}=\frac{\psi_{\mu_{\mcl I}}(q_1)+\lambda\psi_\nu(q_2)}{1\vee\lambda}.
\end{equation*}
The proof of~\cite[Theorem~2.7]{mourrat2021nonconvex} therefore yields that, for every $(s,q)\in\R_+\times\mcl Q_2^2$,
\begin{equation}\label{e.padded_lower}
    \widehat f(s,q)\leq\liminf_{N\to\infty}\widehat F_N(s,q),
\end{equation}
where $\widehat f:\R_+\times\mcl Q_2^2\to\R$ is the solution to~\eqref{e.main.hj} with $\widehat\psi$ in place of $\psi$ as the initial condition. In~\cite[Section~3]{mourrat2021nonconvex}, this solution is defined as the limit of the solutions to finite-dimensional approximations of the equation; by the arguments of~\cite[Proposition~4.13 and Theorem~4.7(1)]{chen2022hamilton} (see also~\cite[Remark~7.7]{chen2026ising}), it is the unique Lipschitz viscosity solution in the sense of Definition~\ref{d.vs}, with $\widehat\psi$ in place of $\psi$.

We now fix $(t,q)\in\R_+\times\mcl Q_2^2$. By~\eqref{e.padding_identity} with $s=ct$, we have $\bar F_N(ct/c_N,q)=c_N\widehat F_N(ct,q)$, while the time regularity in~\eqref{e.F_lip} and the convergence of $c_N$ to $c$ imply that $\bar F_N(t,q)-\bar F_N(ct/c_N,q)$ tends to zero as $N$ tends to infinity. Combining these two observations with~\eqref{e.padded_lower}, we obtain
\begin{equation*}
    \liminf_{N\to\infty}\bar F_N(t,q)=c\,\liminf_{N\to\infty}\widehat F_N(ct,q)\geq c\,\widehat f(ct,q).
\end{equation*}
It remains to identify the right-hand side of this inequality. Since the nonlinearity in~\eqref{e.main.hj} is homogeneous of degree two, we see directly from the definition of viscosity solutions that the function $(t',q')\mapsto c\,\widehat f(ct',q')$ is a Lipschitz viscosity solution to~\eqref{e.main.hj} with initial condition $c\,\widehat\psi=\psi$.\footnote{In Definition~\ref{d.vs}, this rescaling replaces the regularization $\bar\xi$ by $a\mapsto c^2\,\bar\xi(a/c)$, which is also an admissible regularization since $c\geq1$; we recall that the solution does not depend on this choice.} By uniqueness, this function is $f$; in particular, we have $c\,\widehat f(ct,q)=f(t,q)$, and the proof of~\eqref{e.HJ_lower} is complete.
\end{proof}

We now turn to the converse bound. As explained in Section~\ref{s.intro.outline}, we show that the upper limit of the free energy is bounded from above by the value of an explicit functional at one of its critical points. In the course of the argument, we will need to vary the number of spins in the second layer independently of~$N$. For all integers $N,M\geq1$, we denote by $\bar F_{N,M}:\R_+\times\mcl Q_1^2\to\R$ the function defined as $\bar F_N$ in~\eqref{e.H^t,q_N=}--\eqref{e.F_N(t,q)=}, but with $M$ in place of $M_N$ throughout. In other words, $\bar F_{N,M}$ is the enriched free energy of the model with $N$ Ising spins and $M$ spins sampled from~$\nu$, still normalized by~$N$. In particular, we have $\bar F_N=\bar F_{N,M_N}$, and the estimate~\eqref{e.F_lip} is valid for $\bar F_{N,M}$, with $M$ in place of $M_N$.

Some of the results of~\cite{chen2024ms} that we rely on are stated for multi-species models in which the overlaps and the free energy are normalized by the total number of spins. In order to pass from one convention to the other, we set $\widetilde\xi(r_1,r_2):=(1+\lambda)r_1r_2$, and for all integers $N,M\geq1$, we let $\widetilde F_{N,M}:\R_+\times\mcl Q_1^2\to\R$ denote the enriched free energy of the two-species model of~\cite{chen2024ms} with $N$ spins sampled from $\mu_{\mcl I}$, with $M$ spins sampled from $\nu$, and with covariance function $\widetilde\xi$. Explicitly, writing $K:=N+M$, the Gaussian Hamiltonian of this model has covariance
\begin{equation*}
    K\,\widetilde\xi\Ll(\frac{\sigma\cdot\sigma'}{K},\frac{v\cdot v'}{K}\Rr)=\frac{(1+\lambda)N}{K}\,\E\Ll[H_N(\sigma,v)H_N(\sigma',v')\Rr],
\end{equation*}
the fields indexed by the cascade are the same as in~\eqref{e.E[WW]=}, and the free energy is normalized by $K$ instead of $N$. The enriched Hamiltonian of this model at $(t,q)$ thus has the same law as the Hamiltonian in~\eqref{e.H^t,q_N=} at $(t',q)$, with $t':=(1+\lambda)Nt/K$, and we obtain that, for every $(t,q)\in\R_+\times\mcl Q_1^2$,
\begin{equation}\label{e.total_size_normalization}
    \widetilde F_{N,M}(t,q)=\frac{N}{N+M}\,\bar F_{N,M}\Ll(\frac{(1+\lambda)N}{N+M}\,t,q\Rr).
\end{equation}
The factor $1+\lambda$ in the definition of $\widetilde\xi$ is chosen so that the rescaling of the time variable in~\eqref{e.total_size_normalization} becomes trivial in the regime of interest, in which $M/N$ is close to $\lambda$.

For every $t>0$ and $q\in\mcl Q_1^2$, we define $\sP_{t,q}:\mcl Q_\infty^2\to\R$ by
\begin{equation}\label{e.sP_lambda,t,q}
    \sP_{t,q}(p):=\psi_{\mu_{\mcl I}}(q_1+tp_2)+\lambda\psi_\nu(q_2+tp_1)-t\la p_1,p_2\ra_{L^2},
\end{equation}
and we say that $p\in\mcl Q_\infty^2$ is a \textbf{critical point} associated with $(t,q)$ if
\begin{equation}\label{e.critical_rel}
    p_1=\partial_q\psi_{\mu_{\mcl I}}(q_1+tp_2)
    \qquad\text{and}\qquad
    p_2=\lambda\partial_q\psi_\nu(q_2+tp_1).
\end{equation}
We recall from Lemma~\ref{l.psi^vec_smooth} that the derivatives appearing in~\eqref{e.critical_rel} are well-defined on $\mcl Q_1$. Formally, the relations in~\eqref{e.critical_rel} express that the derivative of the functional $\sP_{t,q}$ vanishes at $p$, which explains the terminology.

\begin{proposition}[Upper bound by a critical value]\label{p.crit_pt_bdd_ms}
For every $t>0$ and $q\in\mcl Q_1^2$, there exists a critical point $p\in\mcl Q_{\infty,\leq1}\times\mcl Q_{\infty,\leq\lambda}$ associated with $(t,q)$ such that
\begin{equation}\label{e.critical_value}
    \limsup_{N\to\infty}\bar F_N(t,q)\leq\sP_{t,q}(p).
\end{equation}
\end{proposition}

\begin{proof}
For multi-species models in which every species consists of centered Ising spins, a bound of this form is proved in~\cite[Proposition~8.1]{chen2026ising}. The argument given there applies with very little change to the present setting, in which the spins of the second layer are sampled from~$\nu$. It consists in first treating the case in which the ratio $\lambda$ is rational, by grouping the spins into blocks and viewing the model as a vector spin glass, for which the bound is obtained by means of cavity computations; and in then deducing the general case by approximation. In this argument, the assumption that the spins are Ising spins is only used to ensure that the self-overlap is constant, which simplifies the cavity computations; we explain in Step~1 below the modification that is required in the absence of this property. For the reader's convenience, and since our normalization differs from that of~\cite{chen2026ising}, we give the details of the rest of the argument. We fix $t>0$ throughout the proof.

\smallskip
\noindent\emph{Step 1: The case of rational $\lambda$ and bounded paths.}
We assume that $\lambda=a/b$ for some integers $a,b\geq1$ and that $q\in\mcl Q_\infty^2$, and we show that the statement of the proposition holds in this case. To do so, we group the spins into blocks made of~$b$ Ising spins and of $a$ spins sampled from~$\nu$, and we view the model as a vector spin glass whose spins are these blocks. We set $I_1:=\{1,\dots,b\}$, $I_2:=\{b+1,\dots,b+a\}$, and $\mu:=\mu_{\mcl I}^{\otimes b}\otimes\nu^{\otimes a}$, which is a probability measure on $\{-1,1\}^b\times[-1,1]^a\subseteq\R^{a+b}$. We define $\xi^\Vec:\R^{a+b}\to\R$ by
\begin{equation}\label{e.block_maps}
    \xi^\Vec(y):=\frac1b\Ll(\sum_{i\in I_1}y_i\Rr)\Ll(\sum_{i\in I_2}y_i\Rr),
\end{equation}
and we let $L:\R^2\to\R^{a+b}$ be the linear map such that $(Lx)_i=x_1$ for $i\in I_1$ and $(Lx)_i=x_2$ for $i\in I_2$. For a pair of paths $q=(q_1,q_2)$, we write $Lq$ for the path $s\mapsto L(q(s))$. For every integer $n\geq1$, we consider configurations $\tau=(\tau_{ik})_{i\in I_1\cup I_2,\,1\leq k\leq n}\in(\{-1,1\}^b\times[-1,1]^a)^n$, sampled according to the product measure $\mu^{\otimes n}$. Given independent standard Gaussian random variables $(g_{ik,j\ell})_{i\in I_1,\,j\in I_2,\,1\leq k,\ell\leq n}$, we set
\begin{equation}\label{e.vector_Hamiltonian}
    H_n^\Vec(\tau):=\frac1{\sqrt{bn}}\sum_{\substack{i\in I_1,\,j\in I_2\\ 1\leq k,\ell\leq n}}g_{ik,j\ell}\,\tau_{ik}\tau_{j\ell}.
\end{equation}
For two configurations $\tau$ and $\tau'$, we denote by $R(\tau,\tau')\in\R^{a+b}$ the vector with coordinates $R_i(\tau,\tau'):=n^{-1}\sum_{k=1}^n\tau_{ik}\tau'_{ik}$, so that
\begin{equation*}
    \E\Ll[H_n^\Vec(\tau)H_n^\Vec(\tau')\Rr]=n\,\xi^\Vec\big(R(\tau,\tau')\big).
\end{equation*}
For every $q^\Vec\in\mcl Q_\infty^{a+b}$, conditionally on $\fR$, we let $W_n^{q^\Vec}$ be a centered Gaussian process indexed by $(\tau,\alpha)$, independent of $H_n^\Vec$, with covariance
\begin{equation*}
    \E\Ll[W_n^{q^\Vec}(\tau,\alpha)W_n^{q^\Vec}(\tau',\alpha')\mid\fR\Rr]=n\sum_{i\in I_1\cup I_2}q_i^\Vec(\alpha\wedge\alpha')R_i(\tau,\tau').
\end{equation*}
In analogy with~\eqref{e.H^t,q_N=}--\eqref{e.F_N(t,q)=}, we define
\begin{multline}\label{e.vector_free_energy}
    \bar F_n^\Vec(t,q^\Vec):=-\frac1n\E\log\iint\exp\bigg(\sqrt{2t}H_n^\Vec(\tau)+\sqrt2W_n^{q^\Vec}(\tau,\alpha)\\
    -nt\,\xi^\Vec\big(R(\tau,\tau)\big)-n\sum_{i\in I_1\cup I_2}q_i^\Vec(1)R_i(\tau,\tau)\bigg)\,\d\mu^{\otimes n}(\tau)\,\d\fR(\alpha).
\end{multline}
This is the enriched free energy of a vector spin glass with spins in $\R^{a+b}$: the single-spin measure~$\mu$ has compact support; the covariance of $H_n^\Vec$ is $n$ times a polynomial of the matrix of overlaps $(n^{-1}\sum_{k=1}^n\tau_{ik}\tau'_{jk})_{i,j\in I_1\cup I_2}$, which here only involves the diagonal entries of this matrix; and the fields indexed by the cascade are those associated with the diagonal matrix-valued path with diagonal entries $q^\Vec$.

We now state the bound that we need for this vector spin glass. For every $q^\Vec=(q_i^\Vec)_{i\in I_1\cup I_2}\in\mcl Q_\infty^{a+b}$, we set
\begin{equation}\label{e.vector_psi}
    \psi^\Vec(q^\Vec):=\sum_{i\in I_1}\psi_{\mu_{\mcl I}}(q_i^\Vec)+\sum_{i\in I_2}\psi_\nu(q_i^\Vec),
\end{equation}
and we denote by $\partial_q\psi^\Vec(q^\Vec)$ the family of paths whose coordinate indexed by $i$ is $\partial_q\psi_{\mu_{\mcl I}}(q_i^\Vec)$ if $i\in I_1$, and $\partial_q\psi_\nu(q_i^\Vec)$ if $i\in I_2$. As for~\eqref{e.psi=sumlambdapsi}, the fact that the cascade recursion separates over independent spin coordinates implies that $\bar F_n^\Vec(0,\cdot)=\psi^\Vec$ for every $n\geq1$. We also write $\nabla\xi^\Vec$ for the gradient of $\xi^\Vec$, and we set $\theta^\Vec(y):=y\cdot\nabla\xi^\Vec(y)-\xi^\Vec(y)$ for every $y\in\R^{a+b}$. We claim that there exists $p^\Vec=(p_i^\Vec)_{i\in I_1\cup I_2}\in\mcl Q_\infty^{a+b}$ such that
\begin{equation}\label{e.vector_critical}
    p^\Vec=\partial_q\psi^\Vec\Ll(Lq+t\nabla\xi^\Vec(p^\Vec)\Rr)
\end{equation}
and
\begin{equation}\label{e.vector_upper}
    \limsup_{n\to\infty}\bar F_n^\Vec(t,Lq)\leq\psi^\Vec\Ll(Lq+t\nabla\xi^\Vec(p^\Vec)\Rr)-t\int_0^1\theta^\Vec(p^\Vec(s))\,\d s.
\end{equation}
If $\nu$ were equal to $\mu_{\mcl I}$, then this vector spin glass would be a special case of the model studied in~\cite[Sections~4 to~6]{chen2026ising}, and the existence of such a path $p^\Vec$ would be~\cite[Theorem~6.1]{chen2026ising}. The proof of this theorem applies essentially verbatim to the present model. The only difference is that, since the coordinates of the spins indexed by $I_2$ are now sampled from~$\nu$, the self-overlap $R(\tau,\tau)$ is no longer constant; as in~\cite[Section~6]{chen2025free}, one then needs to add to the Hamiltonian a small perturbation that forces the self-overlap to concentrate, and the cavity computations are otherwise unchanged. This modification is explained in~\cite[Appendix~A]{chen2026ising}, in the more general setting of vector spin glasses whose covariance may depend on the full matrix of overlaps, and whose enriched free energy is therefore parametrized by paths with values in the set of positive semidefinite matrices; see~\cite[Theorem~A.3]{chen2026ising}.

We now express~\eqref{e.vector_critical} and~\eqref{e.vector_upper} in terms of the bipartite model. We set
\begin{equation*}
    p_1:=\frac1b\sum_{i\in I_1}p_i^\Vec,\qquad p_2:=\frac1b\sum_{i\in I_2}p_i^\Vec,\qquad\text{and}\qquad p:=(p_1,p_2)\in\mcl Q_\infty^2.
\end{equation*}
By~\eqref{e.block_maps}, the coordinates of $\nabla\xi^\Vec(p^\Vec)$ indexed by $I_1$ are all equal to $p_2$, and those indexed by $I_2$ are all equal to $p_1$; in other words, we have $Lq+t\nabla\xi^\Vec(p^\Vec)=L(q_1+tp_2,q_2+tp_1)$. The identity~\eqref{e.vector_critical} thus states that $p_i^\Vec=\partial_q\psi_{\mu_{\mcl I}}(q_1+tp_2)$ for every $i\in I_1$, and that $p_i^\Vec=\partial_q\psi_\nu(q_2+tp_1)$ for every $i\in I_2$. In view of the definition of $p_1$ and $p_2$, and since $\lambda=a/b$, this implies that $p$ satisfies~\eqref{e.critical_rel}. That is, the pair $p$ is a critical point associated with $(t,q)$, and by~\eqref{e.bound.der.psi}, it belongs to $\mcl Q_{\infty,\leq1}\times\mcl Q_{\infty,\leq\lambda}$. Moreover, since $\xi^\Vec$ is homogeneous of degree two, we have $\theta^\Vec(p^\Vec)=\xi^\Vec(p^\Vec)=b\,p_1p_2$, and the right-hand side of~\eqref{e.vector_upper} is therefore equal to
\begin{equation*}
    b\,\psi_{\mu_{\mcl I}}(q_1+tp_2)+a\,\psi_\nu(q_2+tp_1)-tb\la p_1,p_2\ra_{L^2}=b\,\sP_{t,q}(p).
\end{equation*}
We have thus shown that
\begin{equation}\label{e.vector_upper_bip}
    \limsup_{n\to\infty}\bar F_n^\Vec(t,Lq)\leq b\,\sP_{t,q}(p).
\end{equation}

We now relate this vector spin glass to the bipartite model. When $N=bn$ and $M=an$, relabeling the coordinates of $(\sigma,v)\in\{-1,1\}^{bn}\times[-1,1]^{an}$ as $(\tau_{ik})_{i\in I_1\cup I_2,\,1\leq k\leq n}$ identifies the Hamiltonian $H_{bn}$ in~\eqref{e.def.HN.bip} with the Hamiltonian $H_n^\Vec$ in~\eqref{e.vector_Hamiltonian}, and the field $W_{bn}^q$ in~\eqref{e.E[WW]=} with $W_n^{Lq}$. Using that $R_i(\tau,\tau)=1$ for every $i\in I_1$, one can check that the correction terms also coincide, and since the two free energies are normalized by $bn$ and by $n$ respectively, we get
\begin{equation}\label{e.block_identity}
    \bar F_{bn,an}(t,q)=\frac1b\,\bar F_n^\Vec(t,Lq).
\end{equation}
For a general sequence $(M_N)_{N\in\N}$ such that $M_N/N$ converges to $a/b$, the identity~\eqref{e.block_identity} is replaced by an approximate version, which we borrow from~\cite{chen2024ms}. In the notation of that paper, the function $\widetilde F_{N,M_N}$ is the free energy of a two-species model with $N+M_N$ spins, whose species proportions converge to $\frac{b}{a+b}$ and $\frac{a}{a+b}$. Moreover, since
\begin{equation*}
    (a+b)\,\widetilde\xi\Ll(\frac{y_1}{a+b},\frac{y_2}{a+b}\Rr)=\frac{y_1y_2}{b},
\end{equation*}
the vector spin glass associated with this model in~\cite[Section~3]{chen2024ms} has covariance function $\xi^\Vec$ as in~\eqref{e.block_maps}, and its enriched free energy at the diagonal path associated with $q$ is $\bar F_n^\Vec(t,Lq)$, with $\bar F_n^\Vec$ as in~\eqref{e.vector_free_energy}. Setting $n_N:=\lceil(N+M_N)/(a+b)\rceil$, we thus obtain from~\cite[Corollary~3.2]{chen2024ms} that
\begin{equation}\label{e.block_approximation}
    \lim_{N\to\infty}\Ll|\widetilde F_{N,M_N}(t,q)-\frac1{a+b}\bar F_{n_N}^\Vec(t,Lq)\Rr|=0.
\end{equation}
Since $n_N$ tends to infinity with $N$, we deduce from~\eqref{e.block_approximation} and~\eqref{e.vector_upper_bip} that
\begin{equation*}
    \limsup_{N\to\infty}\widetilde F_{N,M_N}(t,q)\leq\frac{b}{a+b}\,\sP_{t,q}(p)=\frac{1}{1+\lambda}\,\sP_{t,q}(p).
\end{equation*}
By~\eqref{e.total_size_normalization}, we have
\begin{equation*}
    \bar F_N\Ll(\frac{(1+\lambda)N}{N+M_N}\,t,q\Rr)=\frac{N+M_N}{N}\,\widetilde F_{N,M_N}(t,q).
\end{equation*}
The ratio $(N+M_N)/N$ converges to $1+\lambda$, and by~\eqref{e.F_lip}, the right-hand side of this identity remains bounded, while the difference between its left-hand side and $\bar F_N(t,q)$ tends to zero as $N$ tends to infinity. We thus obtain~\eqref{e.critical_value}.

\smallskip
\noindent\emph{Step 2: Continuity in the number of spins in the second layer.}
We show that, for all integers $N,M,M'\geq1$ and every $q\in\mcl Q_\infty^2$,
\begin{equation}\label{e.rectangular_continuity}
    \Ll|\bar F_{N,M}(t,q)-\bar F_{N,M'}(t,q)\Rr|\leq\frac{|M-M'|}{N}\Ll(t+|q_2|_{L^1}\Rr).
\end{equation}
We assume without loss of generality that $M\leq M'$, and we use the Gaussian interpolation of~\cite[Corollary~2.2]{chen2024ms}, as in the proofs of~\cite[Lemmas~2.3 and~2.4]{chen2024ms}. We interpolate between the Hamiltonian of the model with $M'$ spins in the second layer, and the Hamiltonian obtained from it by replacing $v_j$ by zero for every $j\in\{M+1,\dots,M'\}$. Since $\nu$ is a probability measure, the free energies associated with these two Hamiltonians are $\bar F_{N,M'}(t,q)$ and $\bar F_{N,M}(t,q)$ respectively. By~\cite[Corollary~2.2]{chen2024ms}, the difference $N\bar F_{N,M'}(t,q)-N\bar F_{N,M}(t,q)$ is therefore the integral along the interpolation of the averaged Gibbs expectation of
\begin{equation*}
    \Ll(\frac tN\,\sigma\cdot\sigma'+q_2(\alpha\wedge\alpha')\Rr)\sum_{j=M+1}^{M'}v_jv_j',
\end{equation*}
where $(\sigma,v,\alpha)$ and $(\sigma',v',\alpha')$ denote two independent replicas. Since $|\sigma_i|\leq1$ and $|v_j|\leq1$, the absolute value of this quantity is at most $(M'-M)(t+q_2(\alpha\wedge\alpha'))$. Moreover, by the invariance property of the cascade, the overlap $\alpha\wedge\alpha'$ is uniformly distributed over $[0,1]$ under each of the averaged Gibbs measures along the interpolation; see~\cite[Proposition~4.8]{chen2025free}. The averaged Gibbs expectation of $q_2(\alpha\wedge\alpha')$ is thus equal to $|q_2|_{L^1}$, and we obtain~\eqref{e.rectangular_continuity}.

\smallskip
\noindent\emph{Step 3: The general case.}
Let $q\in\mcl Q_1^2$. Let $(\lambda_k)_{k\in\N}$ be a sequence of positive rational numbers that converges to $\lambda$, and for every $k,N\in\N$, let $M_N^{(k)}:=\max\{1,\lfloor\lambda_kN\rfloor\}$, so that $M_N^{(k)}/N$ converges to $\lambda_k$ as $N$ tends to infinity. For every $k\in\N$, we also let $q^{(k)}\in\mcl Q_\infty^2$ be the pair of truncated paths defined by $q_i^{(k)}(s):=\min\{q_i(s),k\}$, for $i\in\{1,2\}$ and $s\in[0,1)$, so that $q^{(k)}$ converges to $q$ in $L^1$ as $k$ tends to infinity. By Step~1, applied with $(M_N^{(k)})_{N\in\N}$, $\lambda_k$ and $q^{(k)}$ in place of $(M_N)_{N\in\N}$, $\lambda$ and $q$, there exists $p^{(k)}\in\mcl Q_{\infty,\leq1}\times\mcl Q_{\infty,\leq\lambda_k}$ such that
\begin{equation}\label{e.critical_rel_k}
    p_1^{(k)}=\partial_q\psi_{\mu_{\mcl I}}\big(q_1^{(k)}+tp_2^{(k)}\big)
    \qquad\text{and}\qquad
    p_2^{(k)}=\lambda_k\,\partial_q\psi_\nu\big(q_2^{(k)}+tp_1^{(k)}\big),
\end{equation}
and
\begin{equation}\label{e.critical_value_k}
    \limsup_{N\to\infty}\bar F_{N,M_N^{(k)}}\big(t,q^{(k)}\big)\leq\psi_{\mu_{\mcl I}}\big(q_1^{(k)}+tp_2^{(k)}\big)+\lambda_k\psi_\nu\big(q_2^{(k)}+tp_1^{(k)}\big)-t\langle p_1^{(k)},p_2^{(k)}\rangle_{L^2}.
\end{equation}
By~\eqref{e.F_lip} and~\eqref{e.rectangular_continuity}, we have
\begin{equation*}
    \Ll|\bar F_N(t,q)-\bar F_{N,M_N^{(k)}}\big(t,q^{(k)}\big)\Rr|\leq\big|q_1-q_1^{(k)}\big|_{L^1}+\frac{M_N}{N}\big|q_2-q_2^{(k)}\big|_{L^1}+\frac{|M_N-M_N^{(k)}|}{N}\Ll(t+|q_2|_{L^1}\Rr).
\end{equation*}
Letting $N$ tend to infinity, we obtain that $\limsup_{N\to\infty}\bar F_N(t,q)$ is bounded from above by the right-hand side of~\eqref{e.critical_value_k}, up to an error that tends to zero as $k$ tends to infinity. By Lemma~\ref{l.compact_embed_paths}, up to the extraction of a subsequence, we may assume that $p^{(k)}$ converges in $L^1$ to some $p\in\mcl Q_\infty^2$. The paths $q_1^{(k)}+tp_2^{(k)}$ and $q_2^{(k)}+tp_1^{(k)}$ then converge in $L^1$ to $q_1+tp_2$ and $q_2+tp_1$ respectively, and the continuity estimate~\eqref{e.continuity.der.psi} allows us to pass to the limit in~\eqref{e.critical_rel_k}. Hence $p$ is a critical point associated with $(t,q)$, and by~\eqref{e.bound.der.psi}, it belongs to $\mcl Q_{\infty,\leq1}\times\mcl Q_{\infty,\leq\lambda}$. Moreover, since the paths $p_1^{(k)}$ and $p_2^{(k)}$ are uniformly bounded, their convergence in $L^1$ implies that $\langle p_1^{(k)},p_2^{(k)}\rangle_{L^2}$ converges to $\la p_1,p_2\ra_{L^2}$. Using also the Lipschitz estimate~\eqref{e.psi_lip}, we see that the right-hand side of~\eqref{e.critical_value_k} converges to $\sP_{t,q}(p)$ as $k$ tends to infinity, and the proof of~\eqref{e.critical_value} is complete.
\end{proof}

We can now complete the proof of Theorem~\ref{t.main}, by combining the lower bound of Proposition~\ref{p.HJ_bdd_ms} with the upper bound of Proposition~\ref{p.crit_pt_bdd_ms}, the convexity of $\psi_{\mu_{\mcl I}}$, and the variational formula of Theorem~\ref{t.partial_hopf_paths}.

\begin{proof}[Proof of Theorem~\ref{t.main}]
Let $f:\R_+\times\mcl Q_2^2\to\R$ be the unique Lipschitz viscosity solution to~\eqref{e.main.hj}, in the sense of Definition~\ref{d.vs}. By Theorem~\ref{t.partial_hopf_paths}, the function $f$ is given by the right-hand side of~\eqref{e.main.1}. It therefore suffices to show that, for every $(t,q)\in\R_+\times\mcl Q_2^2$,
\begin{equation}\label{e.final_goal}
    \lim_{N\to\infty}\bar F_N(t,q)=f(t,q).
\end{equation}
For $t=0$, this follows from~\eqref{e.psi=sumlambdapsi}, since $f(0,\cdot)=\psi$. We now fix $t>0$ and $q\in\mcl Q_2^2$. By Proposition~\ref{p.crit_pt_bdd_ms}, there exists a critical point $p\in\mcl Q_{\infty,\leq1}\times\mcl Q_{\infty,\leq\lambda}$ associated with $(t,q)$ such that~\eqref{e.critical_value} holds. We set $q_1':=q_1+tp_2$, which belongs to $\mcl Q_2$, and we recall from~\eqref{e.critical_rel} that $p_1=\partial_q\psi_{\mu_{\mcl I}}(q_1')$. The main step of the proof is to show the Fenchel identity
\begin{equation}\label{e.final_fenchel}
    \psi_{\mu_{\mcl I}}^*(p_1)=\la p_1,q'_1\ra_{L^2}-\psi_{\mu_{\mcl I}}(q'_1).
\end{equation}
Let $y\in\mcl Q_2$. For every $\theta\in(0,1]$, the path $q_1'+\theta(y-q_1')=(1-\theta)q_1'+\theta y$ belongs to $\mcl Q_2$, and the convexity of $\psi_{\mu_{\mcl I}}$ given by Proposition~\ref{p.convex.psi.circ} yields
\begin{equation*}
    \frac{\psi_{\mu_{\mcl I}}\big(q_1'+\theta(y-q_1')\big)-\psi_{\mu_{\mcl I}}(q_1')}{\theta}\leq\psi_{\mu_{\mcl I}}(y)-\psi_{\mu_{\mcl I}}(q_1').
\end{equation*}
Letting $\theta$ tend to zero, and using that $\psi_{\mu_{\mcl I}}$ is Gateaux differentiable at $q_1'$ with derivative $p_1$, we obtain
\begin{equation}\label{e.final_convex}
    \psi_{\mu_{\mcl I}}(y)\geq\psi_{\mu_{\mcl I}}(q'_1)+\la p_1,y-q'_1\ra_{L^2},\qquad\forall y\in\mcl Q_2.
\end{equation}
Since $\mcl Q_\infty\subseteq\mcl Q_2$, we deduce from~\eqref{e.final_convex} and the definition of $\psi_{\mu_{\mcl I}}^*$ in~\eqref{e.conjugate.psiI} that
\begin{equation*}
    \psi_{\mu_{\mcl I}}^*(p_1)\leq\la p_1,q'_1\ra_{L^2}-\psi_{\mu_{\mcl I}}(q'_1).
\end{equation*}
For the converse inequality, we cannot directly use $q_1'$ as a test path in~\eqref{e.conjugate.psiI}, since $q_1'$ may be unbounded. For every $j\in\N$, we thus consider the truncated path $y_j:=q'_1\wedge j\in\mcl Q_\infty$, for which~\eqref{e.conjugate.psiI} gives
\begin{equation*}
    \psi_{\mu_{\mcl I}}^*(p_1)\geq\la p_1,y_j\ra_{L^2}-\psi_{\mu_{\mcl I}}(y_j).
\end{equation*}
As $j$ tends to infinity, the path $y_j$ converges to $q_1'$ in $L^1$. Since $0\leq p_1\leq1$, and by~\eqref{e.psi_lip}, the right-hand side of the previous display converges to $\la p_1,q_1'\ra_{L^2}-\psi_{\mu_{\mcl I}}(q_1')$. This completes the proof of~\eqref{e.final_fenchel}.

Since $\la p_1,q_1'\ra_{L^2}=\la p_1,q_1\ra_{L^2}+t\la p_1,p_2\ra_{L^2}$, we can combine the definition of $\sP_{t,q}$ in~\eqref{e.sP_lambda,t,q} with the identity~\eqref{e.final_fenchel} to get
\begin{equation*}
    \sP_{t,q}(p)=\la p_1,q_1\ra_{L^2}+\lambda\psi_\nu(q_2+tp_1)-\psi_{\mu_{\mcl I}}^*(p_1)\leq f(t,q),
\end{equation*}
where the inequality follows from the fact that $p_1$ belongs to $\mcl Q_\infty$ and that $f(t,q)$ is the supremum on the right-hand side of~\eqref{e.main.1}. Using also Proposition~\ref{p.HJ_bdd_ms}, we conclude that
\begin{equation*}
    f(t,q)\stackrel{\eqref{e.HJ_lower}}{\leq}\liminf_{N\to\infty}\bar F_N(t,q)\leq\limsup_{N\to\infty}\bar F_N(t,q)
    \stackrel{\eqref{e.critical_value}}{\leq}\sP_{t,q}(p)\leq f(t,q),
\end{equation*}
and thus that~\eqref{e.final_goal} holds.
\end{proof}

\begin{remark}[Measures with general compact support]\label{r.bounded_support}
Theorem~\ref{t.main} remains valid if the support of $\nu$ is only assumed to be contained in $[-a,a]$, for some $a>0$. To see this, we write $\bar F_N^\nu$ in place of $\bar F_N$ to display the dependence of the free energy on the measure $\nu$, and we denote by $\widehat\nu$ the image of $\nu$ under the mapping $v\mapsto v/a$, which is a probability measure supported on $[-1,1]$. The change of variables $v\mapsto v/a$ in~\eqref{e.F_N(t,q)=} and in~\eqref{e.def.psi.rho} gives
\begin{equation*}
    \bar F_N^\nu(t,q_1,q_2)=\bar F_N^{\widehat\nu}(a^2t,q_1,a^2q_2)
    \qquad\text{and}\qquad
    \psi_\nu(q_2)=\psi_{\widehat\nu}(a^2q_2).
\end{equation*}
Applying Theorem~\ref{t.main} to the measure $\widehat\nu$, we thus obtain that the formula~\eqref{e.main.1} holds for the measure $\nu$ as well. Moreover, in the passage from $\widehat\nu$ to $\nu$, the time derivative of the limit and its derivative with respect to the second path are both multiplied by $a^2$, so the equation in~\eqref{e.main.hj} is unchanged. The only difference is that the derivative of the limit with respect to the second path is now bounded by $\lambda a^2$ instead of $\lambda$; accordingly, in Definition~\ref{d.vs}, the regularization $\bar\xi$ must be chosen to agree with $\xi$ on a ball of radius at least $2\vee(1+\lambda a^2)$.
\end{remark}

\section{A shrinking domain of paths}
\label{s.moving_domain}

In Section~\ref{s.general_perceptron}, we will enrich general perceptron models in such a way that the variance of the pre-activations remains equal to one. This forces the time variable $t$ and the path $q_2$ associated with the pre-activations to satisfy the constraint $q_2\leq1-t$; see~\eqref{e.def.S_j}. The enriched free energy is then expected to satisfy, in the limit of large system size, the same equation as in~\eqref{e.main.hj}, but on a domain that shrinks as $t$ increases, and which is therefore not of the product form $(0,\infty)\times\mcl Q_2^2$ considered in Section~\ref{s.paths}. The goal of this section is to show that this equation satisfies a comparison principle on this domain, with no condition imposed on its moving boundary, and that its solution is still given by the variational formula of Theorem~\ref{t.partial_hopf_paths}; see Proposition~\ref{p.moving_comparison} and Theorem~\ref{t.partial_hopf_moving} below.

To see why no boundary condition is needed on the boundary of the set $\{q_2 \le 1-t\}$, let us denote by $f$ a solution to this equation, let $(t,q) = (t,q_1,q_2) \in (0,\infty) \times \mcl Q_2^2$ be such that $q_2 \le 1-t$, and suppose that $f$ is differentiable at $(t,q)$, with $p=(p_1,p_2) := \dr_q f(t,q)$. The characteristic line from $(t,q)$ that goes backwards in time is given by $s \mapsto (t-s, q_1 + s p_2, q_2 + s p_1)$, for $s \in [0,t]$; see \cite[(3.94)]{HJbook}. In order for this characteristic line not to exit the domain, it suffices that $p_1 \le 1$, for in this case we indeed have $q_2 + s p_1 \le 1-(t-s)$ for all $s \in [0,t]$. In other words, as long as $\dr_{q_1} f \le 1$, characteristic lines may leave the domain through this boundary, but they cannot enter it. In our setting, the bound $\dr_{q_1} f \le 1$ indeed holds, since $\dr_{q_1} f$ represents an average overlap between replicas of the spins, which cannot take values greater than $1$; see \eqref{e.bound.der.psi} or \eqref{e.general_path_derivatives}. 

In order to turn this observation into a rigorous statement, we need a notion of solution that encodes it. As in Definition~\ref{d.vs}, we require the viscosity inequalities to hold at every point of the domain with $t\in(0,1)$, including the points at which $q_2(1)=1-t$. Since the gradient of a test function at a boundary point is not constrained by the bounds satisfied by the derivatives of the solution, this requires extending the nonlinearity beyond the set in which these derivatives take their values. The extension used in Definition~\ref{d.vs} was designed for the cone $\mcl Q_2^2$, and it does not satisfy the inequality that we need at the moving boundary, namely~\eqref{e.moving_normal_compatibility} below. We therefore use a different extension, built so that the outflow property just described also holds at the level of test functions; see Lemma~\ref{l.moving_outflow}. With this notion of solution, the comparison principle follows from the usual doubling of variables, and the main task is to show that the variational formula satisfies the viscosity inequalities up to the boundary.

\subsection{The domain and the notion of solution}
\label{s.moving_definition}

Throughout this section, we write $\cH:=L^2([0,1))$, and for every $c\geq0$, we set $K_c:=\mcl Q_{\infty,\leq c}$. Each of the sets $K_c$ is convex, and it is compact in $\cH$ by Lemma~\ref{l.compact_embed_paths}. We denote by $\omega:\cH\to\R$ the support function of $K_1$, that is,
\begin{equation}\label{e.moving_support}
    \omega(x):=\sup_{p\in K_1}\la p,x\ra_{L^2},\qquad x\in\cH.
\end{equation}
As a supremum of linear functions, the function $\omega$ is positively homogeneous and subadditive. Since $0\in K_1$, and since every $p\in K_1$ takes values in $[0,1]$, we also have
\begin{equation}\label{e.moving_support_bound}
    0\leq\omega(x)\leq|x|_{L^1}\leq|x|_{L^2},\qquad\forall x\in\cH.
\end{equation}
In particular, the function $\omega$ is Lipschitz continuous. Since $K_c=cK_1$, the support function of $K_c$ is $c\,\omega$. We define
\begin{equation}\label{e.moving_domains}
    \widehat{\mathcal D}:=\Ll\{(t,q_1,q_2):t\in[0,1],\ q_1\in\mcl Q_2,\ q_2\in K_{1-t}\Rr\}
    \qquad\text{and}\qquad
    \mathcal D:=\widehat{\mathcal D}\cap(\R\times\mcl Q_\infty^2).
\end{equation}
The condition $q_2\in K_{1-t}$ means that $q_2\in\mcl Q_\infty$ and $q_2\leq1-t$. The set $\widehat{\mathcal D}$ is convex, and it is closed in the Hilbert space $\R\times\cH^2$, since $\mcl Q_2$ is closed in $\cH$ and the constraint $q_2\leq1-t$ is preserved by convergence in $\R\times\cH$. We write $\la\cdot,\cdot\ra$ and $|\cdot|$ for the scalar product and the norm of this space, and we use the same notation for the space $\cH^2$. We often write $z=(t,q)$ for a generic point of $\widehat{\mathcal D}$, with $q=(q_1,q_2)$. Unless otherwise specified, the Lipschitz continuity of a function defined on a subset of $\R\times\cH^2$ is understood with respect to the norm $|\cdot|$. The set $\mathcal D$ is the one on which the enriched free energies of Section~\ref{s.general_perceptron} are defined. Truncating $q_1$ shows that $\mathcal D$ is dense in $\widehat{\mathcal D}$, so every Lipschitz function on $\mathcal D$ admits a unique Lipschitz extension to $\widehat{\mathcal D}$.

We now fix $B\geq0$. In what follows, the sets $K_1$ and $K_B$ are those in which the derivatives of the solution with respect to $q_1$ and to $q_2$ take their values, respectively. We define $\H_B:\cH^2\to\R$ by
\begin{equation}\label{e.moving_H}
    \H_B(a,b):=\min_{p\in K_1}\max_{r\in K_B}\Ll\{\la a,r\ra_{L^2}+\la p,b\ra_{L^2}-\la p,r\ra_{L^2}\Rr\}
    =\min_{p\in K_1}\Ll\{\la p,b\ra_{L^2}+B\,\omega(a-p)\Rr\}.
\end{equation}
The two expressions coincide because $B\,\omega$ is the support function of $K_B$, and the extrema are attained because $K_1$ and $K_B$ are compact and $\omega$ is continuous. The next lemma records the properties of $\H_B$ that we will use: this function extends the bilinear nonlinearity in~\eqref{e.main.hj} from $K_1\times K_B$ to $\cH^2$, and its increments are bounded from above in terms of the support functions of $K_B$ and $K_1$.

\begin{lemma}[Properties of $\H_B$]\label{l.moving_H}
For every $a\in K_1$ and $b\in K_B$, we have
\begin{equation}\label{e.moving_H_agrees}
    \H_B(a,b)=\la a,b\ra_{L^2}.
\end{equation}
For every $a,b,a',b'\in\cH$, we have
\begin{equation}\label{e.moving_H_increment}
    \H_B(a+a',b+b')-\H_B(a,b)\leq B\,\omega(a')+\omega(b'),
\end{equation}
and in particular,
\begin{equation}\label{e.moving_H_lip}
    |\H_B(a,b)-\H_B(a',b')|\leq B|a-a'|_{L^2}+|b-b'|_{L^2}.
\end{equation}
\end{lemma}

\begin{proof}
Let $a\in K_1$ and $b\in K_B$. Choosing $p=a$ in the second expression in~\eqref{e.moving_H} gives $\H_B(a,b)\leq\la a,b\ra_{L^2}$. Conversely, since $b\in K_B$ and $B\,\omega$ is the support function of $K_B$, we have $B\,\omega(a-p)\geq\la b,a-p\ra_{L^2}$ for every $p\in K_1$, and thus $\H_B(a,b)\geq\la a,b\ra_{L^2}$. This proves~\eqref{e.moving_H_agrees}. Let now $a,b,a',b'\in\cH$. For every $p\in K_1$, the subadditivity of $\omega$ and the definition~\eqref{e.moving_support} give
\begin{equation*}
    \la p,b+b'\ra_{L^2}+B\,\omega(a+a'-p)\leq\la p,b\ra_{L^2}+B\,\omega(a-p)+\omega(b')+B\,\omega(a'),
\end{equation*}
and taking the minimum over $p\in K_1$ yields~\eqref{e.moving_H_increment}. The estimate~\eqref{e.moving_H_lip} follows from~\eqref{e.moving_H_increment}, applied in both directions, and from~\eqref{e.moving_support_bound}.
\end{proof}

\begin{definition}[Viscosity solutions on the shrinking domain]\label{d.moving_solution}
Let $B\geq0$, and let
\begin{equation}\label{e.moving_V}
    V:=\Ll\{(t,q_1,q_2)\in\widehat{\mathcal D}:0<t<1\Rr\}.
\end{equation}
A continuous function $u:\widehat{\mathcal D}\to\R$ is a \textbf{viscosity subsolution}, respectively a \textbf{viscosity supersolution}, of
\begin{equation}\label{e.moving_bilinear}
    \partial_tu-\la\partial_{q_1}u,\partial_{q_2}u\ra_{L^2}=0
\end{equation}
\textbf{with derivative bound $B$} if its restriction to $V$ is a viscosity subsolution, respectively a viscosity supersolution, of
\begin{equation}\label{e.moving_regularized}
    \partial_tu-\H_B(\partial_{q_1}u,\partial_{q_2}u)=0
\end{equation}
in the sense of Definition~\ref{d.vs_general}, with ambient Hilbert space $\cH^2$. A \textbf{viscosity solution} of~\eqref{e.moving_bilinear} with derivative bound $B$ is a function that is both a viscosity subsolution and a viscosity supersolution of~\eqref{e.moving_bilinear} with derivative bound $B$. For $u_0:\mcl Q_2\times K_1\to\R$, we say that such a function $u$ has initial condition $u_0$ if $u(0,\cdot)=u_0$ on $\mcl Q_2\times K_1$.
\end{definition}

Let us comment on this definition. The local extrema in there are understood relative to~$V$: for instance, $u-\phi$ has a local maximum at $z\in V$ if $(u-\phi)(z')\leq(u-\phi)(z)$ for every $z'\in V$ sufficiently close to~$z$. Definition~\ref{d.vs_general} involves test functions that are smooth on $V$, and the differential of such a function is required to be unique. This requirement is automatically satisfied here. Indeed, if $\zeta\in\R\times\cH^2$ is the difference between two candidate differentials at a point $z\in V$, then $\la\zeta,z'-z\ra=o(|z'-z|)$ as $z'\in V$ tends to $z$. Since $V$ is convex, we can apply this along the segment joining $z$ to any other point of $V$, and we obtain that $\zeta$ is orthogonal to $V-z$, and thus also to $V-V$. The linear span of $V-V$ contains the vectors $(1,0,0)$, $(0,\one_{[s,s')},0)$ and $(0,0,\one_{[s,s')})$ for every $0\leq s<s'\leq1$, so it is dense in $\R\times\cH^2$, and therefore $\zeta=0$.

For a smooth function $\phi:V\to\R$ and $z\in V$, we write
\begin{equation*}
    D\phi(z):=\Ll(\partial_t\phi(z),\partial_{q_1}\phi(z),\partial_{q_2}\phi(z)\Rr)\in\R\times\cH^2
\end{equation*}
for the differential of $\phi$ at $z$, and we set
\begin{equation}\label{e.moving_F}
    \mathscr F_B(\tau,a,b):=\tau-\H_B(a,b),\qquad(\tau,a,b)\in\R\times\cH^2.
\end{equation}
With this notation, a continuous function $u:\widehat{\mathcal D}\to\R$ is a viscosity subsolution of~\eqref{e.moving_bilinear} with derivative bound $B$ if and only if $\mathscr F_B(D\phi(z))\leq0$ whenever $\phi:V\to\R$ is smooth and $u-\phi$ has a local maximum at $z\in V$; for supersolutions, the condition is that $\mathscr F_B(D\phi(z))\geq0$ whenever $u-\phi$ has a local minimum at $z\in V$.

The key property of this notion of solution is a compatibility between the function $\mathscr F_B$ and the geometry of the domain, which we now state. For every $z\in\widehat{\mathcal D}$, we denote the normal cone to $\widehat{\mathcal D}$ at $z$ by
\begin{equation}\label{e.moving_normal_cone}
    N_{\widehat{\mathcal D}}(z):=\Ll\{n\in\R\times\cH^2:\la n,z'-z\ra\leq0\ \text{ for every }z'\in\widehat{\mathcal D}\Rr\}.
\end{equation}

\begin{lemma}[Outflow property]\label{l.moving_outflow}
Let $B\geq0$ and $z=(t,q_1,q_2)\in V$. For every $p\in K_1$, $r\in K_B$, and $h\in[0,t]$, we have
\begin{equation}\label{e.moving_backward}
    z+h(-1,r,p)\in\widehat{\mathcal D}.
\end{equation}
As a consequence, every $n=(n_t,n_1,n_2)\in N_{\widehat{\mathcal D}}(z)$ satisfies
\begin{equation}\label{e.moving_normal}
    n_t\geq B\,\omega(n_1)+\omega(n_2),
\end{equation}
and for every $\zeta\in\R\times\cH^2$ and $n\in N_{\widehat{\mathcal D}}(z)$, we have
\begin{equation}\label{e.moving_normal_compatibility}
    \mathscr F_B(\zeta+n)\geq\mathscr F_B(\zeta).
\end{equation}
\end{lemma}

\begin{proof}
The point in~\eqref{e.moving_backward} is $(t-h,q_1+hr,q_2+hp)$. We have $t-h\in[0,1]$ and $q_1+hr\in\mcl Q_2$, while the path $q_2+hp$ belongs to $\mcl Q_\infty$ and satisfies $q_2+hp\leq(1-t)+h=1-(t-h)$. This proves~\eqref{e.moving_backward}. Let now $n=(n_t,n_1,n_2)\in N_{\widehat{\mathcal D}}(z)$. Using the point in~\eqref{e.moving_backward} with $h=t$ as the point $z'$ in~\eqref{e.moving_normal_cone}, and dividing by $t>0$, we obtain that $-n_t+\la n_1,r\ra_{L^2}+\la n_2,p\ra_{L^2}\leq0$. Taking the supremum over $r\in K_B$ and $p\in K_1$ yields~\eqref{e.moving_normal}. Finally, for every $\zeta=(\tau,a,b)\in\R\times\cH^2$, the inequalities~\eqref{e.moving_H_increment} and~\eqref{e.moving_normal} give
\begin{equation*}
    \mathscr F_B(\zeta+n)-\mathscr F_B(\zeta)=n_t-\Ll(\H_B(a+n_1,b+n_2)-\H_B(a,b)\Rr)\geq n_t-B\,\omega(n_1)-\omega(n_2)\geq0,
\end{equation*}
which is~\eqref{e.moving_normal_compatibility}.
\end{proof}

The vector $(-1,r,p)$ in~\eqref{e.moving_backward} is the direction of the characteristic line of the equation~\eqref{e.moving_bilinear}, followed backward in time, when the derivatives of the solution with respect to $q_1$ and $q_2$ are $p$ and $r$ respectively. The property~\eqref{e.moving_backward} is thus the precise form of the observation made at the beginning of this section, and the inequality~\eqref{e.moving_normal_compatibility} is its counterpart for test functions. We will use~\eqref{e.moving_normal_compatibility} in the proof of Theorem~\ref{t.partial_hopf_moving}, to verify both viscosity inequalities at boundary points.

\subsection{Comparison and stability}
\label{s.moving_comparison}

Since the viscosity inequalities of Definition~\ref{d.moving_solution} are imposed at every point of $V$, the comparison principle can be proved by doubling the variables, as for equations posed on the whole space, and the boundary requires no specific treatment. We follow the proof of~\cite[Proposition~3.8]{chen2022hamilton}, with one modification that we point out in Remark~\ref{r.moving_comparison}. To compensate for the lack of local compactness of $\cH$, we use Stegall's variational principle~\cite{stegall1978optimization}, in the form stated in~\cite[Theorem~3.7]{chen2022hamilton}: if $\mathcal E$ is a nonempty, closed, bounded and convex subset of a separable Hilbert space $\mathcal X$, and if $J:\mathcal E\to\R$ is upper semicontinuous and bounded from above, then for every $\rho>0$, there exists $\ell\in\mathcal X$ with $|\ell|_{\mathcal X}\leq\rho$ such that the function $J+\la\ell,\cdot\ra_{\mathcal X}$ attains its maximum over $\mathcal E$. The statement in~\cite{chen2022hamilton} is for convex and weakly compact sets; a closed, bounded and convex subset of a Hilbert space is weakly compact.

\begin{proposition}[Comparison principle]\label{p.moving_comparison}
Let $B\geq0$, and let $u,v:\widehat{\mathcal D}\to\R$ be Lipschitz continuous functions such that $u$ is a viscosity subsolution and $v$ is a viscosity supersolution of~\eqref{e.moving_bilinear} with derivative bound $B$, in the sense of Definition~\ref{d.moving_solution}. If $u(0,q)\leq v(0,q)$ for every $q\in\mcl Q_2\times K_1$, then $u\leq v$ on $\widehat{\mathcal D}$.
\end{proposition}

\begin{proof}
Let $L<\infty$ be a Lipschitz constant for both $u$ and $v$. We argue by contradiction, and assume that $u>v$ at some point of $\widehat{\mathcal D}$. Every point of $\widehat{\mathcal D}$ of the form $(1,q_1,0)$ is the limit of the points $(t,q_1,0)$ as $t$ increases to $1$. Since $u$ and $v$ are continuous, and since $u\leq v$ at time $t=0$, there exists $(t_*,q_*)\in V$ such that $u(t_*,q_*)>v(t_*,q_*)$. We fix such a point, and then $\eta,\kappa,c_0\in(0,1]$ such that
\begin{equation}\label{e.moving_gap}
    u(t_*,q_*)-v(t_*,q_*)-\eta t_*-\frac{2\kappa}{1-t_*}\geq3c_0.
\end{equation}
For $\delta,\eps\in(0,1]$ to be chosen below, and for every $(t,q),(s,y)\in\widehat{\mathcal D}$ with $t,s<1$, we set
\begin{equation}\label{e.moving_J}
    J(t,q,s,y):=u(t,q)-v(s,y)-\frac{|t-s|^2+|q-y|^2}{2\eps}-\delta\Ll(|q|^2+|y|^2\Rr)-\eta t-\frac\kappa{1-t}-\frac\kappa{1-s}.
\end{equation}
We also consider linear perturbations of $J$. For $\rho\in(0,\delta]$ to be chosen below, and for $(\alpha,a,\beta,b)\in(\R\times\cH^2)^2$ such that
\begin{equation}\label{e.moving_ell_small}
    |\alpha|,\,|a|,\,|\beta|,\,|b|\leq\rho,
\end{equation}
we set $\ell(t,q,s,y):=\alpha t+\la a,q\ra+\beta s+\la b,y\ra$. Throughout the proof, we impose that
\begin{equation}\label{e.moving_delta_small}
    2\delta|q_*|^2\leq c_0\qquad\text{and}\qquad2\rho(1+|q_*|)\leq c_0,
\end{equation}
so that, by~\eqref{e.moving_gap},
\begin{equation}\label{e.moving_J_star}
    J(t_*,q_*,t_*,q_*)\geq2c_0\qquad\text{and}\qquad(J+\ell)(t_*,q_*,t_*,q_*)\geq c_0.
\end{equation}

\smallskip
\noindent\emph{Step 1: A priori bounds.}
We set $C_L:=L^2+2L+3$, and we show that, for every $(t,q),(s,y)\in\widehat{\mathcal D}$ with $t,s<1$ and every $\ell$ as above, if $(J+\ell)(t,q,s,y)\geq0$, then
\begin{equation}\label{e.moving_doubled_bounds}
    \frac{|t-s|^2+|q-y|^2}{4\eps}+\frac\delta2\Ll(|q|^2+|y|^2\Rr)+\frac\kappa{1-t}+\frac\kappa{1-s}\leq C_L.
\end{equation}
Since $q_2\leq1-t\leq1$, the point $(0,q)$ belongs to $\widehat{\mathcal D}$. Using that $u(0,q)\leq v(0,q)$, and the Lipschitz continuity of $u$ and $v$, we get
\begin{equation}\label{e.moving_initial_ordering}
    u(t,q)-v(s,y)\leq\Ll(u(t,q)-u(0,q)\Rr)+\Ll(v(0,q)-v(s,y)\Rr)\leq L\Ll(t+s+|q-y|\Rr).
\end{equation}
Moreover, we have $L|q-y|\leq\frac{1}{4\eps}|q-y|^2+L^2\eps$, and since $\rho\leq\delta\leq1$, the bound~\eqref{e.moving_ell_small} yields
\begin{equation}\label{e.moving_ell_bound}
    |\ell(t,q,s,y)|\leq\rho\Ll(2+|q|+|y|\Rr)\leq3+\frac\delta2\Ll(|q|^2+|y|^2\Rr).
\end{equation}
Combining these three estimates with the definition of $J$, and using that $t,s,\eps\leq1$, we obtain
\begin{equation*}
    (J+\ell)(t,q,s,y)\leq C_L-\frac{|t-s|^2+|q-y|^2}{4\eps}-\frac\delta2\Ll(|q|^2+|y|^2\Rr)-\frac\kappa{1-t}-\frac\kappa{1-s},
\end{equation*}
and the announced implication follows.

\smallskip
\noindent\emph{Step 2: Construction of a maximizer.}
We set $R:=2\sqrt{C_L/\delta}$ and $\theta:=\kappa/(2C_L)$, and we let
\begin{equation*}
    \mathcal E:=\Ll\{(t,q,s,y):(t,q),(s,y)\in\widehat{\mathcal D},\ t\vee s\leq1-\theta,\ |q|\vee|y|\leq R\Rr\}.
\end{equation*}
This set is closed, bounded and convex in $(\R\times\cH^2)^2$. Since $u$ and $v$ are Lipschitz continuous, the function $J$ is continuous and bounded from above on $\mathcal E$. By~\eqref{e.moving_J_star} and Step~1 with $\ell=0$, the point $(t_*,q_*,t_*,q_*)$ satisfies~\eqref{e.moving_doubled_bounds}, and it therefore belongs to $\mathcal E$. By Stegall's variational principle, there exists $(\alpha,a,\beta,b)$ satisfying~\eqref{e.moving_ell_small} such that, with $\ell$ defined as above, the function $J+\ell$ attains its maximum over $\mathcal E$ at some point $(\bar t,\bar q,\bar s,\bar y)\in\mathcal E$. By~\eqref{e.moving_J_star}, we have
\begin{equation}\label{e.moving_max_positive}
    (J+\ell)(\bar t,\bar q,\bar s,\bar y)\geq c_0.
\end{equation}
Step~1 therefore applies at this point, and~\eqref{e.moving_doubled_bounds} yields
\begin{equation}\label{e.moving_max_bounds}
    |\bar t-\bar s|+|\bar q-\bar y|\leq2\sqrt{2C_L\eps},\qquad|\bar q|+|\bar y|\leq2\sqrt{C_L/\delta},\qquad\text{and}\qquad\bar t\vee\bar s\leq1-\frac\kappa{C_L}.
\end{equation}
More precisely,~\eqref{e.moving_doubled_bounds} gives $|\bar q|^2+|\bar y|^2\leq2C_L/\delta$, so that $|\bar q|\vee|\bar y|<R$, and we also have $\bar t\vee\bar s<1-\theta$. That is, both inequalities introduced in the definition of $\mathcal E$ are strict at the maximizer.

\smallskip
\noindent\emph{Step 3: Positivity of the contact times.}
We show that $\bar t>0$ and $\bar s>0$, provided that
\begin{equation}\label{e.moving_eps_small}
    2L\sqrt{2C_L\eps}\leq\frac{c_0}3\qquad\text{and}\qquad2\rho\Ll(1+\sqrt{C_L/\delta}\Rr)\leq\frac{c_0}3.
\end{equation}
If $\bar t=0$, then the ordering at time zero and the Lipschitz continuity of $v$ give
\begin{equation*}
    u(0,\bar q)-v(\bar s,\bar y)\leq v(0,\bar q)-v(\bar s,\bar y)\leq L\Ll(\bar s+|\bar q-\bar y|\Rr).
\end{equation*}
If $\bar s=0$, then we similarly have $u(\bar t,\bar q)-v(0,\bar y)\leq u(\bar t,\bar q)-u(0,\bar y)\leq L(\bar t+|\bar q-\bar y|)$. In both cases, the definition of $J$, the first inequality in~\eqref{e.moving_ell_bound}, and the bounds~\eqref{e.moving_max_bounds} and~\eqref{e.moving_eps_small} yield
\begin{equation*}
    (J+\ell)(\bar t,\bar q,\bar s,\bar y)\leq L\Ll(|\bar t-\bar s|+|\bar q-\bar y|\Rr)+\rho\Ll(2+|\bar q|+|\bar y|\Rr)\leq2L\sqrt{2C_L\eps}+2\rho\Ll(1+\sqrt{C_L/\delta}\Rr)\leq\frac{2c_0}3,
\end{equation*}
in contradiction with~\eqref{e.moving_max_positive}.

\smallskip
\noindent\emph{Step 4: Conclusion.}
We assume that~\eqref{e.moving_eps_small} holds. By Steps~2 and~3, the points $(\bar t,\bar q)$ and $(\bar s,\bar y)$ belong to $V$, and for every $(t,q)\in V$ sufficiently close to $(\bar t,\bar q)$, the point $(t,q,\bar s,\bar y)$ belongs to $\mathcal E$. The function $(t,q)\mapsto(J+\ell)(t,q,\bar s,\bar y)$ therefore has a local maximum at $(\bar t,\bar q)$ relative to $V$. Up to an additive constant, this function is $u-\phi_1$, where
\begin{equation*}
    \phi_1(t,q):=\frac{|t-\bar s|^2+|q-\bar y|^2}{2\eps}+\delta|q|^2+\eta t+\frac\kappa{1-t}-\alpha t-\la a,q\ra
\end{equation*}
is a smooth function on $V$. Since $u$ is a viscosity subsolution, we obtain that
\begin{equation}\label{e.moving_sub_contact}
    \frac{\bar t-\bar s}\eps+\eta+\frac\kappa{(1-\bar t)^2}-\alpha-\H_B(P)\leq0,
    \qquad\text{where}\quad P:=\frac{\bar q-\bar y}\eps+2\delta\bar q-a.
\end{equation}
Here and below, for $P=(P_1,P_2)\in\cH^2$, we write $\H_B(P):=\H_B(P_1,P_2)$. Similarly, the function $(s,y)\mapsto-(J+\ell)(\bar t,\bar q,s,y)$ has a local minimum at $(\bar s,\bar y)$ relative to $V$, and up to an additive constant, it is equal to $v-\phi_2$, where
\begin{equation*}
    \phi_2(s,y):=-\frac{|\bar t-s|^2+|\bar q-y|^2}{2\eps}-\delta|y|^2-\frac\kappa{1-s}+\beta s+\la b,y\ra
\end{equation*}
is a smooth function on $V$. Since $v$ is a viscosity supersolution, we obtain that
\begin{equation}\label{e.moving_super_contact}
    \frac{\bar t-\bar s}\eps-\frac\kappa{(1-\bar s)^2}+\beta-\H_B(Q)\geq0,
    \qquad\text{where}\quad Q:=\frac{\bar q-\bar y}\eps-2\delta\bar y+b.
\end{equation}
Subtracting~\eqref{e.moving_super_contact} from~\eqref{e.moving_sub_contact}, and using~\eqref{e.moving_H_lip},~\eqref{e.moving_ell_small}, and~\eqref{e.moving_max_bounds}, we get
\begin{align}
    \eta&\leq\H_B(P)-\H_B(Q)+\alpha+\beta\notag\\
    &\leq(B+1)\Ll(2\delta\Ll(|\bar q|+|\bar y|\Rr)+2\rho\Rr)+2\rho\leq4(B+1)\sqrt{C_L\delta}+2(B+2)\rho.\label{e.moving_final_contradiction}
\end{align}
We now choose the parameters. We first choose $\delta\in(0,1]$ such that $2\delta|q_*|^2\leq c_0$ and $4(B+1)\sqrt{C_L\delta}\leq\eta/3$. We then choose $\rho\in(0,\delta]$ such that $2\rho(1+|q_*|)\leq c_0$, $2\rho(1+\sqrt{C_L/\delta})\leq c_0/3$, and $2(B+2)\rho\leq\eta/3$. Finally, we choose $\eps\in(0,1]$ such that $2L\sqrt{2C_L\eps}\leq c_0/3$. With these choices, the conditions~\eqref{e.moving_delta_small} and~\eqref{e.moving_eps_small} are satisfied, and~\eqref{e.moving_final_contradiction} gives $\eta\leq2\eta/3$, which is absurd since $\eta>0$.
\end{proof}

\begin{remark}\label{r.moving_comparison}
In the proof of~\cite[Proposition~3.8]{chen2022hamilton}, the distance between the two contact points is estimated by comparing the value of the doubled function at the maximizer $(\bar t,\bar q,\bar s,\bar y)$ with its values at $(\bar t,\bar q,\bar s,\bar q)$ and at $(\bar t,\bar q,\bar t,\bar y)$. On the shrinking domain, this is not possible in general, since the points $(\bar s,\bar q)$ and $(\bar t,\bar y)$ need not belong to $\widehat{\mathcal D}$. In the proof above, we instead rely on the estimate~\eqref{e.moving_initial_ordering}, in which the only auxiliary point is $(0,q)$; this point always belongs to $\widehat{\mathcal D}$, because the domain shrinks as $t$ increases. The resulting bound on the distance between the contact points in~\eqref{e.moving_max_bounds} is of order $\sqrt\eps$ instead of $\eps$, but this is still sufficient to conclude.
\end{remark}

\begin{proposition}[Stability]\label{p.moving_stability}
Let $B\geq0$, and let $(u_j)_{j\in\N}$ be a sequence of viscosity subsolutions, respectively supersolutions, of~\eqref{e.moving_bilinear} with derivative bound $B$. If $u_j$ converges to a function $u:\widehat{\mathcal D}\to\R$ uniformly on every bounded subset of $\widehat{\mathcal D}$ as $j$ tends to infinity, then $u$ is a viscosity subsolution, respectively supersolution, of~\eqref{e.moving_bilinear} with derivative bound $B$.
\end{proposition}

\begin{proof}
The function $u$ is continuous, since it is a locally uniform limit of continuous functions. We only consider the case of subsolutions, the case of supersolutions being similar. Let $\phi:V\to\R$ be a smooth function and $z=(t,q)\in V$ be such that $u-\phi$ has a local maximum at $z$. We set $\widetilde\phi(z'):=\phi(z')+|z'-z|^2$, which defines a smooth function on $V$ such that $D\widetilde\phi(z)=D\phi(z)$. We fix $\delta\in(0,t\wedge(1-t))$ sufficiently small that, denoting by $\mathcal E$ the intersection of $\widehat{\mathcal D}$ with the closed ball of radius $\delta$ centered at $z$, we have
\begin{equation}\label{e.moving_strict_max}
    (u-\widetilde\phi)(z')\leq(u-\widetilde\phi)(z)-|z'-z|^2,\qquad\forall z'\in\mathcal E.
\end{equation}
The set $\mathcal E$ is closed, bounded and convex, and it is contained in $V$. We set $\eps_j:=\sup_{\mathcal E}|u_j-u|$, which tends to zero as $j$ tends to infinity. For $j$ sufficiently large, the continuous function $u_j-\widetilde\phi$ is bounded from above on $\mathcal E$, by~\eqref{e.moving_strict_max}, and Stegall's variational principle provides us with $\ell_j\in\R\times\cH^2$ such that $|\ell_j|\leq2^{-j}$ and such that the function $u_j-\widetilde\phi+\la\ell_j,\cdot\ra$ attains its maximum over $\mathcal E$ at some point $z_j\in\mathcal E$. Comparing the values of this function at $z_j$ and at $z$, and using~\eqref{e.moving_strict_max}, we get
\begin{equation*}
    0\leq(u_j-\widetilde\phi)(z_j)-(u_j-\widetilde\phi)(z)+\la\ell_j,z_j-z\ra\leq2\eps_j-|z_j-z|^2+2^{-j}\delta,
\end{equation*}
so that $z_j$ converges to $z$ as $j$ tends to infinity. In particular, for $j$ sufficiently large, the point $z_j$ lies in the open ball of radius $\delta$ centered at $z$, and the function $u_j-(\widetilde\phi-\la\ell_j,\cdot\ra)$ therefore has a local maximum at $z_j$ relative to $V$. Since $u_j$ is a viscosity subsolution, we obtain that
\begin{equation*}
    \mathscr F_B\Ll(D\widetilde\phi(z_j)-\ell_j\Rr)\leq0.
\end{equation*}
Letting $j$ tend to infinity, and using that $\mathscr F_B$ and $D\widetilde\phi$ are continuous and that $D\widetilde\phi(z)=D\phi(z)$, we conclude that $\mathscr F_B(D\phi(z))\leq0$, as desired.
\end{proof}

\subsection{The variational formula}
\label{s.moving_formula}

We now show that the variational formula of Theorems~\ref{t.hopf.half} and~\ref{t.partial_hopf_paths} remains valid on the shrinking domain. We state the result in the formulation of Remark~\ref{r.general_g}, in which the convex part of the initial condition is described through its convex conjugate. That is, we give ourselves a convex function $g$ on $\cH$ whose effective domain $\operatorname{dom}g:=\{p\in\cH:g(p)<+\infty\}$ is contained in $K_1$; we also give ourselves a function $\psi_2$ on $K_1$; and we consider the initial condition $(q_1,q_2)\mapsto g^*(q_1)+\psi_2(q_2)$. The function $\psi_2$ is not assumed to be convex. Our only assumption on $\psi_2$ is the one-sided bound~\eqref{e.moving_pattern_slope} below, which is a way to state that the derivative of $\psi_2$ takes values in $K_B$ without assuming that $\psi_2$ is differentiable. Indeed, if $\psi_2$ is the restriction to $K_1$ of a continuously differentiable function on $\cH$ whose gradient takes values in $K_B$, then~\eqref{e.moving_pattern_slope} follows by integrating along the segment from $x$ to $y$, since $B\,\omega$ is the support function of $K_B$; conversely, we will see in the proof of Theorem~\ref{t.partial_hopf_moving} that every function satisfying~\eqref{e.moving_pattern_slope} is a uniform limit of functions of this form. By~\eqref{e.moving_support_bound}, the condition~\eqref{e.moving_pattern_slope} implies that $\psi_2$ is Lipschitz continuous with respect to the $L^1$ norm, with Lipschitz constant $B$. It also implies that $\psi_2(x)\leq\psi_2(y)$ whenever $y-x\in\mcl Q_2^*$, since we have $\omega(x-y)=0$ in this case; compare with the properties of the functions $\psi_\rho$ stated in Lemma~\ref{l.psi^vec_smooth}.

The assumption $\operatorname{dom}g\subseteq K_1$ is the counterpart, for the first path, of the assumption on~$\psi_2$ discussed below: since $g^*$ is a supremum of affine functions whose slopes lie in $K_1$, we have $g^*(y)-g^*(x)\leq\omega(y-x)$ for every $x,y\in\cH$, which is a way to state that the derivative of $g^*$ takes values in $K_1$. This bound is essential: as explained at the beginning of this section, it ensures that the characteristic lines move the variable $q_2$ at speed at most one. In contrast with the bound $B$ on the derivative with respect to $q_2$, the bound $1$ on the derivative with respect to~$q_1$ cannot be relaxed, since it must match the speed at which the domain shrinks.

\begin{theorem}[Variational formula on the shrinking domain]\label{t.partial_hopf_moving}
Let $B\geq0$, let $g:\cH\to\R\cup\{+\infty\}$ be a proper, lower semicontinuous and convex function such that $\operatorname{dom}g\subseteq K_1$, and let $\psi_2:K_1\to\R$ be such that
\begin{equation}\label{e.moving_pattern_slope}
    \psi_2(y)-\psi_2(x)\leq B\,\omega(y-x),\qquad\forall x,y\in K_1.
\end{equation}
For every $(t,q_1,q_2)\in\widehat{\mathcal D}$, let
\begin{equation}\label{e.partial_hopf_moving}
    u(t,q_1,q_2):=\sup_{p\in K_1}\Ll\{\la p,q_1\ra_{L^2}+\psi_2(q_2+tp)-g(p)\Rr\}.
\end{equation}
The supremum in~\eqref{e.partial_hopf_moving} is attained. The function $u$ is the unique Lipschitz viscosity solution of~\eqref{e.moving_bilinear} with derivative bound $B$, in the sense of Definition~\ref{d.moving_solution}, such that
\begin{equation}\label{e.moving_initial}
    u(0,q_1,q_2)=g^*(q_1)+\psi_2(q_2),\qquad\forall(q_1,q_2)\in\mcl Q_2\times K_1,
\end{equation}
where $g^*(q_1):=\sup_{p\in\cH}\{\la p,q_1\ra_{L^2}-g(p)\}$. Moreover, for every $(t,q),(s,y)\in\widehat{\mathcal D}$, we have
\begin{equation}\label{e.moving_formula_lip}
    |u(t,q)-u(s,y)|\leq|q_1-y_1|_{L^1}+B|q_2-y_2|_{L^1}+B|t-s|.
\end{equation}
\end{theorem}

The right-hand side of~\eqref{e.partial_hopf_moving} does not involve $B$. Theorem~\ref{t.partial_hopf_moving} thus shows in particular that the viscosity solution with derivative bound $B$ and initial condition~\eqref{e.moving_initial} does not depend on the choice of $B$, as long as~\eqref{e.moving_pattern_slope} holds.

\begin{proof}[Proof of Theorem~\ref{t.partial_hopf_moving}]
For every $z=(t,q_1,q_2)\in\widehat{\mathcal D}$ and $p\in K_1$, the path $q_2+tp$ belongs to $\mcl Q_\infty$ and is bounded by $(1-t)+t=1$. It thus belongs to $K_1$, and the right-hand side of~\eqref{e.partial_hopf_moving} is well-defined. For every $p\in\operatorname{dom}g$ and $z=(t,q_1,q_2)\in\widehat{\mathcal D}$, we set
\begin{equation}\label{e.moving_branch}
    w_p(z):=\la p,q_1\ra_{L^2}+\psi_2(q_2+tp)-g(p),
\end{equation}
so that $u=\sup_{p\in\operatorname{dom}g}w_p$. The structure of the proof is the same as that of the proof of Theorem~\ref{t.hopf.half}, with the following differences: the space is infinite-dimensional, which we handle using the compactness of $K_1$; the contact points can be on the boundary of the domain, which we handle using Lemma~\ref{l.moving_outflow}; and the approximation of $\psi_2$ by smooth functions must preserve the condition~\eqref{e.moving_pattern_slope}.

\smallskip
\noindent\emph{Step 1: Attainment, Lipschitz continuity, and the initial condition.}
By~\eqref{e.moving_pattern_slope} and~\eqref{e.moving_support_bound}, we have
\begin{equation}\label{e.moving_Phi_lip}
    |\psi_2(y)-\psi_2(x)|\leq B|y-x|_{L^1},\qquad\forall x,y\in K_1.
\end{equation}
We fix $z=(t,q_1,q_2)\in\widehat{\mathcal D}$. By~\eqref{e.moving_Phi_lip}, the function $p\mapsto\la p,q_1\ra_{L^2}+\psi_2(q_2+tp)$ is continuous on $K_1$ for the topology of $\cH$. Since $g$ is lower semicontinuous and is not identically $+\infty$ on $K_1$, and since $K_1$ is compact, we deduce that $u(z)$ is finite, that the supremum in~\eqref{e.partial_hopf_moving} is attained, and that the set of maximizers
\begin{equation*}
    \mcl A(z):=\Ll\{p\in\operatorname{dom}g:w_p(z)=u(z)\Rr\}
\end{equation*}
is nonempty and compact. For every $p\in K_1$ and $(t,q),(s,y)\in\widehat{\mathcal D}$, we have $|\la p,q_1-y_1\ra_{L^2}|\leq|q_1-y_1|_{L^1}$, and by~\eqref{e.moving_Phi_lip},
\begin{equation*}
    |\psi_2(q_2+tp)-\psi_2(y_2+sp)|\leq B|q_2-y_2|_{L^1}+B|t-s|.
\end{equation*}
Each of the functions $w_p$ therefore satisfies the estimate~\eqref{e.moving_formula_lip}, and taking the supremum over $p$ yields~\eqref{e.moving_formula_lip} for $u$. In particular, the function $u$ is Lipschitz continuous. Finally, since $\operatorname{dom}g\subseteq K_1$, the identity~\eqref{e.moving_initial} follows from~\eqref{e.partial_hopf_moving} at $t=0$.

\smallskip
\noindent\emph{Step 2: The subsolution inequality when $\psi_2$ is smooth.}
In this step and the next one, we assume in addition that $\psi_2$ is the restriction to $K_1$ of a continuously differentiable function on $\cH$, which we still denote by $\psi_2$, and that $\nabla\psi_2(x)\in K_B$ for every $x\in\cH$. For every $p\in\operatorname{dom}g$, the formula~\eqref{e.moving_branch} then defines a continuously differentiable function $w_p$ on $\R\times\cH^2$, whose differential at $z'=(t',q_1',q_2')$ is
\begin{equation}\label{e.moving_branch_gradient}
    Dw_p(z')=\Ll(\la p,\nabla\psi_2(q_2'+t'p)\ra_{L^2},\,p,\,\nabla\psi_2(q_2'+t'p)\Rr)\in\R\times\cH^2,
\end{equation}
and the right-hand side of~\eqref{e.moving_branch_gradient} depends continuously on $(p,z')\in K_1\times\R\times\cH^2$. We fix $z=(t,q_1,q_2)\in V$, and for every $p\in\mcl A(z)$, we set
\begin{equation*}
    b_p:=\nabla\psi_2(q_2+tp)\in K_B\qquad\text{and}\qquad\gamma_p:=Dw_p(z)=\Ll(\la p,b_p\ra_{L^2},p,b_p\Rr).
\end{equation*}
Since $p\in K_1$ and $b_p\in K_B$, the identity~\eqref{e.moving_H_agrees} gives
\begin{equation}\label{e.moving_branch_solution}
    \mathscr F_B(\gamma_p)=\la p,b_p\ra_{L^2}-\H_B(p,b_p)=0,\qquad\forall p\in\mcl A(z).
\end{equation}
Let now $\phi:V\to\R$ be a smooth function such that $u-\phi$ has a local maximum at $z$, and let $p\in\mcl A(z)$. Since $w_p\leq u$ and $w_p(z)=u(z)$, the function $w_p-\phi$ also has a local maximum at $z$ relative to $V$. Let $z'\in\widehat{\mathcal D}$. Since $\widehat{\mathcal D}$ is convex and $t\in(0,1)$, the point $z+\eps(z'-z)$ belongs to $V$ for every $\eps\in[0,1)$, and thus $(w_p-\phi)(z+\eps(z'-z))\leq(w_p-\phi)(z)$ for every $\eps>0$ sufficiently small. Dividing by $\eps$ and letting $\eps$ tend to zero, we obtain that $\la\gamma_p-D\phi(z),z'-z\ra\leq0$. Since $z'\in\widehat{\mathcal D}$ is arbitrary, this means that $\gamma_p-D\phi(z)\in N_{\widehat{\mathcal D}}(z)$. By~\eqref{e.moving_normal_compatibility} and~\eqref{e.moving_branch_solution}, we conclude that
\begin{equation*}
    \mathscr F_B(D\phi(z))\leq\mathscr F_B(\gamma_p)=0,
\end{equation*}
which is the subsolution inequality.

\smallskip
\noindent\emph{Step 3: The supersolution inequality when $\psi_2$ is smooth.}
We keep the additional assumption and the notation of Step~2. We first show that
\begin{equation}\label{e.moving_active_monotonicity}
    \la p-p',b_p-b_{p'}\ra_{L^2}\geq0,\qquad\forall p,p'\in\mcl A(z).
\end{equation}
Let $p\in\mcl A(z)$ and $p'\in\operatorname{dom}g$. For every $\eps\in(0,1]$, the point $p_\eps:=p+\eps(p'-p)$ belongs to $\operatorname{dom}g$, by convexity, and $w_{p_\eps}(z)\leq u(z)=w_p(z)$. Together with the convexity inequality $g(p_\eps)\leq g(p)+\eps(g(p')-g(p))$, this gives
\begin{equation*}
    \eps\la p'-p,q_1\ra_{L^2}+\psi_2\big(q_2+tp+\eps t(p'-p)\big)-\psi_2(q_2+tp)\leq g(p_\eps)-g(p)\leq\eps\big(g(p')-g(p)\big).
\end{equation*}
Dividing by $\eps$ and letting $\eps$ tend to zero, we obtain the analogue of~\eqref{e.partial_subgradient}, namely
\begin{equation}\label{e.moving_subgradient}
    g(p')-g(p)\geq\la q_1+tb_p,p'-p\ra_{L^2}.
\end{equation}
For $p,p'\in\mcl A(z)$, we add~\eqref{e.moving_subgradient} to the same inequality with $p$ and $p'$ interchanged. Since $t>0$, this yields~\eqref{e.moving_active_monotonicity}.

We set $\Gamma_z:=\{\gamma_p:p\in\mcl A(z)\}$. This set is compact, since it is the image of the compact set $\mcl A(z)$ under a continuous mapping. Its closed convex hull $\overline{\operatorname{co}}\,\Gamma_z$ is also compact. Indeed, it is closed, and it is totally bounded: for every $\eps>0$, we can cover $\Gamma_z$ by finitely many balls of radius $\eps$, and then every point of the convex hull of $\Gamma_z$ is within distance $\eps$ of the convex hull of the centers of these balls, which is compact. We next show that
\begin{equation}\label{e.moving_active_covariance}
    \mathscr F_B(\bar\gamma)\geq0,\qquad\forall\bar\gamma\in\overline{\operatorname{co}}\,\Gamma_z.
\end{equation}
Since $\mathscr F_B$ is continuous, it suffices to consider the case when $\bar\gamma=\sum_{i=1}^m\alpha_i\gamma_{p_i}$ for some $m\in\N$, $p_1,\dots,p_m\in\mcl A(z)$, and $\alpha_1,\dots,\alpha_m\geq0$ with $\sum_{i=1}^m\alpha_i=1$. Since $K_1$ and $K_B$ are convex, we have $\sum_{i=1}^m\alpha_ip_i\in K_1$ and $\sum_{i=1}^m\alpha_ib_{p_i}\in K_B$, and the identity~\eqref{e.moving_H_agrees} thus gives
\begin{align*}
    \mathscr F_B(\bar\gamma)
    =\sum_{i=1}^m\alpha_i\la p_i,b_{p_i}\ra_{L^2}-\Big\langle\sum_{i=1}^m\alpha_ip_i,\sum_{j=1}^m\alpha_jb_{p_j}\Big\rangle_{L^2}=\frac12\sum_{i,j=1}^m\alpha_i\alpha_j\la p_i-p_j,b_{p_i}-b_{p_j}\ra_{L^2}
    \stackrel{\eqref{e.moving_active_monotonicity}}{\geq}0.
\end{align*}

Let now $\phi:V\to\R$ be a smooth function such that $u-\phi$ has a local minimum at $z$. We will show that
\begin{equation}\label{e.moving_lower_jet}
    D\phi(z)\in\overline{\operatorname{co}}\,\Gamma_z+N_{\widehat{\mathcal D}}(z).
\end{equation}
This implies the supersolution inequality. Indeed, if $D\phi(z)=\bar\gamma+n$ with $\bar\gamma\in\overline{\operatorname{co}}\,\Gamma_z$ and $n\in N_{\widehat{\mathcal D}}(z)$, then~\eqref{e.moving_normal_compatibility} and~\eqref{e.moving_active_covariance} give
\begin{equation*}
    \mathscr F_B(D\phi(z))\geq\mathscr F_B(\bar\gamma)\geq0.
\end{equation*}
As a first step towards~\eqref{e.moving_lower_jet}, we show that
\begin{equation}\label{e.moving_directional}
    \la D\phi(z),z'-z\ra\leq\max_{\gamma\in\Gamma_z}\la\gamma,z'-z\ra,\qquad\forall z'\in\widehat{\mathcal D}.
\end{equation}
We fix $z'\in\widehat{\mathcal D}$ and write $v:=z'-z$. For every $h\in(0,1)$, the point $z+hv$ belongs to $V$, and we pick $p_h\in\mcl A(z+hv)$. Since $u(z+hv)=w_{p_h}(z+hv)$ and $u(z)\geq w_{p_h}(z)$, the local minimum of $u-\phi$ at $z$ and the mean value theorem give, for every $h>0$ sufficiently small,
\begin{equation}\label{e.moving_mvt}
    \phi(z+hv)-\phi(z)\leq u(z+hv)-u(z)\leq w_{p_h}(z+hv)-w_{p_h}(z)=h\la Dw_{p_h}(z+\theta_hhv),v\ra
\end{equation}
for some $\theta_h\in(0,1)$. Since $K_1$ is compact, there is a sequence $h_k\downarrow0$ such that $p_{h_k}$ converges in $\cH$ to some $p\in K_1$. We verify that $p\in\mcl A(z)$. Writing $z+h_kv=(t_k,q_{1,k},q_{2,k})$, the continuity of $\psi_2$ on $K_1$ and of $u$ gives
\begin{equation*}
    g(p_{h_k})=\la p_{h_k},q_{1,k}\ra_{L^2}+\psi_2(q_{2,k}+t_kp_{h_k})-u(z+h_kv)\xrightarrow[k\to\infty]{}\la p,q_1\ra_{L^2}+\psi_2(q_2+tp)-u(z),
\end{equation*}
and the lower semicontinuity of $g$ yields $g(p)\leq\la p,q_1\ra_{L^2}+\psi_2(q_2+tp)-u(z)$. That is, we have $p\in\operatorname{dom}g$ and $w_p(z)\geq u(z)$, and thus $p\in\mcl A(z)$. By the continuity of the right-hand side of~\eqref{e.moving_branch_gradient}, the vector $Dw_{p_{h_k}}(z+\theta_{h_k}h_kv)$ converges to $\gamma_p$ as $k$ tends to infinity. Dividing~\eqref{e.moving_mvt} by $h_k$ and letting $k$ tend to infinity, we obtain that $\la D\phi(z),v\ra\leq\la\gamma_p,v\ra$, which proves~\eqref{e.moving_directional}.

We now derive~\eqref{e.moving_lower_jet} from~\eqref{e.moving_directional} by a separation argument. The set $\overline{\operatorname{co}}\,\Gamma_z+N_{\widehat{\mathcal D}}(z)$ is convex, and it is closed, as the sum of a compact set and a closed set. If $D\phi(z)$ did not belong to this set, then there would exist $v\in\R\times\cH^2$ such that
\begin{equation*}
    \la D\phi(z),v\ra>\sup\Ll\{\la\gamma+n,v\ra:\gamma\in\overline{\operatorname{co}}\,\Gamma_z,\ n\in N_{\widehat{\mathcal D}}(z)\Rr\}.
\end{equation*}
Since $N_{\widehat{\mathcal D}}(z)$ is a cone, this would imply that
\begin{equation}\label{e.moving_separation}
    \la n,v\ra\leq0\quad\text{for every }n\in N_{\widehat{\mathcal D}}(z),\qquad\text{and}\qquad\la D\phi(z),v\ra>\max_{\gamma\in\Gamma_z}\la\gamma,v\ra.
\end{equation}
Let $\mathcal T_z:=\{h(z'-z):z'\in\widehat{\mathcal D},\ h\geq0\}$, which is a convex cone. By~\eqref{e.moving_directional} and by homogeneity, we have $\la D\phi(z),w\ra\leq\max_{\gamma\in\Gamma_z}\la\gamma,w\ra$ for every $w\in\mathcal T_z$. Both sides of this inequality depend continuously on $w$, so the inequality remains valid for every $w$ in the closure of $\mathcal T_z$. By the second part of~\eqref{e.moving_separation}, the vector $v$ thus does not belong to the closure of $\mathcal T_z$. Since this closure is a closed convex cone, there would then exist $n\in\R\times\cH^2$ such that $\la n,v\ra>0$ and $\la n,w\ra\leq0$ for every $w\in\mathcal T_z$. The latter property means that $n\in N_{\widehat{\mathcal D}}(z)$, and we have reached a contradiction with the first part of~\eqref{e.moving_separation}. This completes the proof of~\eqref{e.moving_lower_jet}, and thus of the supersolution inequality.

\smallskip
\noindent\emph{Step 4: Conclusion.}
We now remove the additional assumption on $\psi_2$ made in Steps~2 and~3. Let $\psi_2:K_1\to\R$ be a function satisfying~\eqref{e.moving_pattern_slope}. We first extend it to $\cH$, by setting
\begin{equation}\label{e.moving_pattern_extension}
    \widetilde\psi_2(x):=\inf_{y\in K_1}\Ll\{\psi_2(y)+B\,\omega(x-y)\Rr\},\qquad x\in\cH.
\end{equation}
We fix $y_0\in K_1$. For every $x\in\cH$ and $y\in K_1$, the assumption~\eqref{e.moving_pattern_slope} and the subadditivity of $\omega$ give
\begin{equation*}
    \psi_2(y)+B\,\omega(x-y)\geq\psi_2(y_0)-B\,\omega(y_0-y)+B\,\omega(x-y)\geq\psi_2(y_0)-B\,\omega(y_0-x),
\end{equation*}
so that $\widetilde\psi_2(x)$ is finite. Using~\eqref{e.moving_pattern_slope} once more, we see that the infimum in~\eqref{e.moving_pattern_extension} is attained at $y=x$ when $x\in K_1$, and thus that $\widetilde\psi_2=\psi_2$ on $K_1$. The subadditivity of $\omega$ also implies that
\begin{equation}\label{e.moving_extension_slope}
    \widetilde\psi_2(x+x')-\widetilde\psi_2(x)\leq B\,\omega(x'),\qquad\forall x,x'\in\cH.
\end{equation}
In particular, by~\eqref{e.moving_support_bound}, the function $\widetilde\psi_2$ is Lipschitz continuous on $\cH$, with Lipschitz constant $B$.

We next regularize $\widetilde\psi_2$, using the mappings $\iota_j$ and $\pi_j$ defined in~\eqref{e.dyadic_maps}. We recall that $\la\iota_jy,x\ra_{L^2}=\la y,\pi_jx\ra_j$ for every $y\in E_j$ and $x\in\cH$, and we denote by $|\cdot|_j$ the norm associated with the scalar product $\la\cdot,\cdot\ra_j$, so that $|\iota_jy|_{L^2}=|y|_j$. For every $j\in\N$, let $\chi_j$ be a smooth probability density on $E_j$ supported in $\{y\in E_j:|y|_j\leq2^{-j}\}$, and let
\begin{equation}\label{e.moving_pattern_mollified}
    \psi_{2,j}(x):=\int_{E_j}\widetilde\psi_2\big(\iota_j(\pi_jx-y)\big)\chi_j(y)\,\d y,\qquad x\in\cH.
\end{equation}
The function $\psi_{2,j}$ is the composition of the continuous linear mapping $\pi_j$ with the convolution of the continuous function $\widetilde\psi_2\circ\iota_j$ and the smooth density $\chi_j$. It is therefore continuously differentiable on $\cH$. We now show that its gradient takes values in $K_B$. The mapping $\iota_j\pi_j$ is self-adjoint, and it maps $K_1$ into itself, since averaging a path over dyadic intervals preserves its monotonicity and the bounds that it satisfies. For every $x'\in\cH$, we thus have
\begin{equation*}
    \omega(\iota_j\pi_jx')=\sup_{p\in K_1}\la\iota_j\pi_jp,x'\ra_{L^2}\leq\omega(x').
\end{equation*}
Combining this with~\eqref{e.moving_extension_slope} and~\eqref{e.moving_pattern_mollified}, we obtain that
\begin{equation}\label{e.moving_mollified_slope}
    \psi_{2,j}(x+x')-\psi_{2,j}(x)\leq B\,\omega(x'),\qquad\forall x,x'\in\cH.
\end{equation}
Replacing $x'$ by $hx'$ in~\eqref{e.moving_mollified_slope}, with $h>0$, dividing by $h$, and letting $h$ tend to zero, we get
\begin{equation*}
    \la\nabla\psi_{2,j}(x),x'\ra_{L^2}\leq B\,\omega(x')=\sup_{r\in K_B}\la r,x'\ra_{L^2},\qquad\forall x,x'\in\cH.
\end{equation*}
Since $K_B$ is closed and convex, the separating hyperplane theorem then ensures that $\nabla\psi_{2,j}(x)\in K_B$ for every $x\in\cH$. The function $\psi_{2,j}$ thus satisfies the additional assumption of Steps~2 and~3, and by~\eqref{e.moving_mollified_slope}, its restriction to $K_1$ satisfies~\eqref{e.moving_pattern_slope}. Moreover, since $\widetilde\psi_2$ is $B$-Lipschitz and coincides with $\psi_2$ on $K_1$, we have
\begin{equation}\label{e.moving_pattern_convergence}
    \sup_{x\in K_1}|\psi_{2,j}(x)-\psi_2(x)|\leq B\sup_{x\in K_1}|\iota_j\pi_jx-x|_{L^2}+2^{-j}B\xrightarrow[j\to\infty]{}0.
\end{equation}
The convergence to zero follows from~\eqref{e.dyadic_convergence} and the compactness of $K_1$, since the mappings $\iota_j\pi_j$ are contractions of $\cH$.

For every $j\in\N$, let $u_j:\widehat{\mathcal D}\to\R$ be defined by~\eqref{e.partial_hopf_moving} with $\psi_{2,j}$ in place of $\psi_2$. By Steps~1 to~3, the function $u_j$ is a viscosity solution of~\eqref{e.moving_bilinear} with derivative bound $B$. Since $q_2+tp\in K_1$ for every $(t,q_1,q_2)\in\widehat{\mathcal D}$ and $p\in K_1$, comparing~\eqref{e.partial_hopf_moving} for $\psi_2$ and for $\psi_{2,j}$ gives
\begin{equation*}
    \sup_{\widehat{\mathcal D}}|u_j-u|\leq\sup_{K_1}|\psi_{2,j}-\psi_2|,
\end{equation*}
which tends to zero as $j$ tends to infinity, by~\eqref{e.moving_pattern_convergence}. By Proposition~\ref{p.moving_stability}, the function $u$ is therefore a viscosity solution of~\eqref{e.moving_bilinear} with derivative bound $B$. By Step~1, it is Lipschitz continuous and satisfies~\eqref{e.moving_initial} and~\eqref{e.moving_formula_lip}. Finally, the uniqueness of Lipschitz viscosity solutions with a given initial condition follows from Proposition~\ref{p.moving_comparison}.
\end{proof}

\subsection{Additional unconstrained variables}
\label{s.moving_scalar}

In Section~\ref{s.ext}, we will use a version of Theorem~\ref{t.partial_hopf_moving} in which each of the two groups of variables contains an additional real variable that is not subject to any constraint. We now state this result precisely and explain how to adapt the arguments of this section. The two groups of variables are $(q_1,x_1)$ and $(q_2,x_2)$, with $(t,q_1,q_2)\in\widehat{\mathcal D}$ and $x_1,x_2\in\R$. We set $E:=\cH\times\R$, and we denote by $\la\cdot,\cdot\ra_{L^2\times\R}$ the scalar product on this space, that is, $\la(a,\alpha),(b,\beta)\ra_{L^2\times\R}:=\la a,b\ra_{L^2}+\alpha\beta$. Writing $\partial_{(q_i,x_i)}u:=(\partial_{q_i}u,\partial_{x_i}u)$, we consider the equation
\begin{equation}\label{e.moving_scalar_bilinear}
    \partial_tu-\la\partial_{(q_1,x_1)}u,\partial_{(q_2,x_2)}u\ra_{L^2\times\R}=0.
\end{equation}
We write $(t,q,x)$ for a generic point of $\widehat{\mathcal D}^\ext:=\widehat{\mathcal D}\times\R^2$, with $q=(q_1,q_2)$ and $x=(x_1,x_2)$, and we identify $\widehat{\mathcal D}^\ext$ with a closed convex subset of the Hilbert space $\R\times E^2$ through the mapping $(t,q,x)\mapsto(t,(q_1,x_1),(q_2,x_2))$. We also set $V^\ext:=V\times\R^2$. We fix $B,R_1,R_2\geq0$, and we replace the sets $K_1$ and $K_B$ by the compact and convex subsets of $E$ defined by
\begin{equation}\label{e.moving_scalar_sets}
    K_1^\ext:=K_1\times[-R_1,R_1]\qquad\text{and}\qquad K_B^\ext:=K_B\times[-R_2,R_2].
\end{equation}
The support functions of these two sets are, for every $(y,\xi)\in E$,
\begin{equation}\label{e.moving_scalar_support}
    \omega_1^\ext(y,\xi):=\omega(y)+R_1|\xi|\qquad\text{and}\qquad\omega_B^\ext(y,\xi):=B\,\omega(y)+R_2|\xi|.
\end{equation}
We define $\H^\ext:E^2\to\R$ as in~\eqref{e.moving_H}, that is,
\begin{equation}\label{e.moving_scalar_H}
    \H^\ext(a,b):=\min_{p\in K_1^\ext}\max_{r\in K_B^\ext}\Ll\{\la a,r\ra_{L^2\times\R}+\la p,b\ra_{L^2\times\R}-\la p,r\ra_{L^2\times\R}\Rr\}.
\end{equation}
In analogy with Definition~\ref{d.moving_solution}, we say that a continuous function $u:\widehat{\mathcal D}^\ext\to\R$ is a viscosity subsolution, respectively supersolution, of~\eqref{e.moving_scalar_bilinear} with derivative bounds $(B,R_1,R_2)$ if its restriction to $V^\ext$ is a viscosity subsolution, respectively supersolution, of
\begin{equation*}
    \partial_tu-\H^\ext\Ll(\partial_{(q_1,x_1)}u,\partial_{(q_2,x_2)}u\Rr)=0
\end{equation*}
in the sense of Definition~\ref{d.vs_general}, with ambient Hilbert space $E^2$. A viscosity solution is a function that is both a viscosity subsolution and a viscosity supersolution.

\begin{proposition}[Variational formula with additional variables]\label{p.moving_scalar}
Let $B,R_1,R_2\geq0$, and let $g:E\to\R\cup\{+\infty\}$ be a proper, lower semicontinuous and convex function such that $\operatorname{dom}g\subseteq K_1^\ext$. Let $\psi_2:K_1\times\R\to\R$ be such that, for every $(y,\xi),(y',\xi')\in K_1\times\R$,
\begin{equation}\label{e.moving_scalar_slope}
    \psi_2(y',\xi')-\psi_2(y,\xi)\leq B\,\omega(y'-y)+R_2|\xi'-\xi|.
\end{equation}
For every $(t,q,x)\in\widehat{\mathcal D}^\ext$, let
\begin{equation}\label{e.partial_hopf_moving_scalar}
    u(t,q,x):=\sup_{(p,m)\in K_1^\ext}\Ll\{\la p,q_1\ra_{L^2}+mx_1+\psi_2(q_2+tp,x_2+tm)-g(p,m)\Rr\}.
\end{equation}
The supremum in~\eqref{e.partial_hopf_moving_scalar} is attained. The function $u$ is the unique Lipschitz viscosity solution of~\eqref{e.moving_scalar_bilinear} with derivative bounds $(B,R_1,R_2)$ such that, for every $(q,x)\in\mcl Q_2\times K_1\times\R^2$,
\begin{equation}\label{e.moving_scalar_initial}
    u(0,q,x)=g^*(q_1,x_1)+\psi_2(q_2,x_2),
\end{equation}
where $g^*(q_1,x_1):=\sup_{(p,m)\in E}\{\la p,q_1\ra_{L^2}+mx_1-g(p,m)\}$. Moreover, for every $(t,q,x)$ and $(t',q',x')$ in $\widehat{\mathcal D}^\ext$, we have
\begin{multline}\label{e.moving_scalar_lip}
    |u(t,q,x)-u(t',q',x')|\leq|q_1-q_1'|_{L^1}+R_1|x_1-x_1'|\\
    +B|q_2-q_2'|_{L^1}+R_2|x_2-x_2'|+(B+R_1R_2)|t-t'|.
\end{multline}
\end{proposition}

\begin{proof}
The proof is the same as that of Theorem~\ref{t.partial_hopf_moving}, including the preliminary results on which it relies, up to the following modifications. Throughout, the space $E$ replaces $\cH$, the sets $K_1^\ext$ and $K_B^\ext$ replace $K_1$ and $K_B$, and the functions $\omega_1^\ext$ and $\omega_B^\ext$ replace $\omega$ and $B\,\omega$.

The proof of Lemma~\ref{l.moving_H} only uses that $K_1$ and $K_B$ are compact and convex, and that $\omega$ and $B\,\omega$ are their support functions. The function $\H^\ext$ therefore coincides with $\la\cdot,\cdot\ra_{L^2\times\R}$ on $K_1^\ext\times K_B^\ext$, and for every $a,b,a',b'\in E$, it satisfies
\begin{equation*}
    \H^\ext(a+a',b+b')-\H^\ext(a,b)\leq\omega_B^\ext(a')+\omega_1^\ext(b').
\end{equation*}
Since the sets $K_1^\ext$ and $K_B^\ext$ are bounded, their support functions are bounded by a constant multiple of the norm of $E$, and the function $\H^\ext$ is thus Lipschitz continuous. Concerning Lemma~\ref{l.moving_outflow}, since the variables $x_1$ and $x_2$ are not constrained, for every $(t,q,x)\in V^\ext$, $(p,m)\in K_1^\ext$, $(r,m')\in K_B^\ext$, and $h\in[0,t]$, we have
\begin{equation*}
    \Ll(t-h,q_1+hr,q_2+hp,x_1+hm',x_2+hm\Rr)\in\widehat{\mathcal D}^\ext.
\end{equation*}
The proof of Lemma~\ref{l.moving_outflow} thus applies without change, and it shows that the analogue of~\eqref{e.moving_normal_compatibility} holds for the function $(\tau,a,b)\mapsto\tau-\H^\ext(a,b)$ and the normal cones to $\widehat{\mathcal D}^\ext$ at the points of $V^\ext$. The proofs of Propositions~\ref{p.moving_comparison} and~\ref{p.moving_stability} only use that the domain is closed and convex, that the time variable ranges in $[0,1]$, that the point obtained by replacing $t$ by $0$ in any given point of the domain also belongs to the domain, and that the nonlinearity is Lipschitz continuous. They therefore also apply in the present setting, provided that the norms of $\cH^2$ and $\R\times\cH^2$ are replaced throughout by those of $E^2$ and $\R\times E^2$, so that the quadratic terms in~\eqref{e.moving_J} also involve the variables $x$, and that the factor $B+1$ in~\eqref{e.moving_final_contradiction} is replaced by a Lipschitz constant of $\H^\ext$.

We now turn to the proof of Theorem~\ref{t.partial_hopf_moving} itself. In Step~1, the estimate~\eqref{e.moving_Phi_lip} is replaced by the fact that, by~\eqref{e.moving_scalar_slope} and~\eqref{e.moving_support_bound}, the function $\psi_2$ is Lipschitz continuous on $K_1\times\R$, with Lipschitz constants $B$ in the first variable, with respect to the $L^1$ norm, and $R_2$ in the second variable. Using also the compactness of $K_1^\ext$, we obtain that the supremum in~\eqref{e.partial_hopf_moving_scalar} is attained, and that~\eqref{e.moving_scalar_initial} and~\eqref{e.moving_scalar_lip} hold; the last term in~\eqref{e.moving_scalar_lip} comes from the fact that $|m|\leq R_1$ for every $(p,m)\in K_1^\ext$. Steps~2 and~3 are unchanged, under the additional assumption that $\psi_2$ is the restriction to $K_1\times\R$ of a continuously differentiable function on $E$ whose gradient takes values in $K_B^\ext$. In Step~4, we extend $\psi_2$ to $E$ by setting, for every $(y,\xi)\in E$,
\begin{equation*}
    \widetilde\psi_2(y,\xi):=\inf_{(y',\xi')\in K_1\times\R}\Ll\{\psi_2(y',\xi')+\omega_B^\ext(y-y',\xi-\xi')\Rr\}.
\end{equation*}
The argument given below~\eqref{e.moving_pattern_extension} shows that this function is finite, that it coincides with $\psi_2$ on $K_1\times\R$, and that its increments are bounded from above by $\omega_B^\ext$, as in~\eqref{e.moving_extension_slope}. We then set, for every $j\in\N$ and $(y,\xi)\in E$,
\begin{equation*}
    \psi_{2,j}(y,\xi):=\int_{E_j\times\R}\widetilde\psi_2\big(\iota_j(\pi_jy-y'),\xi-\xi'\big)\chi_j(y',\xi')\,\d y'\,\d\xi',
\end{equation*}
where $\chi_j$ is now a smooth probability density on $E_j\times\R$ supported in $\{(y',\xi'):|y'|_j^2+|\xi'|^2\leq4^{-j}\}$. As in Step~4, the function $\psi_{2,j}$ is continuously differentiable on $E$, its gradient takes values in $K_B^\ext$, and
\begin{equation*}
    \sup_{(y,\xi)\in K_1\times\R}|\psi_{2,j}(y,\xi)-\psi_2(y,\xi)|\leq B\sup_{y\in K_1}|\iota_j\pi_jy-y|_{L^2}+2^{-j}(B+R_2)\xrightarrow[j\to\infty]{}0.
\end{equation*}
The functions defined by~\eqref{e.partial_hopf_moving_scalar} with $\psi_{2,j}$ in place of $\psi_2$ therefore converge to $u$ uniformly on $\widehat{\mathcal D}^\ext$, and we conclude as before.
\end{proof}

\section{Towards general perceptron models}
\label{s.general_perceptron}

We now return to the general perceptron models discussed in Section~\ref{s.intro.perceptron}, and explain how the approach of this paper is expected to extend to them. We consider five models in turn: the Ising perceptron with a general activation function, as in \cite{openai2026ising}, in Section~\ref{s.gen.ising}; the spherical perceptron, as in \cite{aronow2026negative, openai2026spherical}, in Section~\ref{s.spherical}; Ising and spherical perceptrons with i.i.d.\ random activation functions, as in~\cite{fuSohn2026perceptrons}, in Section~\ref{s.random_activations}; an Ising perceptron in which the patterns carry labels and are sampled from a mixture of two Gaussian distributions, in Section~\ref{s.ext}; and the same model with spherical spins, in Section~\ref{s.ext.spherical}. In the last two models, the patterns are no longer standard Gaussian vectors: each pattern is the sum of a standard Gaussian vector and of a random shift, whose direction is the same for all the patterns, and whose sign is given by the label of the pattern.

For each of the first four models, the discussion follows the same four steps. First, we define an enriched free energy, which depends on a time variable $t$ and on a pair of paths $q=(q_1,q_2)$ ranging in the shrinking domain of Section~\ref{s.moving_domain}, and which reduces to the free energy of interest at $t=1$ and $q=0$. Second, a computation based on Gaussian integration by parts shows that the enriched free energy satisfies the Hamilton--Jacobi equation~\eqref{e.moving_bilinear}, up to an error term that is small as soon as the overlap between two replicas concentrates conditionally on the overlap of the cascade. Third, we compute the enriched free energy at $t=0$, and observe that it is the sum of a convex function of $q_1$, which only depends on the reference measure of the spins, and of a function of $q_2$, which only depends on the activation function, and which is not convex in general. Fourth, since this is exactly the structure required by Theorem~\ref{t.partial_hopf_moving}, this theorem identifies the unique viscosity solution to the equation with this initial condition, and represents it by a variational formula. This solution is our candidate for the limit of the enriched free energy, and its value at $t=1$ and $q=0$ is our candidate for the limit free energy of the model. For the fourth model, the initial condition does not have the required convexity property, and an additional step is needed, in which the spin configurations are first sorted according to their magnetization. The fifth model is treated differently: sorting the spin configurations according to their projection onto the direction of the shifts, and using the rotational invariance of the sphere, we reduce it to the spherical perceptron with random activation functions of Section~\ref{s.random_activations}.

The status of these steps is as follows. For the first three models, the third and fourth steps are rigorous: the results of Section~\ref{s.moving_domain} provide the comparison principle and the variational formula for the equation on the shrinking domain, and we verify below that the initial conditions satisfy their assumptions. On the other hand, we do not prove that the enriched free energy converges to the candidate limit. This would require to show analogues of Propositions~\ref{p.HJ_bdd_ms} and~\ref{p.crit_pt_bdd_ms}; in particular, the cavity arguments underlying the proof of Proposition~\ref{p.crit_pt_bdd_ms} would need to be adjusted. Under various assumptions on the activation functions, the limit free energies of the Ising and spherical perceptrons were identified in~\cite{openai2026ising, aronow2026negative, openai2026spherical}, while the models with random activation functions considered in Section~\ref{s.random_activations} were covered in \cite{fuSohn2026perceptrons}. Whenever the results of~\cite{aronow2026negative, fuSohn2026perceptrons, openai2026ising, openai2026spherical} apply, the candidate formulas that we obtain coincide with the formulas proved there. For the fifth model, the reduction to the model of Section~\ref{s.random_activations} is rigorous, so that the candidate formula for this model is a consequence of the one obtained there for the spherical perceptron with random activation functions.

Throughout this section, the parameter $\lambda$, the sequence $(M_N)_{N\in\N}$, the Gaussian random variables $(g_{ij})_{i,j\in\N}$ and the cascade $\fR$ are as in Section~\ref{s.intro.main}, and we use the notation of Section~\ref{s.moving_domain}: we write $\cH=L^2([0,1))$ and $K_c=\mcl Q_{\infty,\leq c}$, the function $\omega$ is the support function of $K_1$ defined in~\eqref{e.moving_support}, and the domains $\mathcal D\subseteq\widehat{\mathcal D}$ are defined in~\eqref{e.moving_domains}. In each of Sections~\ref{s.gen.ising} to~\ref{s.ext}, we denote by $\bar F_N$ the enriched free energy of the model under consideration, and by $f$ the candidate for its limit. The meaning of these two symbols thus changes from one subsection to the next, and it also differs from the one they had in Sections~\ref{s.intro} to~\ref{s.bounds}, where they referred to the bipartite model; the relation between the two settings is explained in Remark~\ref{r.gen.bipartite}.

\subsection{The Ising perceptron}
\label{s.gen.ising}

Let $U:\R\to\R$ be a twice continuously differentiable function with bounded first and second derivatives. We do not assume that $U$ is bounded or concave; for instance, for every $s\geq0$ and every probability measure $\nu$ with compact support, the activation function $u_{s,\nu}$ defined in~\eqref{e.u_t,nu=} satisfies our assumption. We consider the Ising perceptron with activation function $U$, that is, the model in~\eqref{e.def.perc.general} with $P_N=\mu_{\mcl I}^{\otimes N}$ and $M=M_N$, where we recall that $\mu_{\mcl I}=\frac12(\delta_{-1}+\delta_1)$. Its free energy is
\begin{equation}\label{e.general_free_energy}
    -\frac1N\E\log\int\exp\Ll(\sum_{j=1}^{M_N}U\Ll(\frac1{\sqrt N}\sum_{i=1}^Ng_{ij}\sigma_i\Rr)\Rr)\,\d\mu_{\mcl I}^{\otimes N}(\sigma).
\end{equation}

\smallskip
\noindent\textbf{The enriched free energy.}
We recall that the domain $\mathcal D$ defined in~\eqref{e.moving_domains} is the set of triples $(t,q_1,q_2)$ with $t\in[0,1]$, $q_1\in\mcl Q_\infty$, and $q_2\in\mcl Q_\infty$ such that $q_2\leq1-t$, and we write $q=(q_1,q_2)$. Conditionally on $\fR$, let $(V_i^{q_1})_{i\in\N}$ and $(W_j^{q_2})_{j\in\N}$ be independent centered Gaussian processes indexed by $\supp\fR$, independent of $(g_{ij})$, with covariances
\begin{equation*}
    \E\Ll[V_i^{q_1}(\alpha)V_i^{q_1}(\alpha')\mid\fR\Rr]=q_1(\alpha\wedge\alpha')\qquad\text{and}\qquad\E\Ll[W_j^{q_2}(\alpha)W_j^{q_2}(\alpha')\mid\fR\Rr]=q_2(\alpha\wedge\alpha').
\end{equation*}
Let also $(z_j)_{j\in\N}$ be independent standard Gaussian random variables, independent of all the other random variables. For every $(t,q)\in\mathcal D$, $\sigma\in\{-1,1\}^N$, $\alpha\in\supp\fR$ and $j\in\N$, we set
\begin{equation}\label{e.def.S_j}
    S_j:=\sqrt{\frac tN}\sum_{i=1}^Ng_{ij}\sigma_i+W_j^{q_2}(\alpha)+\sqrt{1-t-q_2(1)}\,z_j.
\end{equation}
The first term in~\eqref{e.def.S_j} is the pre-activation of the original model multiplied by $\sqrt t$, and the second term is a field indexed by the cascade, as in the bipartite model. The third term compensates for the variations of the first two: conditionally on $\fR$, and for every fixed $(\sigma,\alpha)$, the random variable $S_j$ is a centered Gaussian with variance $t+q_2(1)+(1-t-q_2(1))=1$, irrespective of the choice of $(t,q)\in\mathcal D$. The constraint $q_2\leq1-t$ in the definition of $\mathcal D$ is precisely what makes this possible. We define the enriched free energy $\bar F_N:\mathcal D\to\R$ by
\begin{equation}\label{e.general_enrichment}
    \bar F_N(t,q):=-\frac1N\E\log\iint\E_z\exp\Ll(\sum_{j=1}^{M_N}U(S_j)+\sqrt{2}\sum_{i=1}^NV_i^{q_1}(\alpha)\sigma_i-Nq_1(1)\Rr)\,\d\mu_{\mcl I}^{\otimes N}(\sigma)\,\d\fR(\alpha),
\end{equation}
where $\E_z$ denotes the expectation with respect to the variables $(z_j)_{j\in\N}$ only, and $\E$ denotes the expectation with respect to all the other sources of randomness. The variables $(z_j)$ are thus integrated inside the logarithm, on the same footing as $\sigma$ and $\alpha$. As in~\eqref{e.H^t,q_N=}, the term $Nq_1(1)$ is one half of the variance of the Gaussian field acting on the spins. At $t=1$, the constraint $q_2\leq1-t$ imposes that $q_2=0$, and $\bar F_N(1,0)$ is the free energy in~\eqref{e.general_free_energy}. At the other endpoint $t=0$, the quantities $S_j$ no longer depend on the spins, and the enriched free energy can be computed explicitly; see~\eqref{e.general_initial} below.

\smallskip
\noindent\textbf{The approximate Hamilton--Jacobi equation.}
We denote by $\la\cdot\ra$ the expectation with respect to the Gibbs measure associated with~\eqref{e.general_enrichment}. This is a random probability measure on the set of triples $(\sigma,\alpha,z)$, with $z=(z_j)_{j\leq M_N}$, and we write $(\sigma^1,\alpha^1,z^1)$ and $(\sigma^2,\alpha^2,z^2)$ for two independent random variables sampled according to it, also called replicas. For $\ell\in\{1,2\}$, we write $S_j^\ell$ for the quantity in~\eqref{e.def.S_j} evaluated at $(\sigma^\ell,\alpha^\ell,z^\ell)$, and we consider the three overlaps
\begin{equation*}
    R_{12}^{\alpha}:=\alpha^1\wedge\alpha^2,\qquad R_{12}^{\sigma}:=\frac1N\sigma^1\cdot\sigma^2,\qquad\text{and}\qquad R_{12}^{U'}:=\frac1N\sum_{j=1}^{M_N}U'(S_j^1)U'(S_j^2).
\end{equation*}
Under the measure $\E\la\cdot\ra$, the overlap $R_{12}^\alpha$ is uniformly distributed on $[0,1]$; see~\cite[Proposition~4.8]{chen2025free}. For every integrable random variable $X$, we write $\E\la X\mid R_{12}^\alpha\ra$ for the conditional expectation of $X$ given $R_{12}^\alpha$ under this measure, and $\E\la X\mid R_{12}^{\alpha}=s\ra$ for its value at $s\in[0,1]$, so that
\begin{equation*}
    \E\la X\ra=\int_0^1\E\la X\mid R_{12}^{\alpha}=s\ra\,\d s.
\end{equation*}
By Gaussian integration by parts, we have, for almost every $s\in[0,1]$,
\begin{equation}\label{e.general_path_derivatives}
    \partial_{q_1}\bar F_N(s)=\E\la R_{12}^{\sigma}\mid R_{12}^{\alpha}=s\ra\qquad\text{and}\qquad\partial_{q_2}\bar F_N(s)=\frac12\E\la R_{12}^{U'}\mid R_{12}^{\alpha}=s\ra,
\end{equation}
as well as
\begin{equation}\label{e.general_time_derivative}
    \partial_t\bar F_N=\frac12\E\la R_{12}^{U'}R_{12}^{\sigma}\ra.
\end{equation}
Here, as in~\eqref{e.main.hj}, the derivative $\partial_{q_i}\bar F_N$ at a given point $(t,q)$ is an element of $L^2([0,1))$, and we write $\partial_{q_i}\bar F_N(s)$ for its value at $s$, keeping the dependence on $(t,q)$ implicit. In the derivation of~\eqref{e.general_path_derivatives} and~\eqref{e.general_time_derivative}, the integrations by parts also produce terms that involve $U''(S_j)$ and $U'(S_j)^2$ evaluated at a single replica, and these terms cancel out. This cancellation comes from the fact that the variance of $S_j$ does not depend on $(t,q)$, and it is the reason for the introduction of the variables $(z_j)$ in~\eqref{e.def.S_j}. Since $R_{12}^{\alpha}$ is uniformly distributed, the identities in~\eqref{e.general_path_derivatives} yield
\begin{equation*}
    \la\partial_{q_1}\bar F_N,\partial_{q_2}\bar F_N\ra_{L^2}=\frac12\E\la R_{12}^{U'}\,\E\la R_{12}^{\sigma}\mid R_{12}^{\alpha}\ra\ra,
\end{equation*}
and combining this with~\eqref{e.general_time_derivative}, we obtain
\begin{equation}\label{e.general_enriched_hj}
    \partial_t\bar F_N-\la\partial_{q_1}\bar F_N,\partial_{q_2}\bar F_N\ra_{L^2}=\frac12\E\la R_{12}^{U'}\Ll(R_{12}^{\sigma}-\E\la R_{12}^{\sigma}\mid R_{12}^{\alpha}\ra\Rr)\ra.
\end{equation}
Since $|R_{12}^{U'}|\leq\frac{M_N}{N}\|U'\|_\infty^2$, the Cauchy--Schwarz inequality shows that the absolute value of the right-hand side of~\eqref{e.general_enriched_hj} is at most
\begin{equation}\label{e.general_concentration_error}
    \frac{M_N}{2N}\|U'\|_\infty^2\,\Delta_N^{1/2},\qquad\text{where}\quad\Delta_N:=\E\la\Ll(R_{12}^{\sigma}-\E\la R_{12}^{\sigma}\mid R_{12}^{\alpha}\ra\Rr)^2\ra.
\end{equation}
The quantity $\Delta_N$ measures the fluctuations of the spin overlap $R_{12}^\sigma$ conditionally on the overlap of the cascade $R_{12}^\alpha$; cf.\ \cite[(6.79)]{HJbook} or \cite[(2.23)]{mourrat2021nonconvex}. Due to the synchronization mechanism of \cite{pan.multi, pan.potts, pan.vec}, we expect this term to be small for generic $q$; see also \cite{engler2026synchronization}. The identity~\eqref{e.general_enriched_hj} thus leads us to expect that the function $\bar F_N$ converges, as $N$ tends to infinity, to a solution $f$ of the equation
\begin{equation}\label{e.general_path_hj}
    \partial_tf-\la\partial_{q_1}f,\partial_{q_2}f\ra_{L^2}=0.
\end{equation}
This is the equation~\eqref{e.moving_bilinear} studied in Section~\ref{s.moving_domain}, and we interpret it in the sense of Definition~\ref{d.moving_solution}. This definition involves a parameter $B$, which bounds the derivative of the solution with respect to~$q_2$. The appropriate choice is suggested by~\eqref{e.general_path_derivatives}: the derivatives $\partial_{q_1}\bar F_N$ and $\partial_{q_2}\bar F_N$ are bounded by $1$ and by $\frac{M_N}{2N}\|U'\|_\infty^2$ respectively, and the latter quantity converges to
\begin{equation}\label{e.gen.B}
    B:=\frac\lambda2\,\|U'\|_\infty^2.
\end{equation}

\smallskip
\noindent\textbf{The initial condition.}
At $t=0$, the quantities $S_j$ in~\eqref{e.def.S_j} do not depend on $\sigma$, and the exponential in~\eqref{e.general_enrichment} becomes a product of $N$ factors, each depending on one spin $\sigma_i$ and on the field $V_i^{q_1}$, and of $M_N$ factors, each depending on one variable $z_j$ and on the field $W_j^{q_2}$. As for~\eqref{e.psi=sumlambdapsi}, the additivity property of the cascade~\cite[Lemma~4.11]{chen2024ms} then gives
\begin{equation}\label{e.general_initial}
    \bar F_N(0,q)=\psi_{\mu_{\mcl I}}(q_1)+\frac{M_N}{N}\mathcal A(q_2),
\end{equation}
where $\psi_{\mu_{\mcl I}}$ is the cascade transform of the Ising measure defined in~\eqref{e.def.psi.rho}, and where, for every $q_2\in K_1$, we set
\begin{equation}\label{e.def.A}
    \mathcal A(q_2):=-\E\log\int\E_{z_1}\exp\Ll(U\Ll(W_1^{q_2}(\alpha)+\sqrt{1-q_2(1)}\,z_1\Rr)\Rr)\,\d\fR(\alpha).
\end{equation}
We call $\mathcal A$ the pattern functional associated with $U$. In view of~\eqref{e.general_initial}, the initial condition for the equation~\eqref{e.general_path_hj} must be
\begin{equation}\label{e.no_signal_hj}
    f(0,q)=\psi_{\mu_{\mcl I}}(q_1)+\lambda\mathcal A(q_2).
\end{equation}
The equation~\eqref{e.general_path_hj} with initial condition~\eqref{e.no_signal_hj} has the same structure as the equation~\eqref{e.main.hj} of the bipartite model: the nonlinearity is bilinear, and the initial condition is the sum of a function of $q_1$ and of a function of $q_2$. The first function does not depend on the activation function; it is the same as in~\eqref{e.main.hj}, and it is convex by Proposition~\ref{p.convex.psi.circ}. The second function is the only place where the activation function enters, and in general it is not convex. For instance, when $U=u_{s,\nu}$ for some $s\geq0$, integrating over $z_1$ in~\eqref{e.def.A} gives $\mathcal A(q_2)=\psi_\nu(sq_2)$, and we have seen in Remark~\ref{r.compare_hopf} that $\psi_\nu$ need not be convex.

The only property of the pattern functional that we need is the following one-sided bound on its increments. As discussed before the statement of Theorem~\ref{t.partial_hopf_moving}, this bound expresses that the derivative of $\mathcal A$ is an increasing path taking values in $[0,L^2/2]$, where $L$ is the Lipschitz constant of $U$, without requiring that $\mathcal A$ be differentiable. This is consistent with the second identity in~\eqref{e.general_path_derivatives}.

\begin{lemma}[Increments of the pattern functional]\label{l.pattern_slope}
Let $U:\R\to\R$ be a Lipschitz continuous function, with Lipschitz constant $L$, and let $\mathcal A$ be defined by~\eqref{e.def.A}. For every $q,q'\in K_1$, we have
\begin{equation}\label{e.pattern_slope}
    \mathcal A(q')-\mathcal A(q)\leq\frac{L^2}2\,\omega(q'-q),
\end{equation}
where $\omega$ is the support function of $K_1$ defined in~\eqref{e.moving_support}. In particular, the function $\psi_2:=\lambda\mathcal A$ satisfies the assumption~\eqref{e.moving_pattern_slope} of Theorem~\ref{t.partial_hopf_moving} with $B=\lambda L^2/2$, and the function $\mathcal A$ is Lipschitz continuous with respect to the $L^1$ norm, with Lipschitz constant $L^2/2$.
\end{lemma}

\begin{proof}
The last assertion follows from~\eqref{e.pattern_slope} and~\eqref{e.moving_support_bound}. We decompose the proof of~\eqref{e.pattern_slope} into three steps.

\smallskip
\noindent\emph{Step 1: Step paths and smooth activation functions.}
We first assume that $U$ is twice continuously differentiable with bounded first and second derivatives, so that $\|U'\|_\infty\leq L$, and we consider step paths. Let $k\geq1$ be an integer, let $0=s_0<s_1<\cdots<s_k=1$ and $0\leq q_0\leq\cdots\leq q_{k-1}\leq1$, and let $q:=\sum_{i=0}^{k-1}q_i\one_{[s_i,s_{i+1})}\in K_1$. We also set $q_{-1}:=0$ and $q_k:=1$. For such a path, the right-hand side of~\eqref{e.def.A} can be computed by a recursive Gaussian integration; see~\cite[Theorem~5.25]{HJbook}. To state this, we let $Z$ be a standard Gaussian random variable, and for every $s\in[0,1]$, $v\geq0$ and every Lipschitz function $h:\R\to\R$, we define $T_{s,v}h:\R\to\R$ by
\begin{equation*}
    T_{s,v}h(x):=\frac1s\log\E\exp\Ll(s\,h\big(x+\sqrt v\,Z\big)\Rr)\quad\text{if }s>0,
    \qquad\text{and}\qquad T_{0,v}h(x):=\E\Ll[h\big(x+\sqrt v\,Z\big)\Rr].
\end{equation*}
Setting $X_{k+1}:=U$ and then, successively for $i=k,k-1,\ldots,0$,
\begin{equation}\label{e.pattern_recursion}
    X_i:=T_{s_i,\,q_i-q_{i-1}}X_{i+1},
\end{equation}
we have $\mathcal A(q)=-X_0(0)$. The step $i=k$ of this recursion is the average over $z_1$ in~\eqref{e.def.A}, since $s_k=1$ and $q_k-q_{k-1}=1-q(1)$, and the step $i=0$ is a plain Gaussian average, since $s_0=0$.

We will use two classical properties of the operator $T_{s,v}$, valid for every function $h$ that is twice continuously differentiable with bounded first and second derivatives. For every $x\in\R$, we denote by $\E_x^{s,v,h}$ the expectation with respect to the probability measure whose density with respect to the law of $Z$ is proportional to $\exp(s\,h(x+\sqrt v\,Z))$. The first property is that the derivative of the mapping $h\mapsto T_{s,v}h(x)$ in the direction of a bounded continuous function $\eta$ is $\E_x^{s,v,h}[\eta(x+\sqrt v\,Z)]$, and that, similarly,
\begin{equation}\label{e.pattern_T_gradient}
    (T_{s,v}h)'(x)=\E_x^{s,v,h}\Ll[h'\big(x+\sqrt v\,Z\big)\Rr].
\end{equation}
The second property is that, for every $v>0$ and $x\in\R$,
\begin{equation}\label{e.pattern_T_pde}
    \partial_vT_{s,v}h(x)=\frac12\Ll((T_{s,v}h)''+s\big((T_{s,v}h)'\big)^2\Rr)(x)=\frac12\,\E_x^{s,v,h}\Ll[\Ll(h''+s(h')^2\Rr)\big(x+\sqrt v\,Z\big)\Rr].
\end{equation}
Indeed, when $s>0$, the function $w(v,x):=\exp(s\,T_{s,v}h(x))=\E\Ll[\exp\Ll(s\,h\big(x+\sqrt v\,Z\big)\Rr)\Rr]$ solves the heat equation $\partial_vw=\frac12\partial_x^2w$, so each of the three quantities in~\eqref{e.pattern_T_pde} is equal to $\partial_x^2w/(2s\,w)$ evaluated at $(v,x)$: for the second one, we express the derivatives of $w$ in terms of those of $T_{s,v}h$, and for the third one, we compute $\partial_x^2w$ by differentiating twice under the expectation. When $s=0$, the function $(v,x)\mapsto T_{0,v}h(x)$ itself solves the heat equation, and the argument is the same. 

Let $(\xi_i)_{-1\leq i\leq k}$ be the Markov chain started from $\xi_{-1}:=0$ and such that, for every $i\in\{0,\ldots,k\}$ and conditionally on $\xi_{i-1}$, the random variable $\xi_i$ has the law of $\xi_{i-1}+\sqrt{q_i-q_{i-1}}\,Z$ under the probability measure associated with $\E^{s_i,\,q_i-q_{i-1},\,X_{i+1}}_{\xi_{i-1}}$. By~\eqref{e.pattern_T_gradient} and~\eqref{e.pattern_recursion}, we have $X_i'(\xi_{i-1})=\E[X_{i+1}'(\xi_i)\mid\xi_{i-1}]$ for every $i\in\{0,\ldots,k\}$. In other words, setting $m_i:=X_{i+1}'(\xi_i)$, the sequence $(m_i)_{-1\leq i\leq k}$ is a martingale. Since $m_k=U'(\xi_k)$, this martingale is bounded by $L$, and by Jensen's inequality,
\begin{equation}\label{e.pattern_martingale}
    0\leq\E\Ll[m_0^2\Rr]\leq\E\Ll[m_1^2\Rr]\leq\cdots\leq\E\Ll[m_{k-1}^2\Rr]\leq L^2.
\end{equation}

We now assume that
\begin{equation}\label{e.pattern_strict}
    0<q_0<q_1<\cdots<q_{k-1}<1.
\end{equation}
The set of vectors $(q_0,\ldots,q_{k-1})$ satisfying~\eqref{e.pattern_strict} is an open subset of $\R^k$, and on this set, all the variance parameters $q_j-q_{j-1}$, $j\in\{0,\ldots,k\}$, appearing in the recursion~\eqref{e.pattern_recursion} are positive. Since the heat kernel depends smoothly on its variance parameter away from zero, the quantity $\mathcal A(q)$ is a continuously differentiable function of $(q_0,\ldots,q_{k-1})$ on this set. We fix $i\in\{0,\ldots,k-1\}$ and compute the derivative of $\mathcal A(q)$ with respect to $q_i$. Increasing $q_i$ increases the variance parameter $q_i-q_{i-1}$ of the step $i$ of the recursion~\eqref{e.pattern_recursion}, and it decreases by the same amount the variance parameter $q_{i+1}-q_i$ of the step $i+1$. By the first equality in~\eqref{e.pattern_T_pde}, the derivative of $X_{i+1}$ with respect to the variance parameter of the step $i+1$ is $\frac12(X_{i+1}''+s_{i+1}(X_{i+1}')^2)$. By the second equality in~\eqref{e.pattern_T_pde}, applied with $h=X_{i+1}$, the derivative of $X_i(x)$ with respect to the variance parameter of the step $i$ is
\begin{equation*}
    \E_x^{s_i,\,q_i-q_{i-1},\,X_{i+1}}\Ll[\eta\big(x+\sqrt{q_i-q_{i-1}}\,Z\big)\Rr],\qquad\text{with}\quad\eta:=\frac12\Ll(X_{i+1}''+s_i(X_{i+1}')^2\Rr),
\end{equation*}
and by the first of the two properties recalled above, this is also the derivative of $X_i(x)$ when $X_{i+1}$ is perturbed in the direction of $\eta$. The combined effect of these two variations on $X_i$ is therefore the same as that of perturbing $X_{i+1}$ in the direction of
\begin{equation*}
    \frac12\Ll(X_{i+1}''+s_i(X_{i+1}')^2\Rr)-\frac12\Ll(X_{i+1}''+s_{i+1}(X_{i+1}')^2\Rr)=-\frac{s_{i+1}-s_i}2\,(X_{i+1}')^2.
\end{equation*}
The terms involving second derivatives have cancelled out. Propagating this perturbation through the steps $i,i-1,\ldots,0$ of the recursion, using the first of the two properties recalled above at each step, we obtain
\begin{equation}\label{e.pattern_first_variation}
    \partial_{q_i}\mathcal A(q)=\frac{s_{i+1}-s_i}{2}\,\E\Ll[m_i^2\Rr].
\end{equation}

Let now $q'$ be another step path defined on the same partition $(s_i)_{0\leq i\leq k}$, and whose values $(q_i')_{0\leq i<k}$ also satisfy the strict inequalities in~\eqref{e.pattern_strict}, that is, $0<q_0'<\cdots<q_{k-1}'<1$. For every $a\in[0,1]$, the path $q(a):=q+a(q'-q)$ is a step path in $K_1$ defined on this partition, its values also satisfy the strict inequalities in~\eqref{e.pattern_strict}, and we denote by $(m_i(a))_{-1\leq i\leq k}$ the martingale constructed above for this path. By~\eqref{e.pattern_first_variation} and the chain rule, we have
\begin{equation*}
    \frac{\d}{\d a}\mathcal A(q(a))=\sum_{i=0}^{k-1}(q_i'-q_i)\,\frac{s_{i+1}-s_i}2\,\E\Ll[m_i(a)^2\Rr]=\la c(a),q'-q\ra_{L^2},
\end{equation*}
where
\begin{equation*}
    c(a):=\frac12\sum_{i=0}^{k-1}\E\Ll[m_i(a)^2\Rr]\one_{[s_i,s_{i+1})}.
\end{equation*}
By~\eqref{e.pattern_martingale}, the path $c(a)$ belongs to $K_{L^2/2}$. Since $\frac{L^2}2\,\omega$ is the support function of $K_{L^2/2}$, the derivative above is at most $\frac{L^2}2\,\omega(q'-q)$, and integrating over $a\in[0,1]$ yields~\eqref{e.pattern_slope} for every pair of step paths that are defined on a common partition and whose values satisfy the strict inequalities in~\eqref{e.pattern_strict}.

\smallskip
\noindent\emph{Step 2: General paths.}
We keep the assumption that $U$ is smooth, and let $q,q'\in K_1$. With the notation of~\eqref{e.dyadic_maps}, for every $j\in\N$, we let $e_j\in E_j$ be the vector with coordinates $(e_j)_k:=(k-\frac12)2^{-j}$, $k\in\{1,\ldots,2^j\}$, and we set
\begin{equation*}
    q^{(j)}:=\iota_j\Ll((1-2^{-j})\pi_jq+2^{-j}e_j\Rr)\qquad\text{and}\qquad q'^{(j)}:=\iota_j\Ll((1-2^{-j})\pi_jq'+2^{-j}e_j\Rr).
\end{equation*}
The coordinates of $\pi_jq$ and of $\pi_jq'$ are increasing and take values in $[0,1]$, while those of $e_j$ are strictly increasing and take values in $(0,1)$. The paths $q^{(j)}$ and $q'^{(j)}$ are therefore step paths defined on the same dyadic partition, their values satisfy the strict inequalities in~\eqref{e.pattern_strict}, and by Step~1, the inequality~\eqref{e.pattern_slope} holds for this pair of paths. Moreover, since the dyadic averages $\iota_j\pi_jq$ and $\iota_j\pi_jq'$ converge to $q$ and $q'$ in $L^1$ as $j$ tends to infinity, so do the paths $q^{(j)}$ and $q'^{(j)}$. By~\eqref{e.moving_support_bound} and the subadditivity of $\omega$, the function $\omega$ is continuous with respect to the $L^1$ norm. The functional $\mathcal A$ is also continuous on $K_1$ with respect to the $L^1$ norm; this can be seen by a Gaussian interpolation between the fields associated with two paths, as for the function $\psi_\rho$ in Lemma~\ref{l.psi^vec_smooth}. We can therefore pass to the limit in the inequality~\eqref{e.pattern_slope} for $q^{(j)}$ and $q'^{(j)}$, and we obtain that it also holds for $q$ and $q'$.

\smallskip
\noindent\emph{Step 3: Lipschitz activation functions.}
Finally, let $U$ be Lipschitz continuous with Lipschitz constant $L$. For every $\eps>0$, let $U_\eps$ be the convolution of $U$ with a smooth probability density supported in $[-\eps,\eps]$. The function $U_\eps$ is twice continuously differentiable with bounded first and second derivatives, it is Lipschitz continuous with Lipschitz constant $L$, and $\sup_{x\in\R}|U_\eps(x)-U(x)|\leq L\eps$. The quantity in~\eqref{e.def.A} changes by at most $L\eps$ when $U$ is replaced by $U_\eps$, uniformly over $q_2\in K_1$. By Steps~1 and~2, the inequality~\eqref{e.pattern_slope} holds for the pattern functional associated with $U_\eps$, and letting $\eps$ tend to zero, we conclude that it also holds for $\mathcal A$.
\end{proof}

\smallskip
\noindent\textbf{The candidate limit.}
We can now identify the solution of the equation~\eqref{e.general_path_hj} with initial condition~\eqref{e.no_signal_hj}. We recall that the function $\psi_{\mu_{\mcl I}}^*$ is defined in~\eqref{e.conjugate.psiI}.

\begin{proposition}[Candidate limit for the Ising perceptron]\label{p.gen.candidate}
Let $B$ be as in~\eqref{e.gen.B}. For every $(t,q)\in\widehat{\mathcal D}$, let
\begin{equation}\label{e.no_signal_candidate}
    f(t,q):=\sup_{p\in\mcl Q_{\infty,\leq1}}\Ll\{\la p,q_1\ra_{L^2}+\lambda\mathcal A(q_2+tp)-\psi_{\mu_{\mcl I}}^*(p)\Rr\}.
\end{equation}
The supremum in~\eqref{e.no_signal_candidate} is attained. The function $f$ is the unique Lipschitz viscosity solution of~\eqref{e.general_path_hj} with derivative bound $B$, in the sense of Definition~\ref{d.moving_solution}, such that $f(0,q)=\psi_{\mu_{\mcl I}}(q_1)+\lambda\mathcal A(q_2)$ for every $q\in\mcl Q_2\times K_1$. Moreover, for every $(t,q),(t',q')\in\widehat{\mathcal D}$, we have
\begin{equation}\label{e.gen.candidate_lip}
    |f(t,q)-f(t',q')|\leq|q_1-q_1'|_{L^1}+B|q_2-q_2'|_{L^1}+B|t-t'|.
\end{equation}
\end{proposition}

We recall that, for every $(t,q)\in\widehat{\mathcal D}$ and $p\in\mcl Q_{\infty,\leq1}$, we have $q_2+tp\leq(1-t)+t=1$, so that the path $q_2+tp$ belongs to $K_1$ and the right-hand side of~\eqref{e.no_signal_candidate} is well-defined.

\begin{proof}[Proof of Proposition~\ref{p.gen.candidate}]
We apply Theorem~\ref{t.partial_hopf_moving} with $\psi_2:=\lambda\mathcal A$, and with the function $g:\cH\to\R\cup\{+\infty\}$ defined by $g:=\psi_{\mu_{\mcl I}}^*$ on $K_1$ and $g:=+\infty$ on $\cH\setminus K_1$. By Lemma~\ref{l.pattern_slope}, the function $\psi_2$ satisfies the assumption~\eqref{e.moving_pattern_slope}. By the definition in~\eqref{e.conjugate.psiI}, the function $\psi_{\mu_{\mcl I}}^*$ is a supremum of functions of $p$ that are affine and continuous for the topology of $\cH$. Since $K_1$ is convex and closed in $\cH$, the function $g$ is convex and lower semicontinuous. It is also proper, since $\psi_{\mu_{\mcl I}}^*(0)=0$, as observed in Step~1 of the proof of Theorem~\ref{t.partial_hopf_paths}. Theorem~\ref{t.partial_hopf_moving} therefore applies, and it yields all the announced properties, provided that we identify the initial condition in~\eqref{e.moving_initial}; that is, we need to show that, for every $q_1\in\mcl Q_2$,
\begin{equation}\label{e.gen.biconjugate}
    \sup_{p\in\mcl Q_{\infty,\leq1}}\Ll\{\la p,q_1\ra_{L^2}-\psi_{\mu_{\mcl I}}^*(p)\Rr\}=\psi_{\mu_{\mcl I}}(q_1).
\end{equation}
This identity is contained in Theorem~\ref{t.partial_hopf_paths}. Indeed, this theorem states that the right-hand side of~\eqref{e.main.1}, in which the supremum can be restricted to $\mcl Q_{\infty,\leq1}$, is equal to $\psi_{\mu_{\mcl I}}(q_1)+\lambda\psi_\nu(q_2)$ at $t=0$, and at $t=0$, the term $\lambda\psi_\nu(q_2)$ can be taken out of the supremum.
\end{proof}

In view of the formal derivation of the equation~\eqref{e.general_path_hj} and of the identity~\eqref{e.general_initial}, we expect that, for every $(t,q)\in\mathcal D$,
\begin{equation}\label{e.gen.convergence}
    \lim_{N\to\infty}\bar F_N(t,q)\stackrel{?}{=}f(t,q),
\end{equation}
where $f$ is the function defined in~\eqref{e.no_signal_candidate}. Evaluating this at $(t,q)=(1,0)$, we obtain the following prediction for the limit of the free energy in~\eqref{e.general_free_energy}:
\begin{equation}\label{e.general_free_energy_candidate}
    \lim_{N\to\infty}\bar F_N(1,0)\stackrel{?}{=}f(1,0)=\sup_{p\in\mcl Q_{\infty,\leq1}}\Ll\{\lambda\mathcal A(p)-\psi_{\mu_{\mcl I}}^*(p)\Rr\}.
\end{equation}

\smallskip
\noindent\textbf{Comparison with the formula of~\cite{openai2026ising}.}
When $U$ is bounded, the limit of the free energy in~\eqref{e.general_free_energy} was identified in~\cite[Theorem~2.3]{openai2026ising}, and the result coincides with the right-hand side of~\eqref{e.general_free_energy_candidate}. To see this, we first note that in~\cite{openai2026ising}, the free energy is defined without the minus sign in front of the logarithm. Accordingly, the variational formula obtained there is an infimum over paths $p\in\mcl Q_{\infty,\leq1}$, and the functional to be minimized is the sum of two terms. The first term is $\lambda$ times the pattern functional of~\cite[Definition~2.1]{openai2026ising}. For step paths, this functional is defined by the Gaussian recursion~\eqref{e.pattern_recursion}, including the final Gaussian average with variance $1-p(1)$, which corresponds to the average over $z_1$ in~\eqref{e.def.A}; it is therefore equal to $-\mathcal A(p)$. The second term is the Ising entropy of~\cite[Definition~2.2]{openai2026ising}, which is the supremum over increasing step paths $h$ of $\ell(h)+\frac12\la h,p\ra_{L^2}$, where, in our notation, $\ell(h)=-\psi_{\mu_{\mcl I}}(h/2)$. The change of variables $h=2r$ shows that this entropy is equal to $\psi_{\mu_{\mcl I}}^*(p)$. The recipe described here thus recovers the formula of~\cite{openai2026ising}, and interprets it as the value at $(t,q)=(1,0)$ of the solution to a Hamilton--Jacobi equation with the structure of Theorem~\ref{t.hopf.half}. The same formula was also obtained in~\cite{fuSohn2026perceptrons} for Lipschitz activation functions; see Section~\ref{s.random_activations}.

\begin{remark}[Back to the bipartite model]\label{r.gen.bipartite}
Let $\nu$ be a probability measure supported in $[-1,1]$, let $s\geq0$, and suppose that $U=u_{s,\nu}$, with $u_{s,\nu}$ as in~\eqref{e.u_t,nu=}. For clarity, we denote here by $\bar F_N^\bip$ the enriched free energy of the bipartite model defined in~\eqref{e.F_N(t,q)=}, and by $f^\bip$ its limit, given by Theorem~\ref{t.main}. Inserting the definition of $u_{s,\nu}$ into~\eqref{e.general_enrichment}, integrating out the variables $(z_j)$, and using that the process $\sqrt s\,W_j^{q_2}$ has the same law as $W_j^{sq_2}$, we find that, for every $(t,q)\in\mathcal D$,
\begin{equation}\label{e.gen.bipartite}
    \bar F_N(t,q_1,q_2)=\bar F_N^\bip(st,q_1,sq_2).
\end{equation}
In this computation, the integration over the variables $(z_j)$ turns the term $-s|v|^2$ coming from the definition of $u_{s,\nu}$ into $-s(t+q_2(1))|v|^2$, where $v$ denotes the layer of $\nu$-distributed spins of the bipartite model; this is the correction term in~\eqref{e.H^t,q_N=}, evaluated at $(st,sq_2)$ in place of $(t,q_2)$. The variables $(z_j)$ thus play here the role that this correction term plays in the bipartite model. Similarly, since $\mathcal A(q_2)=\psi_\nu(sq_2)$ in this case, comparing~\eqref{e.no_signal_candidate} with~\eqref{e.main.1} gives $f(t,q_1,q_2)=f^\bip(st,q_1,sq_2)$. For these activation functions, the convergence~\eqref{e.gen.convergence} is therefore a consequence of Theorem~\ref{t.main}. The identity~\eqref{e.gen.bipartite} also clarifies the relation between the time variables of the two settings. In the bipartite model, the time variable is the parameter $s$ of the activation function $u_{s,\nu}$, which plays the role of an inverse temperature. In the present section, the activation function is fixed, and the time variable $t\in[0,1]$ is the fraction of the variance of the quantity $S_j$ in~\eqref{e.def.S_j} that comes from the pre-activation of the original model.
\end{remark}

\subsection{The spherical perceptron}
\label{s.spherical}

We now replace the Ising spins by spherical spins. We keep the activation function $U$ of Section~\ref{s.gen.ising}, and we denote by $\mu_{\mathbb{S}_N}$ the uniform probability measure on the sphere $\Sph_N:=\{\sigma\in\R^N:|\sigma|^2=N\}$. With $S_j$ as in~\eqref{e.def.S_j}, we define the enriched free energy $\bar F_N:\mathcal D\to\R$ by
\begin{equation}\label{e.sph_enrichment}
    \bar F_N(t,q):=-\frac1N\E\log\iint\E_z\exp\Ll(\sum_{j=1}^{M_N}U(S_j)+\sqrt{2}\sum_{i=1}^NV_i^{q_1}(\alpha)\sigma_i-Nq_1(1)\Rr)\,\d\mu_{\mathbb{S}_N}(\sigma)\,\d\fR(\alpha).
\end{equation}
This is the definition in~\eqref{e.general_enrichment} with $\mu_{\mathbb{S}_N}$ in place of $\mu_{\mcl I}^{\otimes N}$, and $\bar F_N(1,0)$ is the free energy of the spherical perceptron with activation function $U$, that is, the quantity in~\eqref{e.def.perc.general} with $P_N=\mu_{\mathbb{S}_N}$ and $M=M_N$.

\smallskip
\noindent\textbf{The approximate Hamilton--Jacobi equation.}
The derivation of the identity~\eqref{e.general_enriched_hj} in Section~\ref{s.gen.ising} only uses one property of the Ising spins, namely the identity $|\sigma|^2=N$. This identity ensures that the variance of $S_j$ is equal to one, and that the term $-Nq_1(1)$ in~\eqref{e.general_enrichment} is one half of the variance of the Gaussian field acting on the spins. Since it also holds on $\Sph_N$, the identities~\eqref{e.general_path_derivatives},~\eqref{e.general_time_derivative} and~\eqref{e.general_enriched_hj} remain valid for the function $\bar F_N$ in~\eqref{e.sph_enrichment}, provided that the replicas are now sampled from the Gibbs measure associated with~\eqref{e.sph_enrichment}. The candidate limit $f$ is therefore again expected to solve the equation~\eqref{e.general_path_hj}, in the sense of Definition~\ref{d.moving_solution} with $B$ as in~\eqref{e.gen.B}, and the error term in the equation satisfied by $\bar F_N$ is still controlled by the quantity $\Delta_N$ in~\eqref{e.general_concentration_error}.

\smallskip
\noindent\textbf{The initial condition.}
At $t=0$, the additivity property of the cascade still separates the contribution of the patterns from that of the spins, and we get
\begin{equation}\label{e.sph_initial}
    \bar F_N(0,q)=\psi_{\sph,N}(q_1)+\frac{M_N}{N}\mathcal A(q_2),
\end{equation}
where $\mathcal A$ is the pattern functional defined in~\eqref{e.def.A}, and where
\begin{equation}\label{e.def.psi.sph.N}
    \psi_{\sph,N}(q_1):=-\frac1N\E\log\iint\exp\Ll(\sqrt2\sum_{i=1}^NV_i^{q_1}(\alpha)\sigma_i-Nq_1(1)\Rr)\,\d\mu_{\mathbb{S}_N}(\sigma)\,\d\fR(\alpha).
\end{equation}
The part of the initial condition that depends on $q_2$ is thus the same as for Ising spins, and only the part that depends on $q_1$ is modified. Since $\mu_{\mathbb{S}_N}$ is not a product measure, the quantity $\psi_{\sph,N}(q_1)$ does not reduce to a computation involving a single spin, and it depends on $N$. Its limit
\begin{equation*}
    \psi_\sph(q_1):=\lim_{N\to\infty}\psi_{\sph,N}(q_1)
\end{equation*}
exists, and it is a convex function of $q_1$; see~\cite[Proposition~2.11 and Lemma~2.9]{chen2026spherical}. This is the spherical counterpart of Proposition~\ref{p.convex.psi.circ}. Moreover, as will be explained just below, we have for every $q_1\in\mcl Q_\infty$ that
\begin{equation}\label{e.psi.sph.dual}
    \psi_\sph(q_1)=\sup_{p\in\mcl Q_{\infty,\leq1}}\Ll\{\la p,q_1\ra_{L^2}-\mathsf S(p)\Rr\},
\end{equation}
where the functional $\mathsf S:\mcl Q_{\infty,\leq1}\to[0,+\infty]$ is defined by
\begin{equation}\label{e.def.S.sph}
    \mathsf S(p):=\frac12\int_0^1\Ll(\frac1{D_p(v)}-\frac1{1-v}\Rr)\d v,\qquad\text{with}\quad D_p(v):=1-v-\int_0^1\big(p(u)-v\big)_+\,\d u,
\end{equation}
and with the understanding that $1/0=+\infty$. Denoting by $\gamma_p(v):=|\{u\in[0,1):p(u)\leq v\}|$ the distribution function of the path $p$, we can also write $D_p(v)=\int_v^1\gamma_p(w)\,\d w$. The identity~\eqref{e.psi.sph.dual} is a consequence of~\cite[Lemmas~2.2--2.3]{ali2026free}. Indeed, let $\psi_\sph^*$ denote the convex conjugate of $\psi_\sph$, defined as in~\eqref{e.conjugate.psiI}. By~\cite[Lemma~2.3]{ali2026free}, we have $\psi_\sph^*(p)=\mathsf S(p)$ for every $p\in\mcl Q_{\infty,\leq1}$ such that $p(1)<1$. In view of the definition of $\psi_\sph^*$, this implies that the inequality
\begin{equation*}
    \la p,q_1\ra_{L^2}-\mathsf S(p)\leq\psi_\sph(q_1)
\end{equation*}
holds for every such~$p$, and by~\cite[Lemma~2.2]{ali2026free}, it is an equality for at least one such~$p$. It thus remains to verify that this inequality is also valid when $p(1)=1$. For every $\theta\in(0,1)$, since $p\geq0$, we have $D_{\theta p}\geq D_p$, and thus $\mathsf S(\theta p)\leq\mathsf S(p)$. Applying the inequality to the path $\theta p$ therefore gives $\theta\la p,q_1\ra_{L^2}-\mathsf S(p)\leq\psi_\sph(q_1)$, and we conclude by letting $\theta$ tend to~$1$.

\smallskip
\noindent\textbf{The candidate limit.}
The representation~\eqref{e.psi.sph.dual} provides us with an explicit convex function whose convex conjugate is $\psi_\sph$. We extend $\mathsf S$ to $\cH$ by setting $\mathsf S:=+\infty$ outside of $K_1$. For every $v\in[0,1)$, the mapping $p\mapsto D_p(v)$ is concave and continuous on $K_1$, and it takes values in $[0,1-v]$. Since the function $x\mapsto1/x$ is convex and decreasing on $[0,\infty)$, the integrand in~\eqref{e.def.S.sph} is a nonnegative, convex and lower semicontinuous function of $p$, and by Fatou's lemma, the functional $\mathsf S$ is convex and lower semicontinuous. It is also proper, since $\mathsf S(0)=0$. We can therefore apply Theorem~\ref{t.partial_hopf_moving} with $g:=\mathsf S$, with $\psi_2:=\lambda\mathcal A$, which satisfies the assumption~\eqref{e.moving_pattern_slope} by Lemma~\ref{l.pattern_slope}, and with $B$ as in~\eqref{e.gen.B}. We obtain that the function $f:\widehat{\mathcal D}\to\R$ defined by
\begin{equation}\label{e.sph_candidate}
    f(t,q):=\sup_{p\in\mcl Q_{\infty,\leq1}}\Ll\{\la p,q_1\ra_{L^2}+\lambda\mathcal A(q_2+tp)-\mathsf S(p)\Rr\}
\end{equation}
is the unique Lipschitz viscosity solution of~\eqref{e.general_path_hj} with derivative bound $B$ such that, for every $q\in\mcl Q_2\times K_1$,
\begin{equation}\label{e.sph_candidate_initial}
    f(0,q)=\psi_\sph(q_1)+\lambda\mathcal A(q_2),
\end{equation}
where we have extended the identity~\eqref{e.psi.sph.dual} to every $q_1\in\mcl Q_2$ by continuity. In view of~\eqref{e.sph_initial}, this function is our candidate for the limit of $\bar F_N$, and at $(t,q)=(1,0)$, it yields the prediction
\begin{equation}\label{e.sph_free_energy_candidate}
    \lim_{N\to\infty}\bar F_N(1,0)\stackrel{?}{=}\sup_{p\in\mcl Q_{\infty,\leq1}}\Ll\{\lambda\mathcal A(p)-\mathsf S(p)\Rr\}.
\end{equation}

\smallskip
\noindent\textbf{Comparison with the formulas of~\cite{aronow2026negative, openai2026spherical}.}
For activation functions satisfying the assumptions of this subsection and bounded above, the limit of $\bar F_N(1,0)$ was identified in~\cite[Part~I, Theorem~1.1]{aronow2026negative}. The opposite sign convention for the free energy is used there, and their formula is therefore an infimum. Their order parameter is the distribution function $\gamma_p$ of our path $p$. Under this identification, their Crisanti--Sommers entropy equals $\mathsf S(p)$, and their pattern functional  equals $-\mathcal A(p)$, so the functional they minimize is $\mathsf S(p)-\lambda\mathcal A(p)$. The restriction $p(1)<1$ in their parametrization can be removed by approximation, and their formula thus agrees with~\eqref{e.sph_free_energy_candidate}. 
When $U$ is bounded and continuous, the limit of $\bar F_N(1,0)$ was also identified in~\cite[Theorem~1]{openai2026spherical}. The variational problem there is again parametrized by distribution functions $\gamma_p$, with spherical term $\mathsf S(p)$ and a pattern functional expressed through a stochastic control problem whose value is $-\mathcal A(p)$. After accounting for the opposite sign convention, this gives the same formula as in~\eqref{e.sph_free_energy_candidate}.

\subsection{Random activation functions}
\label{s.random_activations}

We now allow the activation function to vary from one pattern to the next, in an i.i.d.\ fashion, and we treat Ising and spherical spins simultaneously. This is the setting of~\cite{fuSohn2026perceptrons}. Let $\pi$ be a probability measure on $\R$, let $(\eta_j)_{j\in\N}$ be i.i.d.\ random variables with law $\pi$, independent of $(g_{ij})$, and let $(\eta,x)\mapsto U_\eta(x)$ be a measurable function from $\R^2$ to $\R$. We call $\eta_j$ the mark of the $j$-th pattern. We set $P_N^{\mathrm I}:=\mu_{\mcl I}^{\otimes N}$ and $P_N^{\mathrm S}:=\mu_{\mathbb{S}_N}$, and for $\kappa\in\{\mathrm I,\mathrm S\}$, we consider the model with free energy
\begin{equation}\label{e.random_activation_model}
    -\frac1N\E\log\int\exp\Ll(\sum_{j=1}^{M_N}U_{\eta_j}\Ll(\frac1{\sqrt N}\sum_{i=1}^Ng_{ij}\sigma_i\Rr)\Rr)\,\d P_N^\kappa(\sigma),
\end{equation}
where, in this subsection, the expectation $\E$ also averages over the marks $(\eta_j)$. The models of Sections~\ref{s.gen.ising} and~\ref{s.spherical} correspond to the case in which $\pi$ is a Dirac mass. We assume that there exists a constant $L<\infty$ such that, for $\pi$-almost every $\eta$, the function $U_\eta$ is twice continuously differentiable with $\|U_\eta'\|_\infty\leq L$ and $\|U_\eta''\|_\infty\leq L$, and that $\int|U_\eta(0)|\,\d\pi(\eta)<\infty$. We explain at the end of this subsection how to remove the smoothness assumption.

\smallskip
\noindent\textbf{The enriched free energy and the approximate Hamilton--Jacobi equation.}
The enriched free energy $\bar F_N:\mathcal D\to\R$ is defined as in~\eqref{e.general_enrichment}, with the function $U_{\eta_j}$ applied to $S_j$ and with the measure $P_N^\kappa$ on the spins, that is,
\begin{equation}\label{e.random_activation_enrichment}
    \bar F_N(t,q):=-\frac1N\E\log\iint\E_z\exp\Ll(\sum_{j=1}^{M_N}U_{\eta_j}(S_j)+\sqrt{2}\sum_{i=1}^NV_i^{q_1}(\alpha)\sigma_i-Nq_1(1)\Rr)\,\d P_N^\kappa(\sigma)\,\d\fR(\alpha),
\end{equation}
with $S_j$ as in~\eqref{e.def.S_j}. The quantity $\bar F_N(1,0)$ is the free energy in~\eqref{e.random_activation_model}. The computations of Section~\ref{s.gen.ising} can be performed conditionally on the marks, and they show that the identities~\eqref{e.general_path_derivatives},~\eqref{e.general_time_derivative} and~\eqref{e.general_enriched_hj} remain valid, with
\begin{equation*}
    R_{12}^{U'}:=\frac1N\sum_{j=1}^{M_N}U'_{\eta_j}(S_j^1)U'_{\eta_j}(S_j^2),
    \qquad\text{so that}\quad|R_{12}^{U'}|\leq\frac{M_N}{N}L^2.
\end{equation*}
The candidate limit $f$ is thus again expected to solve the equation~\eqref{e.general_path_hj}, in the sense of Definition~\ref{d.moving_solution} with $B:=\lambda L^2/2$, and the error term in the equation satisfied by $\bar F_N$ is still controlled by the quantity $\Delta_N$ in~\eqref{e.general_concentration_error}.

\smallskip
\noindent\textbf{The initial condition and the candidate limit.}
For $\pi$-almost every $\eta$, we write $\mathcal A_{U_\eta}$ for the pattern functional defined in~\eqref{e.def.A} with $U_\eta$ in place of $U$, and for every $q_2\in K_1$, we set
\begin{equation}\label{e.def.A.pi}
    \mathcal A_\pi(q_2):=\int\mathcal A_{U_\eta}(q_2)\,\d\pi(\eta).
\end{equation}
Notice that the average over the mark is taken outside of the logarithm in~\eqref{e.def.A}. Conditioning on the marks and using the additivity property of the cascade, we see that
\begin{equation*}
    \bar F_N(0,q)=\psi_{\kappa,N}(q_1)+\frac{M_N}N\mathcal A_\pi(q_2),
\end{equation*}
where $\psi_{\mathrm I,N}:=\psi_{\mu_{\mcl I}}$ and $\psi_{\mathrm S,N}:=\psi_{\sph,N}$. Setting $\psi_{\mathrm I}:=\psi_{\mu_{\mcl I}}$ and $\psi_{\mathrm S}:=\psi_\sph$, the initial condition for the equation~\eqref{e.general_path_hj} is therefore
\begin{equation}\label{e.random_activation_initial}
    f(0,q)=\psi_\kappa(q_1)+\lambda\mathcal A_\pi(q_2).
\end{equation}
Compared with the previous two subsections, only the part of the initial condition that depends on $q_2$ has changed. Averaging the inequality~\eqref{e.pattern_slope} over the marks, we see that the function $\psi_2:=\lambda\mathcal A_\pi$ satisfies the assumption~\eqref{e.moving_pattern_slope} with $B=\lambda L^2/2$. We also set
\begin{equation*}
    \psi_{\mathrm I}^*:=\psi_{\mu_{\mcl I}}^*\qquad\text{and}\qquad\psi_{\mathrm S}^*:=\mathsf S,
\end{equation*}
and we recall from the proof of Proposition~\ref{p.gen.candidate} and from Section~\ref{s.spherical} that, for each $\kappa\in\{\mathrm I,\mathrm S\}$, the function equal to $\psi_\kappa^*$ on $K_1$ and to $+\infty$ elsewhere is proper, lower semicontinuous and convex, and that its convex conjugate is $\psi_\kappa$. By Theorem~\ref{t.partial_hopf_moving}, the function $f:\widehat{\mathcal D}\to\R$ defined by
\begin{equation}\label{e.random_activation_candidate}
    f(t,q):=\sup_{p\in\mcl Q_{\infty,\leq1}}\Ll\{\la p,q_1\ra_{L^2}+\lambda\mathcal A_\pi(q_2+tp)-\psi_\kappa^*(p)\Rr\}
\end{equation}
is therefore the unique Lipschitz viscosity solution of~\eqref{e.general_path_hj} with derivative bound $B$ and with initial condition~\eqref{e.random_activation_initial}. This formula contains~\eqref{e.no_signal_candidate} and~\eqref{e.sph_candidate} as special cases. At $(t,q)=(1,0)$, it yields the prediction
\begin{equation}\label{e.random_activation_free_energy}
    \lim_{N\to\infty}\bar F_N(1,0)\stackrel{?}{=}\sup_{p\in\mcl Q_{\infty,\leq1}}\Ll\{\lambda\mathcal A_\pi(p)-\psi_\kappa^*(p)\Rr\}.
\end{equation}

\smallskip
\noindent\textbf{Comparison with the formulas of~\cite{fuSohn2026perceptrons}.}
The right-hand side of~\eqref{e.random_activation_free_energy} coincides with the formula proved in~\cite[Theorem~1.1]{fuSohn2026perceptrons}, as we now explain. For every $p\in\mcl Q_{\infty,\leq1}$, let $\gamma_p$ be the distribution function of $p$ introduced below~\eqref{e.def.S.sph}. For every Lipschitz function $U:\R\to\R$, the channel functional of~\cite{fuSohn2026perceptrons} is defined by $\Phi_U(p):=\varphi(0,0)$, where $\varphi:[0,1]\times\R\to\R$ is the solution of the Parisi equation
\begin{equation}\label{e.fu_sohn_channel_pde}
    \partial_v\varphi+\frac12\Ll(\partial_x^2\varphi+\gamma_p(v)(\partial_x\varphi)^2\Rr)=0\quad\text{on }[0,1]\times\R,
    \qquad\varphi(1,\cdot)=U.
\end{equation}
When $p$ is a step path, the function $\gamma_p$ is piecewise constant, and on each interval on which $\gamma_p$ is equal to a constant $s$, the equation~\eqref{e.fu_sohn_channel_pde} is solved by the operators $T_{s,v}$ appearing in~\eqref{e.pattern_recursion}; see~\eqref{e.pattern_T_pde}. In particular, we have $\gamma_p=1$ on the interval $[p(1),1]$, and solving the equation on this interval produces the logarithm of the Gaussian average with variance $1-p(1)$ in~\eqref{e.def.A}. Solving~\eqref{e.fu_sohn_channel_pde} backward in $v$ thus reproduces the recursion~\eqref{e.pattern_recursion}, and we obtain that $\Phi_U=-\mathcal A_U$. We note that the paths considered in~\cite{fuSohn2026perceptrons} are defined on $[0,1]$ and take the value $1$ at the endpoint $1$; our path $p$ is the restriction of such a path to $[0,1)$, and $p(1)$ denotes its left limit at $1$. The functionals $\psi_{\mathrm I}^*$ and $\psi_{\mathrm S}^*$ coincide with the spin costs of~\cite{fuSohn2026perceptrons}, which are denoted there by $A_{\mathrm I}$ and $A_{\mathrm S}$ respectively. Indeed, the Ising cost $A_{\mathrm I}$ is defined in the same way as the Ising entropy of~\cite{openai2026ising} discussed in Section~\ref{s.gen.ising}, and it is therefore equal to $\psi_{\mu_{\mcl I}}^*$; and the spherical cost $A_{\mathrm S}$ is defined by the formula~\eqref{e.def.S.sph}. Changing the sign in~\eqref{e.random_activation_free_energy}, we can therefore rewrite this prediction as
\begin{multline}\label{e.fu_sohn_formula}
    \lim_{N\to\infty}\frac1N\E\log\int\exp\Ll(\sum_{j=1}^{M_N}U_{\eta_j}\Ll(\frac1{\sqrt N}\sum_{i=1}^Ng_{ij}\sigma_i\Rr)\Rr)\,\d P_N^\kappa(\sigma)\\
    \stackrel{?}{=}\inf_{p\in\mcl Q_{\infty,\leq1}}\Ll\{A_\kappa(p)+\lambda\int\Phi_{U_\eta}(p)\,\d\pi(\eta)\Rr\},
\end{multline}
and the right-hand side of~\eqref{e.fu_sohn_formula} is the formula of~\cite[Theorem~1.1]{fuSohn2026perceptrons}, in which the ratio $\lambda$ is denoted by $\alpha$. 

\smallskip
\noindent\textbf{On the assumptions and the mode of convergence.}
The result of~\cite{fuSohn2026perceptrons} does not require the activation functions to be smooth: it only assumes that, for $\pi$-almost every $\eta$, the function $U_\eta$ is Lipschitz continuous with Lipschitz constant at most $L$, and that $\int|U_\eta(0)|\,\d\pi(\eta)<\infty$. The smoothness assumption was only used here to perform the integrations by parts leading to~\eqref{e.general_enriched_hj}, and the prediction~\eqref{e.random_activation_free_energy} extends to this more general setting by regularization. Indeed, Lemma~\ref{l.pattern_slope} is valid for Lipschitz activation functions, so the function $f$ in~\eqref{e.random_activation_candidate} is well-defined in this setting, and it is still the solution of the Hamilton--Jacobi equation. Moreover, if $U_\eta^\eps$ denotes the convolution of $U_\eta$ with a smooth probability density supported in $[-\eps,\eps]$, then $\sup_{x\in\R}|U_\eta^\eps(x)-U_\eta(x)|\leq L\eps$. Replacing $U_\eta$ by $U_\eta^\eps$ therefore changes the free energy in~\eqref{e.random_activation_model} by at most $\frac{M_N}{N}L\eps$, and it changes the pattern functional $\mathcal A_\pi$ by at most $L\eps$, uniformly over $K_1$. Hence, if the prediction~\eqref{e.random_activation_free_energy} is valid for the regularized activation functions, then it is also valid for the original ones.

We also point out that~\cite[Theorem~1.1]{fuSohn2026perceptrons} concerns the random variable obtained by removing the expectation in~\eqref{e.random_activation_model}, which is shown there to converge in $L^1$; this implies the convergence of its expectation. Conversely, the fluctuations of this random variable can be controlled independently of the identification of the limit of its expectation. The constants $U_{\eta_j}(0)$ contribute the term $-N^{-1}\sum_{j=1}^{M_N}U_{\eta_j}(0)$, which obeys the law of large numbers. For the remaining part, resampling one pattern together with its mark changes the logarithm of the partition function by a quantity whose second moment is bounded in terms of $L$ only, and the Efron--Stein inequality thus shows that its variance is of order $N^{-1}$. The convergence of the expectation therefore implies the convergence of the random variable in $L^1$, and also in $L^2$ if $\int|U_\eta(0)|^2\,\d\pi(\eta)<\infty$. When the activation functions are uniformly bounded, the bounded differences inequality also gives the almost sure convergence, and the convergence in $L^r$ for every $r<\infty$.

\subsection{Labelled patterns from a Gaussian mixture}
\label{s.ext}

We next consider a model in which each pattern carries a label, and in which the patterns are no longer standard Gaussian vectors: they are sampled from a mixture of two Gaussian distributions, as in the continuous-spin perceptron of~\cite[Section~2]{camilli2026perceptron}. We point out however that, in~\cite{camilli2026perceptron}, the spins are real-valued with a log-concave prior, and the activation functions are assumed to be concave. Under these assumptions, the Gibbs measure is log-concave, which ensures that the overlaps concentrate; in other words, the model is replica-symmetric, and its free energy is described in terms of finitely many scalar parameters. Here, the spins are Ising spins and the activation functions are not assumed to be concave, so that replica symmetry breaking is expected in general, and the variational formula below involves paths rather than scalar parameters. 

Let $(y_j)_{j\in\N}$ be i.i.d.\ random variables with values in $\{-1,1\}$, and let $(\theta_i)_{i\in\N}$ be i.i.d.\ random variables whose law $\mu$ is a compactly supported probability measure on $\R$. These two families are independent of each other and of $(g_{ij})$, and in this subsection, the expectation $\E$ also averages over them. For each $y\in\{-1,1\}$, let $U_y:\R\to\R$ be a twice continuously differentiable function with bounded first and second derivatives. We consider the model with Ising spins whose free energy is
\begin{equation}\label{e.ext_free_energy}
    -\frac1N\E\log\int\exp\Ll(\sum_{j=1}^{M_N}U_{y_j}\Ll(\frac1{\sqrt N}\sum_{i=1}^Ng_{ij}\sigma_i+\frac{y_j}{N}\sum_{i=1}^N\theta_i\sigma_i\Rr)\Rr)\,\d\mu_{\mcl I}^{\otimes N}(\sigma).
\end{equation}
The argument of $U_{y_j}$ in~\eqref{e.ext_free_energy} can be rewritten as $N^{-1/2}\sum_{i=1}^N(g_{ij}+N^{-1/2}y_j\theta_i)\sigma_i$. The $j$-th pattern is thus the vector $(g_{ij})_{i\leq N}+N^{-1/2}y_j\theta$, where $\theta:=(\theta_1,\ldots,\theta_N)$: it is obtained by shifting a standard Gaussian vector by $N^{-1/2}\theta$ or by $-N^{-1/2}\theta$, according to the value of its label. Conditionally on $\theta$, the patterns are therefore independent samples from the mixture of the two Gaussian distributions with identity covariance and with means $N^{-1/2}\theta$ and $-N^{-1/2}\theta$, and the label $y_j$ selects both the component of the mixture from which the $j$-th pattern is sampled and the activation function. When $\mu=\delta_0$, the patterns are not shifted, and we recover the model of Section~\ref{s.random_activations} with Ising spins and with the marks $\eta_j=y_j$; if in addition $U_{-1}=U_1$, then we recover the model of Section~\ref{s.gen.ising}. 

\smallskip
\noindent\textbf{The enriched free energy.}
In addition to the pair of paths $q=(q_1,q_2)$, the enrichment of this model involves a pair of real parameters $x=(x_1,x_2)$: the parameter $x_1$ is an external field acting on the weighted magnetization, and the parameter $x_2$ shifts the pre-activations in the direction of the labels. We set $\mathcal D^\ext:=\mathcal D\times\R^2$. With the same cascade fields $(V_i^{q_1})$ and $(W_j^{q_2})$ and the same Gaussian random variables $(z_j)$ as in Section~\ref{s.gen.ising}, for every $(t,q,x)\in\mathcal D^\ext$, we replace the definition of $S_j$ in~\eqref{e.def.S_j} by
\begin{equation}\label{e.def.S_j.ext}
    S_j:=\sqrt{\frac tN}\sum_{i=1}^Ng_{ij}\sigma_i+\frac{ty_j}{N}\sum_{i=1}^N\theta_i\sigma_i+y_jx_2+W_j^{q_2}(\alpha)+\sqrt{1-t-q_2(1)}\,z_j,
\end{equation}
and we define $\bar F_N:\mathcal D^\ext\to\R$ by
\begin{multline}\label{e.ext_enrichment}
    \bar F_N(t,q,x):=-\frac1N\E\log\iint\E_z\exp\Ll(\sum_{j=1}^{M_N}U_{y_j}(S_j)\Rr.\\
    \Ll.{}+\sqrt{2}\sum_{i=1}^NV_i^{q_1}(\alpha)\sigma_i-Nq_1(1)-x_1\sum_{i=1}^N\theta_i\sigma_i\Rr)\,\d\mu_{\mcl I}^{\otimes N}(\sigma)\,\d\fR(\alpha).
\end{multline}
The quantity $\bar F_N(1,0,0)$ is the free energy in~\eqref{e.ext_free_energy}. When $\mu=\delta_0$, the function $\bar F_N$ does not depend on $x_1$, and at $x_2=0$, it reduces to the enriched free energy in~\eqref{e.random_activation_enrichment} with $\kappa=\mathrm I$ and $\eta_j=y_j$. Notice that in~\eqref{e.def.S_j.ext}, the shift of the patterns is multiplied by $t$, while their Gaussian part is multiplied by $\sqrt t$. This choice ensures that the contribution of the shift to the time derivative of $\bar F_N$ is matched by the product of the derivatives of $\bar F_N$ with respect to $x_1$ and $x_2$, as we now explain.

\smallskip
\noindent\textbf{The approximate Hamilton--Jacobi equation.}
Let $\la\cdot\ra$ denote the expectation with respect to the Gibbs measure associated with~\eqref{e.ext_enrichment}, and let the overlaps $R_{12}^\alpha$ and $R_{12}^\sigma$ be as in Section~\ref{s.gen.ising}. We also set
\begin{equation*}
    R_{12}^{U'}:=\frac1N\sum_{j=1}^{M_N}U'_{y_j}(S_j^1)U'_{y_j}(S_j^2),\qquad R_1^{yU'}:=\frac1N\sum_{j=1}^{M_N}y_jU'_{y_j}(S_j^1),\qquad\text{and}\qquad m_1^\theta:=\frac1N\sum_{i=1}^N\theta_i\sigma_i^1.
\end{equation*}
The identities in~\eqref{e.general_path_derivatives} for the derivatives of $\bar F_N$ with respect to $q_1$ and $q_2$ remain valid, and we also have
\begin{equation}\label{e.ext_x_derivatives}
    \partial_{x_1}\bar F_N=\E\la m_1^\theta\ra\qquad\text{and}\qquad\partial_{x_2}\bar F_N=-\E\la R_1^{yU'}\ra,
\end{equation}
while the time derivative now is
\begin{equation}\label{e.ext_time_derivative}
    \partial_t\bar F_N=\frac12\E\la R_{12}^{U'}R_{12}^{\sigma}\ra-\E\la R_1^{yU'}m_1^\theta\ra.
\end{equation}
The second term on the right-hand side of~\eqref{e.ext_time_derivative} comes from the derivative with respect to $t$ of the second term on the right-hand side of~\eqref{e.def.S_j.ext}. Writing $\partial_{(q_i,x_i)}:=(\partial_{q_i},\partial_{x_i})$, and denoting by $\la\cdot,\cdot\ra_{L^2\times\R}$ the scalar product on $L^2([0,1))\times\R$, as in Section~\ref{s.moving_scalar}, we obtain
\begin{multline}\label{e.ext_enriched_hj}
    \partial_t\bar F_N-\la\partial_{(q_1,x_1)}\bar F_N,\partial_{(q_2,x_2)}\bar F_N\ra_{L^2\times\R}\\
    =\frac12\E\la R_{12}^{U'}\Ll(R_{12}^{\sigma}-\E\la R_{12}^{\sigma}\mid R_{12}^{\alpha}\ra\Rr)\ra-\E\la R_1^{yU'}\Ll(m_1^\theta-\E\la m_1^\theta\ra\Rr)\ra.
\end{multline}
The first term on the right-hand side of~\eqref{e.ext_enriched_hj} is the same as in~\eqref{e.general_enriched_hj}. The second term is new, and it is small as soon as the weighted magnetization $m_1^\theta$ is concentrated around its mean. More precisely, setting $L:=\max_{y\in\{-1,1\}}\|U_y'\|_\infty$, we have $|R_{12}^{U'}|\leq\frac{M_N}NL^2$ and $|R_1^{yU'}|\leq\frac{M_N}NL$, and by the Cauchy--Schwarz inequality, the absolute value of the right-hand side of~\eqref{e.ext_enriched_hj} is at most $\frac{M_N}N(L^2+L)(\Delta_N^\ext)^{1/2}$, where
\begin{equation}\label{e.ext_concentration_error}
    \Delta_N^\ext:=\E\la\Ll(R_{12}^{\sigma}-\E\la R_{12}^{\sigma}\mid R_{12}^{\alpha}\ra\Rr)^2\ra+\E\la\Ll(m_1^\theta-\E\la m_1^\theta\ra\Rr)^2\ra.
\end{equation}
Whenever the quantity $\Delta_N^\ext$ is small, the function $\bar F_N$ is thus an approximate solution of the equation
\begin{equation}\label{e.ext_path_hj}
    \partial_tf-\la\partial_{(q_1,x_1)}f,\partial_{(q_2,x_2)}f\ra_{L^2\times\R}=0.
\end{equation}
This is the equation~\eqref{e.moving_scalar_bilinear} studied in Section~\ref{s.moving_scalar}: the nonlinearity is again bilinear, now with respect to the two groups of variables $(q_1,x_1)$ and $(q_2,x_2)$. By~\eqref{e.general_path_derivatives} and~\eqref{e.ext_x_derivatives}, the derivatives of $\bar F_N$ with respect to $q_1$, $x_1$, $q_2$ and $x_2$ are bounded by $1$, $\|\theta_1\|_\infty$, $\frac{M_N}{2N}L^2$ and $\frac{M_N}NL$ respectively. We therefore interpret the equation~\eqref{e.ext_path_hj} in the sense of Section~\ref{s.moving_scalar}, with the derivative bounds
\begin{equation}\label{e.ext_bounds}
    B:=\frac\lambda2L^2,\qquad R_1:=\|\theta_1\|_\infty,\qquad\text{and}\qquad R_2:=\lambda L.
\end{equation}

\smallskip
\noindent\textbf{The initial condition.}
The enriched free energy at $t=0$ is computed as in the previous subsections. We define $\mathcal A^\ext:K_1\times\R\to\R$ and $\mathcal B:\mcl Q_\infty\times\R\to\R$ by
\begin{align}
    \mathcal A^\ext(q_2,x_2)&:=-\E\log\int\E_{z_1}\exp\Ll(U_{y_1}\Ll(W_1^{q_2}(\alpha)+y_1x_2+\sqrt{1-q_2(1)}\,z_1\Rr)\Rr)\,\d\fR(\alpha),\label{e.def.A.ext}\\
    \mathcal B(q_1,b)&:=-\E\log\int e^{-q_1(1)}\cosh\Ll(\sqrt{2}\,V_1^{q_1}(\alpha)+b\Rr)\,\d\fR(\alpha).\label{e.def.B}
\end{align}
The function $\mathcal A^\ext(\cdot,x_2)$ is the average over the label $y_1$ of the pattern functionals associated with the shifted activation functions $U_{y_1}(\cdot+y_1x_2)$. The function $\mathcal B(\cdot,b)$ is, up to the additive constant $-\log\cosh b$, the cascade transform of the law of an Ising spin tilted by the external field $b$; in particular, we have $\mathcal B(q_1,0)=\psi_{\mu_{\mcl I}}(q_1)$. The additivity property of the cascade gives
\begin{equation}\label{e.ext_initial}
    \bar F_N(0,q,x)=\E\Ll[\mathcal B(q_1,-x_1\theta_1)\Rr]+\frac{M_N}{N}\mathcal A^\ext(q_2,x_2).
\end{equation}
As in the previous subsections, this is the sum of a function of the first group of variables $(q_1,x_1)$ and of a function of the second group of variables $(q_2,x_2)$. However, the first function is no longer convex: for instance, the function $b\mapsto\mathcal B(0,b)=-\log\cosh b$ is concave. Proposition~\ref{p.moving_scalar} therefore does not apply directly to this initial condition.

\smallskip
\noindent\textbf{Fixing the magnetization.}
To overcome this difficulty, we follow~\cite{boursier2026free} and sort the spin configurations according to their magnetization. We call magnetization profile any element $m$ of $L^\infty(\mu;[-1,1])$, that is, any measurable function from $\R$ to $[-1,1]$, two such functions being identified if they coincide $\mu$-almost everywhere. Given a magnetization profile $m$, we consider the model obtained by restricting the integral over $\sigma$ in~\eqref{e.ext_enrichment} to the configurations with magnetization profile $m$; informally, these are the configurations such that, for every $\vartheta$ in the support of $\mu$, the average of the spins $\sigma_i$ over the indices $i$ for which $\theta_i$ is close to $\vartheta$ is close to $m(\vartheta)$. For these configurations, the weighted magnetization $\frac1N\sum_{i=1}^N\theta_i\sigma_i$ is close to
\begin{equation}\label{e.ext_mbar}
    \bar m:=\E\Ll[\theta_1m(\theta_1)\Rr].
\end{equation}
The restriction to a magnetization profile has two effects. First, the derivation of the identity~\eqref{e.ext_enriched_hj} does not depend on the set of spin configurations over which we integrate, so this identity remains valid for the restricted model. Since the weighted magnetization is essentially constant on the set of configurations with magnetization profile $m$, the second term on the right-hand side of~\eqref{e.ext_enriched_hj} is then negligible, and we are left with the same concentration requirement as in the previous subsections. Second, the initial condition becomes convex with respect to $(q_1,x_1)$. To see this, for every $a\in[-1,1]$, we define $\psi_{\mu_{\mcl I},a}:\mcl Q_\infty\to\R$ by
\begin{equation}\label{e.def.psi.I.a}
    \psi_{\mu_{\mcl I},a}(r):=\sup_{b\in\R}\Ll\{\mathcal B(r,b)+ba\Rr\}.
\end{equation}
This function is the cascade transform of Ising spins whose magnetization is constrained to be equal to $a$; see~\cite[Proposition~7.5]{boursier2026free}. It is convex, and $\psi_{\mu_{\mcl I},0}=\psi_{\mu_{\mcl I}}$; see~\cite[Lemma~6.7 and Proposition~6.10]{boursier2026free}. For every magnetization profile $m$, we define $\Psi_m:\mcl Q_\infty\to\R$ and $\Psi_m^*:\mcl Q_\infty\to\R\cup\{+\infty\}$ by
\begin{equation}\label{e.def.Psi_m}
    \Psi_m(r):=\E\Ll[\psi_{\mu_{\mcl I},m(\theta_1)}(r)\Rr]\qquad\text{and}\qquad\Psi_m^*(p):=\sup_{r\in\mcl Q_\infty}\Ll\{\la p,r\ra_{L^2}-\Psi_m(r)\Rr\}.
\end{equation}
The enriched free energy at $t=0$ of the model restricted to the magnetization profile $m$ is then expected to converge to
\begin{equation}\label{e.ext_constrained_initial}
    \Psi_m(q_1)+x_1\bar m+\lambda\mathcal A^\ext(q_2,x_2).
\end{equation}
The function $(q_1,x_1)\mapsto\Psi_m(q_1)+x_1\bar m$ is convex, since $\Psi_m$ is an average of convex functions.

\smallskip
\noindent\textbf{The candidate limit at fixed magnetization.}
We now verify that Proposition~\ref{p.moving_scalar} applies to the initial condition in~\eqref{e.ext_constrained_initial}. We start with the part that depends on $(q_2,x_2)$. For every fixed $x_2\in\R$, Lemma~\ref{l.pattern_slope} applies to each of the two shifted activation functions $U_y(\cdot+yx_2)$, with a Lipschitz constant $L$ that depends neither on $x_2$ nor on $y$. Moreover, changing $x_2$ into $x_2'$ changes each of these two functions by at most $L|x_2'-x_2|$ in the supremum norm, and it thus changes $\mathcal A^\ext(q_2,\cdot)$ by at most the same amount. The function $\psi_2:=\lambda\mathcal A^\ext$ therefore satisfies the assumption~\eqref{e.moving_scalar_slope}, with $B$ and $R_2$ as in~\eqref{e.ext_bounds}. We now turn to the part that depends on $(q_1,x_1)$. For every $b\in\R$, the function $\mathcal B(\cdot,b)$ is a cascade transform up to an additive constant, so by Lemma~\ref{l.psi^vec_smooth}, it is $\mcl Q_2^*$-increasing, and Lipschitz continuous with respect to the $L^1$ norm, with Lipschitz constant $1$. These two properties are preserved by the supremum in~\eqref{e.def.psi.I.a} and by the average in~\eqref{e.def.Psi_m}. In particular, the function $\Psi_m$ extends to $\mcl Q_2$ by continuity. Since $\Psi_m$ is also convex, the arguments of Sections~\ref{s.cone} and~\ref{s.paths} leading to the identity~\eqref{e.gen.biconjugate}, which only rely on these three properties of $\psi_{\mu_{\mcl I}}$, also show that $\Psi_m^*$ is infinite outside of $K_1$, and that, for every $q_1\in\mcl Q_2$,
\begin{equation}\label{e.ext_biconjugate}
    \sup_{p\in\mcl Q_{\infty,\leq1}}\Ll\{\la p,q_1\ra_{L^2}-\Psi_m^*(p)\Rr\}=\Psi_m(q_1).
\end{equation}
We can thus apply Proposition~\ref{p.moving_scalar} with $\psi_2:=\lambda\mathcal A^\ext$, and with the function $g:\cH\times\R\to\R\cup\{+\infty\}$ defined by $g(p,c):=\Psi_m^*(p)$ if $p\in K_1$ and $c=\bar m$, and $g(p,c):=+\infty$ otherwise. This function is proper, lower semicontinuous and convex. Its effective domain is contained in $K_1\times[-R_1,R_1]$, since $|\bar m|\leq\|\theta_1\|_\infty=R_1$, and by~\eqref{e.ext_biconjugate}, its convex conjugate is the function $(q_1,x_1)\mapsto\Psi_m(q_1)+x_1\bar m$. We obtain that the function $f_m:\widehat{\mathcal D}\times\R^2\to\R$ defined by
\begin{equation}\label{e.ext_candidate_m}
    f_m(t,q,x):=x_1\bar m+\sup_{p\in\mcl Q_{\infty,\leq1}}\Ll\{\la p,q_1\ra_{L^2}+\lambda\mathcal A^\ext\Ll(q_2+tp,\,x_2+t\bar m\Rr)-\Psi_m^*(p)\Rr\}
\end{equation}
is the unique Lipschitz viscosity solution of~\eqref{e.ext_path_hj} with derivative bounds $(B,R_1,R_2)$ whose initial condition is the function in~\eqref{e.ext_constrained_initial}. This function is our candidate for the limit of the enriched free energy of the model restricted to the magnetization profile $m$. In the formula~\eqref{e.ext_candidate_m}, the dual variable associated with $x_1$ is frozen at the value $\bar m$, because the initial condition is affine in $x_1$ with slope $\bar m$. Its only effect is to shift the variable $x_2$ by $t\bar m$, in agreement with the definition of $S_j$ in~\eqref{e.def.S_j.ext}, in which the labels are multiplied by the sum of $x_2$ and of $t$ times the weighted magnetization.

\smallskip
\noindent\textbf{The candidate limit.}
It remains to relate the restricted models to the original one. The integral over $\sigma$ in~\eqref{e.ext_enrichment} is the sum of its restrictions to the different magnetization profiles, and on the exponential scale, such a sum is governed by its largest term. Taking into account the minus sign in the definition of the free energy, and in analogy with~\cite[(1.14)]{boursier2026free}, we thus expect the enriched free energy $\bar F_N$ to be asymptotically equal to the infimum over the magnetization profiles of the enriched free energies of the restricted models. This leads to the candidate
\begin{equation}\label{e.ext_candidate}
    \lim_{N\to\infty}\bar F_N(t,q,x)\stackrel{?}{=}f(t,q,x):=\inf_{m\in L^\infty(\mu;[-1,1])}f_m(t,q,x),
\end{equation}
and in particular, at $(t,q,x)=(1,0,0)$, to the following prediction for the limit of the free energy in~\eqref{e.ext_free_energy}:
\begin{equation}\label{e.ext_free_energy_candidate}
    \lim_{N\to\infty}\bar F_N(1,0,0)\stackrel{?}{=}\inf_{m\in L^\infty(\mu;[-1,1])}\ \sup_{p\in\mcl Q_{\infty,\leq1}}\Ll\{\lambda\mathcal A^\ext\Ll(p,\E[\theta_1m(\theta_1)]\Rr)-\Psi_m^*(p)\Rr\}.
\end{equation}
We can perform two consistency checks on~\eqref{e.ext_candidate}. First, at $t=0$, the Legendre duality between the external field and the magnetization~\cite[Lemma~6.7]{boursier2026free} gives
\begin{equation}\label{e.ext_min_over_m}
    \inf_{m\in L^\infty(\mu;[-1,1])}\Ll\{\Psi_m(q_1)+x_1\E[\theta_1m(\theta_1)]\Rr\}=\E\Ll[\mathcal B(q_1,-x_1\theta_1)\Rr],
\end{equation}
so that $f(0,\cdot)$ is indeed the limit of the right-hand side of~\eqref{e.ext_initial}. Second, when $\mu=\delta_0$, a magnetization profile reduces to a number $a\in[-1,1]$, and we have $\bar m=0$ and $\Psi_m=\psi_{\mu_{\mcl I},a}\geq\psi_{\mu_{\mcl I},0}=\psi_{\mu_{\mcl I}}$, so that $\Psi_m^*\leq\psi_{\mu_{\mcl I}}^*$. The infimum in~\eqref{e.ext_candidate} is therefore attained at $m=0$, and at $x_2=0$, we recover the candidate in~\eqref{e.random_activation_candidate}, with $\kappa=\mathrm I$ and with $\pi$ the law of $y_1$.

\subsection{Gaussian mixtures and spherical spins}
\label{s.ext.spherical}

We now consider the model of Section~\ref{s.ext} with spherical spins in place of Ising spins. For this model, the rotational invariance of the sphere allows us to reduce the model to the spherical perceptron with random activation functions of Section~\ref{s.random_activations}. The reduction consists in sorting the spin configurations according to their projection onto the vector $\theta$; this is the spherical counterpart of the sorting of the Ising configurations according to their magnetization profile performed in Section~\ref{s.ext}. Once this projection is fixed, the remaining coordinates of the spin configuration range in a sphere of one dimension less, and they only interact with the components of the patterns that are orthogonal to $\theta$, which are standard Gaussian vectors. The components of the patterns along $\theta$ then only enter the model as parameters of the activation functions, on the same footing as the labels. A similar reduction is used in the proof of~\cite[Theorem~1.6]{fuSohn2026perceptrons}, in the context of a teacher--student model at zero temperature.

We keep the labels $(y_j)_{j\in\N}$, the random variables $(\theta_i)_{i\in\N}$ of law $\mu$, and the vector $\theta=(\theta_1,\ldots,\theta_N)$ of Section~\ref{s.ext}, and we recall that the expectation $\E$ also averages over these random variables. Concerning the activation functions, we only assume here that $U_{-1}$ and $U_1$ are Lipschitz continuous, and we denote by $L$ the larger of their two Lipschitz constants. We write $g_j:=(g_{ij})_{i\leq N}$ for the Gaussian part of the $j$-th pattern, and we consider the model whose free energy is $-\frac1N\E\log Z_N$, where
\begin{equation}\label{e.ext_sph_model}
    Z_N:=\int\exp\Ll(\sum_{j=1}^{M_N}U_{y_j}\Ll(\frac{g_j\cdot\sigma}{\sqrt N}+\frac{y_j}{N}\,\theta\cdot\sigma\Rr)\Rr)\,\d\mu_{\mathbb{S}_N}(\sigma).
\end{equation}
This is the free energy in~\eqref{e.ext_free_energy} with the measure $\mu_{\mathbb{S}_N}$ in place of $\mu_{\mcl I}^{\otimes N}$.

\smallskip
\noindent\textbf{Rotational invariance.}
The law of the family $(g_j)_{j\leq M_N}$ and the measure $\mu_{\mathbb{S}_N}$ are both invariant under the orthogonal transformations of $\R^N$, and the Gaussian vectors $(g_j)_{j\leq M_N}$ are independent of $\theta$ and of the labels. The law of $Z_N$ is therefore not modified if we replace $\theta$ by the vector $|\theta|\,e_1$ in~\eqref{e.ext_sph_model}, where $e_1$ denotes the first vector of the canonical basis of $\R^N$. After this replacement, the second term in the argument of $U_{y_j}$ becomes $y_j\,\Theta_N\,\sigma_1/\sqrt N$, where $\Theta_N:=N^{-1/2}|\theta|$. By the law of large numbers, the random variable $\Theta_N$ converges almost surely to
\begin{equation}\label{e.def.Theta}
    \Theta:=\Ll(\E\Ll[\theta_1^2\Rr]\Rr)^{\frac12}
\end{equation}
as $N$ tends to infinity. We therefore set
\begin{equation}\label{e.ext_sph_hat}
    \widehat Z_N:=\int\exp\Ll(\sum_{j=1}^{M_N}U_{y_j}\Ll(\frac{g_j\cdot\sigma}{\sqrt N}+y_j\,\Theta\,\frac{\sigma_1}{\sqrt N}\Rr)\Rr)\,\d\mu_{\mathbb{S}_N}(\sigma).
\end{equation}
Since $|\sigma_1|\leq\sqrt N$ for every $\sigma\in\Sph_N$, and since the activation functions are Lipschitz continuous with Lipschitz constant $L$, we have
\begin{equation}\label{e.ext_sph_rotation}
    \Ll|\frac1N\E\log Z_N-\frac1N\E\log\widehat Z_N\Rr|\leq\frac{M_N}{N}\,L\,\E\Ll[|\Theta_N-\Theta|\Rr].
\end{equation}
Since $\Theta_N\leq\|\theta_1\|_\infty$, the dominated convergence theorem ensures that the right-hand side of~\eqref{e.ext_sph_rotation} tends to zero as $N$ tends to infinity. It thus suffices to study the quantity $\widehat Z_N$. In particular, the limit of the free energy can only depend on the law $\mu$ through its second moment, and the role played by the magnetization profile in Section~\ref{s.ext} will be played here by a single real parameter.

\smallskip
\noindent\textbf{The model at fixed projection.}
In the definition of $\widehat Z_N$, the vector $\theta$ has been replaced by a multiple of $e_1$, so the projection of a spin configuration onto the direction of the shifts is given by its first coordinate. We now sort the spin configurations according to the value of this coordinate. We assume that $N\geq3$, we set $n:=N-1$, and we write every $\sigma\in\Sph_N$ such that $|\sigma_1|<\sqrt N$ in the form
\begin{equation}\label{e.ext_sph_slice}
    \sigma=\Ll(\sqrt N\,a,\ \sqrt{\frac{N(1-a^2)}{n}}\,\tau\Rr),\qquad\text{with}\quad a\in(-1,1)\quad\text{and}\quad\tau\in\Sph_n.
\end{equation}
Under the measure $\mu_{\mathbb{S}_N}$, the random variables $a$ and $\tau$ are independent, the law of $\tau$ is $\mu_{\mathbb{S}_n}$, and the law of $a$ has the density
\begin{equation}\label{e.ext_sph_density}
    \beta_N(a):=\frac{\Gamma(N/2)}{\sqrt\pi\,\Gamma((N-1)/2)}\,(1-a^2)^{\frac{N-3}2}
\end{equation}
with respect to the Lebesgue measure on $(-1,1)$. We decompose the Gaussian vectors accordingly: for every $j$, we write $g_j=(g_{1j},\widetilde g_j)$, with $\widetilde g_j:=(g_{ij})_{2\leq i\leq N}\in\R^n$. For $\sigma$ as in~\eqref{e.ext_sph_slice}, the argument of $U_{y_j}$ in~\eqref{e.ext_sph_hat} can then be rewritten as
\begin{equation}\label{e.ext_sph_preactivation}
    \frac{g_j\cdot\sigma}{\sqrt N}+y_j\,\Theta\,\frac{\sigma_1}{\sqrt N}=\sqrt{1-a^2}\ \frac{\widetilde g_j\cdot\tau}{\sqrt n}+a\,g_{1j}+a\,\Theta\,y_j.
\end{equation}
For every $a\in(-1,1)$, $y\in\{-1,1\}$ and $\zeta\in\R$, we define the function $U^a_{y,\zeta}:\R\to\R$ by
\begin{equation}\label{e.ext_sph_random_channels}
    U^a_{y,\zeta}(v):=U_y\Ll(\sqrt{1-a^2}\,v+a\,\zeta+a\,\Theta\,y\Rr),
\end{equation}
and we set
\begin{equation}\label{e.ext_sph_slice_partition}
    Z_{N,a}:=\int\exp\Ll(\sum_{j=1}^{M_N}U^a_{y_j,g_{1j}}\Ll(\frac{\widetilde g_j\cdot\tau}{\sqrt n}\Rr)\Rr)\,\d\mu_{\mathbb{S}_n}(\tau).
\end{equation}
By~\eqref{e.ext_sph_preactivation} and the description of the measure $\mu_{\mathbb{S}_N}$ given above, we have
\begin{equation}\label{e.ext_sph_disintegration}
    \widehat Z_N=\int_{-1}^1Z_{N,a}\,\beta_N(a)\,\d a.
\end{equation}

For each fixed $a\in(-1,1)$, the quantity $Z_{N,a}$ is the partition function of a spherical perceptron with random activation functions, as considered in Section~\ref{s.random_activations}. This perceptron has $n$ spins and $M_N$ patterns; its patterns are the standard Gaussian vectors $(\widetilde g_j)_{j\leq M_N}$, and the activation function applied to the $j$-th pattern is $U^a_{y_j,g_{1j}}$. The mark of the $j$-th pattern is thus the pair $(y_j,g_{1j})$. These marks are i.i.d., and they are independent of the vectors $(\widetilde g_j)_{j\leq M_N}$. The fact that they are not real-valued, as was assumed in Section~\ref{s.random_activations}, is immaterial. We stress that the Gaussian random variable $g_{1j}$, which is a coordinate of the $j$-th pattern in the original model, has now become a parameter of the activation function, on the same footing as the label $y_j$. This perceptron falls within the setting considered at the end of Section~\ref{s.random_activations}, which is also that of~\cite[Theorem~1.1]{fuSohn2026perceptrons}. Indeed, the functions in~\eqref{e.ext_sph_random_channels} are Lipschitz continuous with Lipschitz constant at most $L$, and since
\begin{equation*}
    \Ll|U^a_{y,\zeta}(0)\Rr|\leq\max_{y'\in\{-1,1\}}|U_{y'}(0)|+L\Ll(|\zeta|+\Theta\Rr),
\end{equation*}
their values at the origin are integrable with respect to the law of the marks. Moreover, the ratio $M_N/n$ converges to $\lambda$.

It remains to identify the pattern functional of this perceptron. In the notation of Section~\ref{s.random_activations}, it is the function $r\mapsto\E[\mathcal A_{U^a_{y_1,\zeta}}(r)]$, where $\zeta$ is a standard Gaussian random variable independent of $y_1$, and where the expectation is over $y_1$ and $\zeta$. We claim that, for every $r\in K_1$,
\begin{equation}\label{e.ext_sph_channel_identity}
    \E\Ll[\mathcal A_{U^a_{y_1,\zeta}}(r)\Rr]=\mathcal A^\ext\Ll(a^2\one+(1-a^2)r,\,a\Theta\Rr),
\end{equation}
where $\mathcal A^\ext$ is the functional defined in~\eqref{e.def.A.ext}, and where $\one$ denotes the constant path equal to one. Indeed, by the definitions in~\eqref{e.def.A} and~\eqref{e.ext_sph_random_channels}, the left-hand side of~\eqref{e.ext_sph_channel_identity} is
\begin{equation*}
    -\E\log\int\E_{z_1}\exp\Ll(U_{y_1}\Ll(a\,\zeta+\sqrt{1-a^2}\,W_1^{r}(\alpha)+a\,\Theta\,y_1+\sqrt{(1-a^2)(1-r(1))}\,z_1\Rr)\Rr)\,\d\fR(\alpha),
\end{equation*}
where the expectation $\E$ now also averages over $y_1$ and $\zeta$. Conditionally on $\fR$, the process $\alpha\mapsto a\,\zeta+\sqrt{1-a^2}\,W_1^r(\alpha)$ is centered Gaussian with covariance $a^2+(1-a^2)r(\alpha\wedge\alpha')$. It thus has the same law as the process $W_1^{q_2}$ with $q_2:=a^2\one+(1-a^2)r$, and for this choice of $q_2$, we have $1-q_2(1)=(1-a^2)(1-r(1))$. Comparing with~\eqref{e.def.A.ext}, we obtain~\eqref{e.ext_sph_channel_identity}. The path $a^2\one+(1-a^2)r$ appearing in this identity has a simple interpretation in terms of overlaps: if $\sigma^1$ and $\sigma^2$ are two spin configurations of the form~\eqref{e.ext_sph_slice}, with the same parameter $a$ and with $\tau^1$ and $\tau^2$ in place of $\tau$, then
\begin{equation*}
    \frac{\sigma^1\cdot\sigma^2}{N}=a^2+(1-a^2)\,\frac{\tau^1\cdot\tau^2}{n}.
\end{equation*}

For every $a\in(-1,1)$, we set
\begin{equation}\label{e.ext_sph_def_F}
    \mathcal F(a):=\sup_{r\in K_1}\Ll\{\lambda\,\mathcal A^\ext\Ll(a^2\one+(1-a^2)r,\,a\Theta\Rr)-\mathsf S(r)\Rr\},
\end{equation}
where $\mathsf S$ is the functional defined in~\eqref{e.def.S.sph}. By~\eqref{e.ext_sph_channel_identity}, the quantity $\mathcal F(a)$ is the right-hand side of~\eqref{e.random_activation_free_energy} for the spherical perceptron with partition function $Z_{N,a}$. Since normalizing by $N$ instead of $n$ does not affect the limit, the discussion of Section~\ref{s.random_activations} thus leads to the prediction that, for every $a\in(-1,1)$,
\begin{equation}\label{e.ext_sph_slice_limit}
    \lim_{N\to\infty}-\frac1N\E\log Z_{N,a}\stackrel{?}{=}\mathcal F(a).
\end{equation}
As explained there, this prediction coincides with the formula of~\cite[Theorem~1.1]{fuSohn2026perceptrons}.

\smallskip
\noindent\textbf{The candidate limit.}
By~\eqref{e.ext_sph_disintegration}, the quantity $\widehat Z_N$ is the integral over $a\in(-1,1)$ of $Z_{N,a}\,\beta_N(a)$, and by~\eqref{e.ext_sph_density}, we have, for every $a\in(-1,1)$,
\begin{equation}\label{e.ext_sph_entropy}
    \lim_{N\to\infty}-\frac1N\log\beta_N(a)=-\frac12\log(1-a^2).
\end{equation}
The right-hand side of~\eqref{e.ext_sph_entropy} is the entropic cost of the constraint that the projection of the spin configuration be close to $a$. In view of~\eqref{e.ext_sph_slice_limit}, we expect that, on the exponential scale, the function $a\mapsto Z_{N,a}\,\beta_N(a)$ behaves like $\exp\Ll(-N\Ll(\mathcal F(a)-\frac12\log(1-a^2)\Rr)\Rr)$, and that its integral is governed by its largest value. As in Section~\ref{s.ext}, and taking into account the minus sign in the definition of the free energy, we thus expect the free energy associated with $\widehat Z_N$ to converge to the infimum over $a$ of the sum of this cost and of $\mathcal F(a)$. Recalling~\eqref{e.ext_sph_rotation}, we are led to the following prediction for the limit of the free energy of the model:
\begin{equation}\label{e.ext_sph_variational}
    \lim_{N\to\infty}-\frac1N\E\log Z_N\stackrel{?}{=}\inf_{a\in(-1,1)}\Ll\{-\frac12\log(1-a^2)+\sup_{r\in K_1}\Ll\{\lambda\,\mathcal A^\ext\Ll(a^2\one+(1-a^2)r,\,a\Theta\Rr)-\mathsf S(r)\Rr\}\Rr\}.
\end{equation}

The passage from~\eqref{e.ext_sph_slice_limit} to~\eqref{e.ext_sph_variational} can be justified as follows. Since the activation functions are Lipschitz continuous, the identity~\eqref{e.ext_sph_preactivation} and the Cauchy--Schwarz inequality yield that, for every $a,b\in(-1,1)$,
\begin{equation}\label{e.ext_sph_holder}
    \Ll|\frac1N\log Z_{N,a}-\frac1N\log Z_{N,b}\Rr|\leq C_N\sqrt{|a-b|},
\end{equation}
where
\begin{equation*}
    C_N:=\frac{\sqrt2\,L}{N}\Ll(\sqrt{M_N}\,\|\widetilde G\|_{\mathrm{op}}+\sum_{j=1}^{M_N}|g_{1j}|+M_N\Theta\Rr),
\end{equation*}
and where $\|\widetilde G\|_{\mathrm{op}}$ denotes the operator norm of the matrix with $M_N$ rows and $n$ columns whose rows are the vectors $(\widetilde g_j)_{j\leq M_N}$. Since $\E\|\widetilde G\|_{\mathrm{op}}\leq\sqrt{M_N}+\sqrt n$, see for instance~\cite[Theorem~7.3.1]{vershynin2018high}, the expectation of $C_N$ is bounded uniformly over $N$. As explained at the end of Section~\ref{s.random_activations}, the prediction~\eqref{e.ext_sph_slice_limit} is equivalent to the convergence of $-\frac1N\log Z_{N,a}$ to $\mathcal F(a)$ in $L^1$. If this convergence holds for every $a\in(-1,1)$, then the estimate~\eqref{e.ext_sph_holder} ensures that the function $\mathcal F$ is bounded and uniformly continuous, and that
\begin{equation*}
    \lim_{N\to\infty}\E\Ll[\,\sup_{a\in(-1,1)}\Ll|\frac1N\log Z_{N,a}+\mathcal F(a)\Rr|\,\Rr]=0.
\end{equation*}
By~\eqref{e.ext_sph_rotation},~\eqref{e.ext_sph_disintegration} and the Laplace method, this yields~\eqref{e.ext_sph_variational}.

\small
\bibliographystyle{plain}
\newcommand{\noop}[1]{} \def\cprime{$'$}

\end{document}